\documentclass{article}

\usepackage[utf8]{inputenc}
\usepackage[T1]{fontenc}
\usepackage{geometry}
\usepackage{graphicx}
\usepackage{amsmath}
\usepackage{amssymb}
\usepackage{amsthm}
\usepackage{cases}
\usepackage{hyperref}
\usepackage{abstract}
\usepackage{float}
\usepackage{titling}
\usepackage[dvipsnames]{xcolor}
\usepackage[ruled,vlined]{algorithm2e}
\SetCommentSty{textit}

\usepackage{tikz}
\newcommand*\circled[1]{\tikz[baseline=(char.base)]{
            \node[shape=circle,draw,inner sep=1pt] (char) {#1};}}

\usepackage{enumerate}
\usepackage[shortlabels]{enumitem}
\usepackage[
	safeinputenc,
	backend=biber,
	style=alphabetic,
	maxnames = 10,
	minnames = 1,
	isbn=false,
	url=false,
	eprint=false
  ]{biblatex}
\hypersetup{
    colorlinks=true,
    linkcolor=blue,
    filecolor=magenta,      
    urlcolor=cyan,
    citecolor=ForestGreen
}
\newtheorem{Theorem}{Theorem}[section]
\newtheorem{Proposition}{Proposition}[section]
\newtheorem{Lemma}{Lemma}[section]
\newtheorem{Corollary}{Corollary}[section]
\newtheorem{Assumption}{Assumption}

\theoremstyle{definition}

\newtheorem{Definition}{Definition}[section]
\newtheorem{Remark}{Remark}[section]
\newtheorem{Observation}{Observation}[section]

\def\ra{\rangle}
\def\la{\langle}

\begin{document}
\title{Elastoplasticity with softening as a state-dependent sweeping process: non-uniqueness of solutions and emergence of shear bands in lattices of springs}
\author{Ivan Gudoshnikov \thanks{Institute of Mathematics of the Czech Academy of Sciences, \v{Z}itn\'{a} 25, 110 00, Praha 1, Czech Republic, \url{gudoshnikov@math.cas.cz}}}
\thanksmarkseries{arabic}
\date{2025}

\renewcommand{\thefootnote}{\fnsymbol{footnote}} 
\footnotetext{\emph{Key words:} state-dependent sweeping process, quasi-variational inequalities, plasticity with softening, non-uniqueness of solutions, shear band, strain localization.}     
\renewcommand{\thefootnote}{\arabic{footnote}} 

\renewcommand{\thefootnote}{\fnsymbol{footnote}} 
\footnotetext{\emph{Mathematics Subject Classification (2020):} 74C05, 74N30, 47J20, 47J26.}     
\renewcommand{\thefootnote}{\arabic{footnote}} 

\maketitle
\begin{abstract}
Plasticity with softening and fracture mechanics lead to ill-posed mathematical problems due
to the loss of monotonicity. Multiple co-existing solutions are possible when softening elements
are coupled together, and solutions cannot be continued beyond the point of complete degradation of the set of admissible stresses.  We present a state-dependent sweeping process which solves the evolution of elasto-plastic Lattice
Spring Models with arbitrary placement of softening, hardening and perfectly plastic springs.
Using numerical simulations of regular grid lattices with softening we demonstrate the emergence of non-symmetric shear bands with strain localization. At the same time, in toy examples
it is easy to analytically derive multiple co-existing solutions. These solutions correspond to fixed
points in the implicit catch-up algorithm and we observe a discontinuous bifurcation with the exchange
of stability of those fixed points.
\end{abstract}
\tableofcontents
\section{Introduction}
Modeling of plastic deformation process is an important scientific pursuit due to numerous practical applications e.g. in manufacturing, 
 construction and transportation industries. Different materials and phenomena are described by different models and their mathematical complexity varies significantly. 

We focus on rate-independent elastoplasticity models, which can be divided into {\it strain hardening}, {\it perfect plasticity} and {\it strain softening} types, classified by observed stress response to imposed deformation. In this text we develop a finite-dimensional in space and continuous in time 
model, which would cover all these types in a uniform manner. In such a model a transition between the qualitatively different types is due to a quantitative change of parameters, which paves the way to a mathematical study of this bifurcation. Now let us provide an introduction on the three types of plasticity to explain their qualitative differences, as well as the mathematical results and challenges associated with each type. 

At first we illustrate plastic behaviors with a simple example of a one-dimensional elasto-plastic spring. Let the spring be linearly elastic with stiffness $k>0$ near its stress-free state, and impose a prescribed monotonically increasing elongation $l(t)$ on the spring (Fig. \ref{fig:intro_fig_models} a). This and the rest of the models in the paper are quasi-static.

\subsection{Perfect plasticity}

   The simplest to construct is the model of a {\it perfectly plastic} material, with the {\it yield limit} $c_0$ being the only material parameter apart from $k$. At all times the stress $\sigma$ must stay in the {\it set of admissible stresses}: 
\begin{equation}
\sigma\in[-c_0, c_0]
\label{eq:admissible-stresses-pp}
\end{equation}
(we limit the scope of the paper to {\it isotropic} plasticity, which means that the sets of admissible stresses are symmetric with respect to $0$).
After the stress reaches the yield limit, all further elongation contributes to progressive plastic deformation, but the stress remains constant  (Fig. \ref{fig:intro_fig_plasticity} a).

\begin{figure}[h]\center
\includegraphics{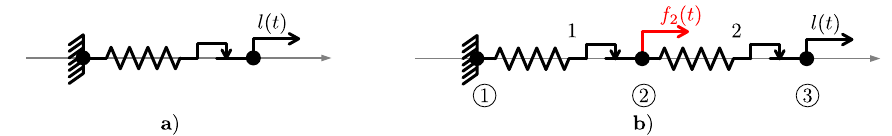}
\caption{
\footnotesize {\bf a)}  An elasto-plastic spring with the displacements of its left and right endpoints prescribed to be $0$ and $l(t)$ respectively. {\bf b)}A model of two elasto-plastic springs connected in series with displacements $0$ and $l(t)$ prescribed for the nodes $\protect\circled{1}$ and $\protect\circled{3}$ respectively, and a prescribed force $f_2(t)$, applied to the node $\protect\circled{2}$.
} 
\label{fig:intro_fig_models}
\end{figure}

\begin{figure}[h]\center
\includegraphics{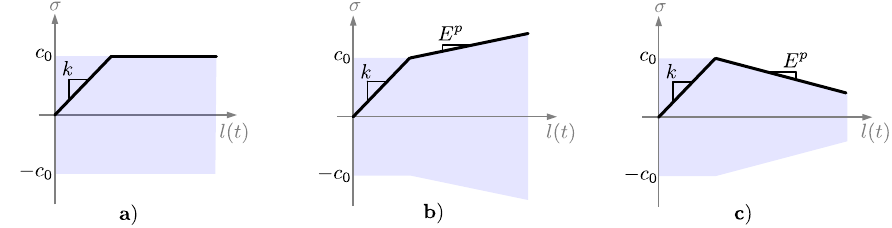}
\caption{
\footnotesize {\bf a)} Elasticity-perfect plasticity: the stress remains constant during plastic deformation, and the set of admissible stresses (the size of the blue area along the $\sigma$ axis) is fixed. {\bf b)}  Elasto-plasticity with hardening: during plastic deformation stress increases under imposed elongation, and  the set of admissible stresses enlarges. {\bf c)} Elasto-plasticity with softening: during plastic deformation stress decreases under imposed elongation, and  the set of admissible stresses contracts.} 
\label{fig:intro_fig_plasticity}
\end{figure}

Such simplicity is, however, deceptive, because the elongation-stress relation (or, in the case of continuous media, the strain-stress relation) is non-injective, which causes plenty of complications. If we connect just two elastic-perfectly plastic springs in a series to make a slightly more complex model (Fig. \ref{fig:intro_fig_models} b), then,  in general, the position of the connection point ($\protect\circled{2}$ in Fig. \ref{fig:intro_fig_models} b) cannot be determined uniquely. 

Even worse, continuum models, which are constructed as elliptic PDEs coupled with the laws of perfect plasticity, generally do not possess a solution for the strain as a function of the space variable even when the applied body force is smoothly distributed. Instead, as it was realized first by J.-J. Moreau \cite{Moreau1976}, strain variable should be a time-dependent measure. This is because in a generic scenario the displacement variable inevitably develops discontinuities on the surfaces of codimension $1$ inside a perfectly plastic body, and such discontinuities are incompatible with the trace theorem for Sobolev functions, see the discussions at \cite[p.~vi]{Temam1985}, \cite[p.~239]{DalMaso2006}. The discontinuity surfaces are called {\it strain localizations}, {\it shear bands} and {\it plastic slips} by various authors. 

While Moreau considered a prototypical example with one spatial dimension (i.e. a rod), in order to describe strain localizations in two- and three-dimensional domains the theory of the space of bounded deformations ($BD$) was developed with major  contributions by P.-M. Suquet \cite{Suquet1978a}, \cite{Suquet1980}, \cite{Suquet1981}, \cite{Suquet1988}, R. Temam and G. Strang \cite{Strang1980}, \cite{TemamStrangBD}, \cite{Temam1985}. More recent treatments of the problems with perfect plasticity include \cite{Ebobisse2004}, \cite{DalMaso2006}, \cite{Demyanov2009}, \cite{Francfort2016}, \cite{Mora2016}. Additionally we suggest our previous work \cite{Gudoshnikov2025rl} focused on the detailed comparison of discrete (such as Fig. \ref{fig:intro_fig_models} b) and continuum (PDE) models with perfect plasticity. It must also be noted that it is possible to workaround the non-uniqueness and discontinuity issues in a model of a thin beam by considering only small displacements, orthogonal to the beam \cite{Krejci2007}.

\subsection{Plasticity with hardening}
Thus  well-posedness requires injective stress-strain relation, and it is such for another class of rate-independent models: plasticity with {\it hardening}. Hardening here means that, upon reaching the yield limit, the stress continues to increase at a slower rate $E^p<k$ along with progressive plastic deformation (see Fig. \ref{fig:intro_fig_plasticity} b), and such rate $E^p$ (called the {\it modulus of plasticity}) is positive --- unlike the case of perfect plasticity. It is common to think of hardening as of an expansion of the set of admissible stresses: in the one-dimensional spring example, instead of constraint \eqref{eq:admissible-stresses-pp} 
we have, say,
\begin{equation}
\sigma \in [-(c_0+s\alpha), c_0+s\alpha],
\label{eq:admissible-stresses-hardening}
\end{equation}
where $s$ is another material constant to characterize plasticity, and $\alpha$ is a time-dependent unknown variable, usually called the {\it internal variable}. In the models with hardening the value of $\alpha$ changes during the plastic deformation, and so does the set of admissible stresses, see again Fig. \ref{fig:intro_fig_plasticity} b) and also see \cite[Sect. 1.2.2, Figs. 1.4B, 1.6]{SimoHughes1998} for more details. 



Plasticity with hardening is the ``nicest'' among the problems of our interest in the sense that the evolution of all unknown variables can be solved uniquely, and, for continuous media in particular, with values in Sobolev spaces. The well-posedness of the problem with hardening attracted many researchers in mathematics and engineering, and the literature on the topic is abundant. For proper mathematical formulations and the proofs of existence and uniqueness of a solution for continuous elastoplastic media with hardening we refer to \cite[Sect. 13.5]{Necas1981} and to the monumental book \cite{Han2012}, which is solely focused on solving this particular type of the problem. The corresponding numerical algorithms are also presented in \cite{SimoHughes1998}. In turn, spatially discrete models (networks of elasto-plastic springs) with hardening are the main topic of the paper \cite{Martins2007truss}.

\subsection{Moreau's sweeping process and evolution variational inequalities}
\label{ssec:sp_intro}
From the mathematical perspective such mechanical problems can be posed as {\it evolution variational inequalities}, which were properly studied since 1970s \cite{Duvant1976_eng} at least. In particular, the problem of plasticity with hardening and the stress problem of perfect plasticity can be expressed mathematically as a variational inequality known as a {\it sweeping process}, devised by Moreau about the same time \cite[Sect. 5.f]{Moreau1973}, \cite{Moreau1977}.

The sweeping process is an abstract problem to find a trajectory of a point $X$, which is subject to one-sided constraint
\begin{equation}
X(t)\in \mathcal{C}(t),
\label{eq:sp_constr}
\end{equation}
where $\mathcal{C}(t)$ is a given time-dependent closed convex set, usually called the {\it moving set}. It can be visualized as a container, placed on top of a small object laying on a rough surface, so that the moving boundary of the container pushes the object upon contact, see \cite[Fig. 1.2, p. 4]{Mielke2015}. Naturally, \eqref{eq:sp_constr} alone is not enough to specify a trajectory of $X$, and the sweeping process is written mathematically as the following inclusion:
\begin{equation}
-\dot X(t)\in N_{\mathcal{C}(t)}(X(t)) \qquad \text{for a.a. }t\in [0,T],
\label{eq:sp_classic}
\end{equation}
where the right-hand side is the cone of outward normals (supporting vectors) to the convex set $\mathcal{C}(t)$ at $X(t)$, defined only for \eqref{eq:sp_constr} as
\begin{equation}
N_{\mathcal{C}(t)}(X(t)) = \left\{Y: \la Y, c-X(t)\ra\leqslant 0 \text{ for all }c\in \mathcal{C}(t)\right\}.
\label{eq:nc_specific_intro}
\end{equation}
Substitute \eqref{eq:nc_specific_intro} into \eqref{eq:sp_classic} to obtain the evolution variational inequality form of the problem:
\begin{equation}
\text{find } X(t)\in \mathcal{C}(t) \quad \text{ such that }\quad \la\dot X(t), c-X(t)\ra\geqslant 0\quad  \text{ for all }\quad  c\in \mathcal{C}(t),\quad\text{ a.a. }\quad t\in[0,T].
\label{eq:intro_evi}
\end{equation}

For a reader-friendly introduction to the theory we refer to \cite[pp. 1--16]{Kunze2000}, from which we will state the existence and unqueness:
\begin{Theorem}
\label{th:sp_classic} Let $\mathcal{H}$ be a Hilbert space, and $\mathcal{C}(t)\subset \mathcal{H}$ be a closed convex nonempty set for all $t\in[0,T]$. Let $\mathcal{C}(t)$ be Lipschitz-continuous with respect to the Hausdorff distance, i.e. there is $L>0$ such that
\[
d_{H} (\mathcal{C}(t_1)-\mathcal{C}(t_2))\leqslant L|t_1-t_2|\qquad \text{for all }t_1, t_2\in[0,T].
\]
Then for every initial value $X(0)\in \mathcal{C}(0)$ there exists a unique solution of the sweeping process \eqref{eq:sp_classic}, which will be Lipschitz-continuous with the same constant $L$. 
\end{Theorem}

When the problems of elastoplasticity are expressed as a sweeping process, Theorem \ref{th:sp_classic} yields existence and uniqueness
\begin{itemize}
\item of the stress solution in problems with perfect plasticity: an interested reader may compare \cite[(4.20)--(4.22), p.~248]{Duvant1976_eng}, \cite[(4.16), p.~326]{Necas1981} to \eqref{eq:intro_evi}, \cite[(18), p.~8]{Gudoshnikov2021ESAIM}, \cite[(48), p.~23]{Gudoshnikov2023preprint} to \eqref{eq:sp_classic}, also see \cite[Sect. 4]{Gudoshnikov2025rl}.
\item of the solution in problems with hardening: compare \cite[(5.11), p.~331]{Necas1981},\cite[(8.16), p.~230]{Han2012} to \eqref{eq:intro_evi}.
\end{itemize}
\subsection{Plasticity with softening}
\label{ssect:introduction_softening}
Finally, let us come back to the example of an elasto-plastic spring and consider the case when the stress reduces along with progressive plastic deformation. This phenomenon is called {\it softening}, and it can be modeled as a contraction of the set of admissible stresses, see Fig. \ref{fig:intro_fig_plasticity} c). Softening is the most difficult case in our list, as the elongation-stress (or strain-stress) relation is ``worse'' than non-injective --- it is non-monotone. 

If just two springs with softening are coupled in series (Fig. \ref{fig:intro_fig_models} b) then, similarly to perfect plasticity, we have multiple solutions for the position of the connection point  \textcircled{\raisebox{-0.9pt}{2}}, but  now the stress values are non-unique as well (we will verify it in Section \ref{sect:ex1}, but also it follows from e.g. \cite{ChenBaker2004}). Written in a variational form, such a system has a non-convex energy increment with multiple isolated minima and a saddle point. 

Predictably, continuous media with the strain-stress relation of the type of Fig. \ref{fig:intro_fig_plasticity} c) at each point of a body make an extremely challenging problem. Consider even the simplified equilibrium problem without inequalities:  for a small time-step during plastic deformation a known subset of points must follow the linear yet decreasing strain-stress relation. Such a problem {\it will not be elliptic}, thus no opportunity to use Lax-Milgram-type theorems (e.g. \cite[Section~1.1]{Gatica2014}). 

The non-convexity in the variational formulation also rules out the standard infinite-dimensional convex optimization (e.g. \cite[Section~II.1, pp.~34--36]{EkelandTemam1999}).  Analytic results on a one-dimensional continuum in the variational formulation can be found in \cite{Bradies1999}, \cite{Lambrecht2003}. For $2$- and $3$-dimensional bodies there exists a rather difficult variational theory by Dal Maso et al. \cite{DalMaso2008Softening1}, \cite{DalMaso2008Softening2}. They have demonstrated the existence of a solution in the space of compatible systems of generalized Young measures, confirmed non-uniqueness of solutions \cite[Remark 5.2]{DalMaso2008Softening2} and the presence of strain localizations \cite[Section 8.2]{DalMaso2008Softening1} in their formulation.

On the other hand, plasticity with softening was researched extensively in mechanical and civil engineering, where it was observed that strain localization due to softening is the initial phase of {\it ductile fracture} \cite{Bazant1986} \cite{Bazant1988} \cite{Gruben2017}.   In other words, ductile fracture follows plastic deformation, and this makes it different from brittle fracture, which may occur abruptly without  any preceding plastic deformation \cite[Chapter 13]{Campbell2008Metallurgy}. Therefore, in order to predict the failure of a ductile material, one must have a model for elasto-plasticity with softening. 


We will not attempt to provide a comprehensive overview of the great amount of engineering publications on the topic, and only mention a few of ample contributions by Z. P. Ba{\v{z}}ant \cite{Bazant1986} \cite{Bazant1988}, M. Jir{\'{a}}sek \cite{Jirasek2003}, \cite{Jirasek2007} and R. de Borst \cite{deBorst1991}, \cite{deBorst1993}, \cite{deBorst2019}, where the reader can find related discussions. Notably, there is a fundamental issue of severe mesh-dependence in numerical simulations due to the lack of ellipticity of the problem. Finite element approximations develop strain localizations to a layer of single elements and such dependence on the size of elements invalidates objectivity of the models. For this reason much of the literature is devoted to non-local modifications of the constitutive law, which would limit strain localizations.

\begin{figure}[H]\center
\includegraphics{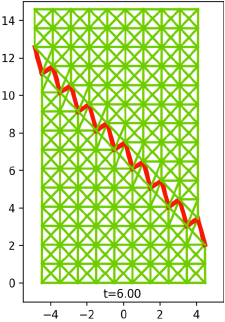}
\caption{\footnotesize A network of springs with softening developed a shear band. Such numerical examples will be presented in Section \ref{sect:LSM-examples}.
} \label{fig:paper1_big_lattice1_softening_solution2}
\end{figure}

Because of such difficulties, we will not consider continuum media in the current text, and focus instead of the discrete models, such as toy models of Fig.~\ref{fig:intro_fig_models} and networks of springs (also known as Lattice Spring Models (LSMs), or, informally, ``trusses''), see Fig.~\ref{fig:paper1_big_lattice1_softening_solution2}. Such spring models have already been used to model fracture \cite{Chen2014}, \cite{Meng2022} and our work is a step towards a rigorous mathematical theory for them.
In this text we extend the sweeping process approach of \cite{Gudoshnikov2023preprint} and formulate the elastoplasticity problem for discrete media in terms of stress and internal variable (``the dual problem'' according to \cite{Han2012})  is such a way that includes not only perfect plasticity, but also hardening and, most importantly, softening. However, this will inevitably require a more general formalism than the ``classical'' sweeping process \eqref{eq:sp_classic}. 

\subsection{State-dependent sweeping process}
As we have mentioned in the beginning of the previous Section \ref{ssect:introduction_softening}, in the case of softening the evolution of stresses may not be unique even in such a simple model as Fig. \ref{fig:intro_fig_models} b).      
This is clearly incompatible with uniqueness of the solution, stated in Theorem \ref{th:sp_classic}, and excludes any expectation that the problem with softening can be converted to the form \eqref{eq:sp_classic}.

Instead, we will adopt a richer framework of the {\it state-dependent sweeping process}, i.e. the modification of \eqref{eq:sp_classic} with the moving set $\mathcal{C}$ not just being a prescribed function of $t$, but also depending on the state variable $X$:
\begin{equation}
-\dot X \in N_{\mathcal{C}(t,X)}(X).
\label{eq:abstract_sp}
\end{equation}

This type of problem was investigated by Kunze and Monteiro Marques in \cite{Kunze1998}, see also \cite[Section~3.3, pp.~23--29]{Kunze2000}, from where we adopt the counterpart of Theorem \ref{th:sp_classic} for the state-dependent sweeping process in finite dimensions:

\begin{Theorem}\label{th:sp_state-dep}
{\bf \cite[Th.~3.3, p.~184]{Kunze1998}, \cite[Th.~6, p.~25]{Kunze2000}} Let $E$ be a Euclidean space, and $\mathcal{C}(t,X)\subset E$ be a closed convex nonempty set for all $t\in[0,T], X\in E$. Let $\mathcal{C}(t,X)$ have Lipschitz-continuous dependence on $t$ and $X$ with respect to the Hausdorff distance, i.e. there are $L_1, L_2>0$ such that
\begin{equation}
d_{H} (\mathcal{C}(t_1,X_1)-\mathcal{C}(t_2,X_2))\leqslant L_1|t_1-t_2|+L_2\|X_1-X_2\|\qquad \text{for all }t_1, t_2\in[0,T], \,X_1,X_2\in E.
\label{eq:abstract-state-dep-set-Lip-continuous}
\end{equation}
If 
\begin{equation}
L_2<1
\label{eq:L2-less-than-one}
\end{equation}
then for every initial value $X(0)$, compatible in the sense that $X(0)\in \mathcal{C}(0, X(0))$, there exists at least one solution of the state-dependent sweeping process \eqref{eq:abstract_sp}. 
\end{Theorem}

Conditions \eqref{eq:abstract-state-dep-set-Lip-continuous}--\eqref{eq:L2-less-than-one} are nowadays the standard requirement for the existence of a solution in this type of problems, generalized to many further modifications of the problem \eqref{eq:abstract_sp} (perturbed by a vector field, with a nonconvex prox-regular $\mathcal{C}$, with a lower regularity in time etc.). But condition \eqref{eq:L2-less-than-one} remains essential for existence, as for $L_2\geqslant 1$ there are simple counterexamples with no solutions (see the references above). 

Moreover, uniqueness of a solution to \eqref{eq:abstract_sp} is not guaranteed at all regardless of the value of $L_2$, see the discussion \cite[p.~112]{BrokateKrejciSchnabel2004}. It is also shown there that uniqueness can be achieved if set $\mathcal{C}$ has a sufficiently smooth boundary. But in the problems of our interest neither uniqueness of a solution, nor a smooth boundary of $\mathcal{C}$ will be present.

As the central result of this text, the problem of quasi-static evolution of a network of springs (including softening, hardening, and perfectly plastic springs) will be converted to the initial value problem, governed by a state-dependent sweeping process \eqref{eq:abstract_sp} (Theorem \ref{th:LSMtoMSP}) and solved numerically. This way, all of the plastic properties of the springs and the way they are connected will be encoded into the set $\mathcal{C}(t,X)$.

To find a numerical solution to the problem \eqref{eq:abstract_sp} we use the {\it implicit catch-up algorithm} of \cite{Kunze1998}. In the algorithm we must take a partition of the time-domain $[0,T]$:
\[
0=t_0<\dots<t_{i-1}<t_{i}<\dots t_k=T,\quad i\in \overline{1,k}
\]
and for each time-step we must find $X_{i+1}\in E$, such that
\begin{equation}
X_{i}- X_{i+1} \in N_{\mathcal{C}(t_{i+1}, X_{i+1})}(X_{i+1}),
\label{eq:discrete-process-abstract}
\end{equation}
which is equivalent to the fixed-point problem, involving a projection onto a closed convex set:
\begin{equation}
X_{i+1} = {\rm proj}(X_{i}, \mathcal{C}(t_{i+1}, X_{i+1})).
\label{eq:fixed-point-abstract}
\end{equation}
To find such fixed points we will numerically run the iterations of the map
\begin{equation}
\widetilde{X}\mapsto  {\rm proj}\left(X_{i}, \mathcal{C}\left(t_{i+1}, \widetilde{X}\right)\right).
\label{eq:map-to-iterate-abstract}
\end{equation}
This method works in practice, yet it raises the questions not only about the {\it existence} of the fixed points, but also about their {\it stability} as equilibria of the discrete dynamical system, defined by the map \eqref{eq:map-to-iterate-abstract}. And, although the conditions of Theorem \ref{th:sp_state-dep} yield the existence (see \cite[Lemma~2.1, p.~182]{Kunze1998}), we will see that they only cover the ``nicest'' case of hardening. Furthermore, stability of such dynamical systems 
is largely unexplored. 

Another way to look at the problem \eqref{eq:discrete-process-abstract} is to convert it to the equivalent form of a {\it quasi-variational inequality} (see \cite[Sect.~11, pp. 237--261]{Baiocchi1984} and references therein), for which there are more numerical methods available, see e.g. \cite{Outrata2013}.

Our particular line of research on modeling softening as a state-dependent sweeping process was initiated in 1990s by M. Brokate and H. Studt, but it was interrupted at the time. One can still find the trace of Studt's PhD thesis (see reference {[10]} in \cite{Kunze1998}), which is confirmed to be never finished and never published. The work was continued by M. Brokate, P. Krej{\v{c}}{\'{i}} and H. Schnabel, with their results in \cite{BrokateKrejciSchnabel2004}, \cite{Schnabel2006} aimed at a particular type of softening, known as Gurson's model \cite{Gurson1977}, \cite{Needleman1984}, \cite{Tvergaard1989}. Brokate et al. proved uniqueness of solutions, yet they consider single-element (``zero-dimensional'') models . In the current work we demonstrate how coupling of just two elements (Fig.~\ref{fig:intro_fig_models} b) with softening leads to non-uniqueness.

The paper is organized as follows:
\begin{itemize}
\item
in the next Section \ref{sect:notation} we briefly remind the reader a few basic definitions and properties, which we will use later.
\item In Section \ref{sect:single_spring} we derive the state-dependent sweeping process for a single spring and show that with different values of parameters we can reproduce all three types of plasticity of Fig. \ref{fig:intro_fig_plasticity}. Such state-dependent sweeping process is two-dimensional, as it represents the evolution of stress and internal variable in the toy model of Fig. \ref{fig:intro_fig_models} a). 
\item In Section \ref{sect:LSMs} we connect $m$ such springs into a lattice according to a given graph and derive the state-dependent sweeping process in $2m$ dimensions to model the evolution of the entire lattice. The dimension $2m$ can be reduced, depending on the underlying graph of the lattice, and this is important for speeding up the numerical algorithm, which we also discuss in that section.
\item In Section \ref{sect:ex1} we focus on a specific example of the lattice, the toy model of Fig. \ref{fig:intro_fig_models} b). We analytically show non-uniqueness of solutions in the corresponding state-dependent sweeping process, which agrees with \cite{ChenBaker2004}, where the same spring model is studied. Noteworthy, for this particular state-dependent sweeping process we can directly observe how the implicit catch-up algotithm can find unique and non-unique solutions. Also we visualize the moving set $\mathcal{C}(t,X)$ to understand how such multiple solution are possible in the case of softening.
\item In Section \ref{sect:LSM-examples} we provide several examples of more complex Lattice Spring Models (similar to Fig.~\ref{fig:paper1_big_lattice1_softening_solution2}) to compare the cases of hardening, perfect plasticity and softening. We observe the formation of shear bands and non-uniqueness of solutions, especially in the latter case. 
\item Finally, Section \ref{sect:conclusions} contains a concluding summary and a discussion on open questions. 
\end{itemize}

\section{Notation and technical preliminaries}
\label{sect:notation}
By $0_n$ we denote the zero vector in $\mathbb{R}^n$ and $I_{n \times n}$ is the $n \times n$ identity matrix. We use the following notation for the normal cone in a subspace of $\mathbb{R}^n$ (see \cite[p. 15]{Rockafellar1970}, \cite[p. 3]{Kunze2000}, \cite[p. 125]{Bauschke2011} for a general definition of the normal cone $N_C(x)$ a Hilbert space):
\begin{Definition} Given a vector space $E\subset \mathbb{R}^n$, for an $n\times n$ symmetric positive definite matrix $S$ we can define an inner product on $E$
\begin{equation}
(x,y) \mapsto \la x,y\ra_S:= x^{\top} S y \qquad \text{for any }x,y \in E.
\label{eq:abstract_ip_via_a_matrix}
\end{equation}
Now, for a closed convex nonempty set $\mathcal{C}\subset E$ and a point $x\in E$ we define the {\it outward normal cone to $\mathcal{C}$ at $x$ in sense of inner product $\la\cdot, \cdot\ra_S$} as 
\[
N^S_{\mathcal{C}}(x):=\{y\in E: y^{\top}S(c-x)\leqslant 0 \text{ for all }c\in \mathcal{C}\},
\]
If $S=I$ we simply write $N_{\mathcal{C}}(x)$.
\label{def:normal_cone_in_E}
\end{Definition}
\begin{Proposition} In the setting of Definition \ref{def:normal_cone_in_E} it follows from the definition that for any $y\in E$
\begin{equation}
N^S_{\mathcal{C}}(x+y)=N^S_{\mathcal{C}-y}(x),
\label{eq:nc_basic_property1}
\end{equation}
\begin{equation}
SN_{\mathcal{C}}(x) = N^{S^{-1}}_{\mathcal{C}}(x).
\label{eq:nc_basic_property2}
\end{equation}
\begin{equation}
y-x\in N_{\mathcal{C}}^S(x) \qquad \Longleftrightarrow \qquad x = {\rm proj}^S(y, \mathcal{C})
\label{eq:proj_vs_normal_cone}
\end{equation}
where 
\[
{\rm proj}^S(y, \mathcal{C}):=\underset{c\in \mathcal{C}}{\rm arg\, min}\,(y-c)^{\top}S(y-c).
\]
\label{prop:nc_basic_property}
\end{Proposition}
The following fact may appear trivial, however, we verify it for the sake of completeness of the proofs.
\begin{Lemma} 
\label{lemma:apa-normal-cones}
Let $E$ be an Euclidean space with the inner product $\langle \cdot, \cdot\rangle_E$. Let us be given a family of $m$ convex polytopes $\mathcal{C}_i, i\in \mathbb{N}$.  Consider the space $E^m$ with the inner product 
\[
(x_i)_{i\in\overline{1,m}},\,(y_i)_{i\in\overline{1,m}}\mapsto \left\la (x_i)_{i\in\overline{1,m}},(y_i)_{i\in\overline{1,m}} \right\ra_{E^m}=\sum_{i=1}^m\la x_i,y_i\ra_E
\]
and the set
\[
\mathcal{C}:= \left\{(x_i)_{i\in\overline{1,m}}\in E^m: x_i\in \mathcal{C}_i \text{ for all }i\in \overline{1,m}\right\}= \mathcal{C}_1\times \mathcal{C}_{2}\times\dots \times \mathcal{C}_{m}.
\]
Then for any $x^*=(x^*_i)_{i\in \overline{1,m}}\in \mathcal{C}$
\begin{equation}
N_{\mathcal{C}}(x^*) = N_{\mathcal{C}_1}(x^*_1)\times N_{\mathcal{C}_2}(x^*_2) \times \dots \times N_{\mathcal{C}_m}(x^*_m).
\label{eq:apa-normal-cones}
\end{equation}
\end{Lemma}
\noindent{\bf Proof.} For each $i\in \overline{1,m}$ let 
\[
{\rm C}_i:= \left\{(x_j)_{j\in \overline{1,m}}\in E^m: x_i\in \mathcal{C}_i\right\} = E\times\dots \times E\times\mathcal{C}_{i}\times E \times \dots \times E.
\]
From the definition of the normal cone we can derive that
\begin{equation}
N_{{\rm C}_i}(x^*)=\{0\} \times \dots \times \{0\}\times N_{\mathcal{C}_i}(x^*_i)\times \{0\} \times \dots \times \{0\}.
\label{eq:apa-normal-cones2}
\end{equation}
Observe that $\mathcal{C}=\bigcap_{i=1}^m {\rm C}_i$, therefore by the additivity of the normal cones for polyhedra
\[
N_{\mathcal{C}}(x^*)=\sum_{i=1}^m N_{{\rm C}_i}(x^*),
\]
from which \eqref{eq:apa-normal-cones} follows because of \eqref{eq:apa-normal-cones2}. $\blacksquare$

\section{Single spring element: constitutive equations, the corresponding sweeping process and its analysis}
\label{sect:single_spring}
\subsection{An elasto-plastic spring element with the state-dependent normality rule}
\label{ssect:single-spring-model}
\begin{figure}[h]\center
\includegraphics{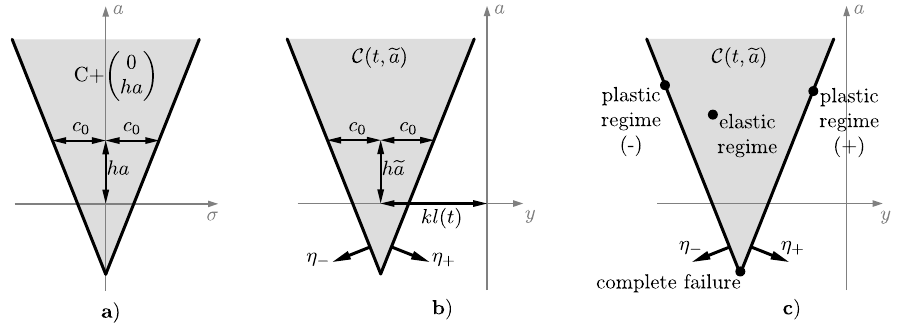}
\caption{
\footnotesize {\bf a)} Set ${\rm C}+\begin{pmatrix}0\\ h a\end{pmatrix}$ in plastic flow rule \eqref{eq:plastic_flow_single},\eqref{eq:tilde_C_single}. {\bf b)} The corresponding moving set \eqref{eq:tilde_C_single2} in the implicit sweeping process \eqref{eq:sp_single2}.
{\bf c)} Four possible regimes of the single spring model.} \label{fig:lin_soft_yielding_and_sweeping}
\end{figure}

We begin with a mathematical description of a single spring, so that in the later sections many such springs could be assembled to a lattice. At each moment of time the state of the spring is characterized by the following scalar variables: the {\it total elongation} $x$, the {\it elastic elongation} $\varepsilon$, the  {\it plastic elongation} $p$, the {\it stress} $\sigma$ and the {\it damage variable} $a$. The latter is also called the {\it internal variable}. The usual constitutive laws are:
\begin{align}
\text{additive decomposition: }&& x&=\varepsilon + p,&\label{eq:additive_decomposition_single}\\
\text{Hooke's law: }&& \sigma&= k \varepsilon,& \label{eq:Hooke's_law_single}
\end{align}
where parameter $k>0$ is the stiffness of the spring. We assume the following plastic flow rule to hold for a.a. $t$ from a given time-interval $[0,T]$:
\begin{equation}
\frac{d}{dt}\begin{pmatrix}p\\ -a\end{pmatrix}\in N_{{\rm C}+\begin{pmatrix}0\\ h a\end{pmatrix}}\begin{pmatrix}\sigma \\ a\end{pmatrix},
\label{eq:plastic_flow_single}
\end{equation}
In the current text we focus on the linear case and put  ${\rm C}$ as 
\begin{equation}
{\rm C}:=\left\{\begin{pmatrix}\sigma\\ \alpha\end{pmatrix}\in \mathbb{R}^2: - (c_{0}+s\alpha)\leqslant \sigma\leqslant c_{0}+s\alpha \right\}.
\label{eq:tilde_C_single}
\end{equation}
We use $\alpha$ to distinguish the dummy variable in the definition of ${\rm C}$ from the actual state variable $a$ (they also differ by a translation on $ha$, so that they have different meanings). The corresponding state-dependent set in \eqref{eq:plastic_flow_single} is illustrated in Fig.~\ref{fig:lin_soft_yielding_and_sweeping} a). Here we assume that $h\in \mathbb{R}$, $s \geqslant 0$ and $c_{0}>0$ are the parameters of the spring which describe its behavior in the plastic regime. We call $h, s$ and $c_0$ the {\it coefficient of state-dependent feedback}, the {\it geometric slope} and the {\it initial yield stress} respectively. 

The admissible stresses of the spring are defined via the ``horizontal'' section of Fig.~\ref{fig:lin_soft_yielding_and_sweeping} a) at a specified value of $a$:
\[
\left\{\begin{pmatrix}\sigma\\a\end{pmatrix}:\sigma\in \mathbb{R}\right\}\cap\left({\rm C}+\begin{pmatrix}0\\ ha\end{pmatrix}\right),
\]
cf. \cite[Fig.~1.6]{SimoHughes1998}, i.e. with ${\rm C}$ given as \eqref{eq:tilde_C_single} we have the set of admissible stresses
\[
\sigma\in [- (c_{0}+s(1-h)a,\, c_{0}+s(1-h)a].
\]

\subsection{A state-dependent sweeping process for a single spring}
Before we proceed to the construction of lattices we shall explain the behavior of a single spring, described by \eqref{eq:additive_decomposition_single}--\eqref{eq:tilde_C_single}. To such a spring we apply a displacement loading (see Fig. \ref{fig:intro_fig_models} a)
\begin{equation}
x = l(t).
\label{eq:loading_single}
\end{equation}
From \eqref{eq:additive_decomposition_single}, \eqref{eq:Hooke's_law_single} and \eqref{eq:loading_single} we get that
\[
kp = k l(t) - \sigma.
\]
Then it follows from \eqref{eq:plastic_flow_single} and \eqref{eq:nc_basic_property2} that
\[
\frac{d}{dt}\begin{pmatrix}kl(t)-\sigma \\ -a\end{pmatrix}=\frac{d}{dt}\begin{pmatrix}k p\\ -a\end{pmatrix}=\mathbb{K}\,\frac{d}{dt}\begin{pmatrix}p\\ -a\end{pmatrix}\in \mathbb{K}\, N_{{\rm C}+\begin{pmatrix}0\\ h a\end{pmatrix}}\begin{pmatrix}\sigma \\ a\end{pmatrix}= N^{\mathbb{K}^{-1}}_{{\rm C}+\begin{pmatrix}0\\ h a\end{pmatrix}}\begin{pmatrix}\sigma \\ a\end{pmatrix}
\]
where 
\[
\mathbb{K}=\begin{pmatrix}k & 0 \\ 0 & 1\end{pmatrix}.
\]
Denote 
\begin{equation}
y := \sigma - k l(t)
\label{eq:y_change_var_single}
\end{equation}
and deduce that the single spring model \eqref{eq:additive_decomposition_single}--\eqref{eq:loading_single} is equivalent to the following inclusion in $\mathbb{R}^2$:
\begin{equation}
-\frac{d}{dt}\begin{pmatrix}y\\a\end{pmatrix}\in N^{\mathbb{K}^{-1}}_{{\rm C}+\begin{pmatrix}-k l(t)\\ h a\end{pmatrix}} \begin{pmatrix}y \\a\end{pmatrix}.
\label{eq:sp_single}
\end{equation}
Inclusion \eqref{eq:sp_single} is a state-dependent sweeping process.
\subsection{Analysis of the plastic behavior of the single spring}
Now we can analyze the behavior of the single spring by looking at the process \eqref{eq:sp_single}. For convenience, we first rewrite it in the compact form:
\begin{equation}
-\frac{d}{dt}\begin{pmatrix}y\\a\end{pmatrix}\in N^{\mathbb{K}^{-1}}_{\mathcal{C}(t, a)} \begin{pmatrix}y \\a\end{pmatrix}.
\label{eq:sp_single2}
\end{equation}
The following is the definition of state-dependent moving set $\mathcal{C}(t, a)$, in which we must distinguish between the parameter $\widetilde{a}$ and the dummy variables $a, \alpha$.
\begin{multline}
\mathcal{C}(t, \widetilde{a}):=\left\{\begin{pmatrix}\sigma - k l(t)\\ \alpha+h \widetilde{a}\end{pmatrix}:\begin{pmatrix}\sigma\\ \alpha\end{pmatrix}\in \mathbb{R}^2,\, - (c_{0}+s\alpha)\leqslant \sigma\leqslant c_{0}+s\alpha \right\}=\\
=\left\{\begin{pmatrix}y\\ a\end{pmatrix}\in \mathbb{R}^2: - (c_{0}+s(a-h\widetilde{a}))\leqslant y+k l(t)\leqslant c_{0}+s(a-h\widetilde{a}) \right\}=\\
=\left\{\begin{pmatrix}y\\ a\end{pmatrix}\in \mathbb{R}^2: \begin{array}{r}-(y+k l(t))- s(a-h\widetilde{a})\leqslant c_{0},\\(y+k l(t))- s(a-h\widetilde{a})\leqslant c_{0}\phantom{,}\end{array}\right\}=\\
=\left\{\begin{pmatrix}y\\ a\end{pmatrix}\in \mathbb{R}^2: \begin{array}{r}\eta_{-}^{\top}\mathbb{K}^{-1}\begin{pmatrix}y+k l(t)\\ a-h\widetilde{a}\end{pmatrix}\leqslant c_{0},\\\eta_{+}^{\top}\mathbb{K}^{-1} \begin{pmatrix}y+k l(t)\\ a-h\widetilde{a}\end{pmatrix}\leqslant c_{0}\phantom{,}\end{array}\right\},
\label{eq:tilde_C_single2}
\end{multline}
where
\begin{equation}
\eta_{-}:= \begin{pmatrix}-k\\ -s \end{pmatrix},\qquad \eta_{+}:= \begin{pmatrix}k\\ -s \end{pmatrix},
\label{eq:normals_single}
\end{equation}
see Fig. \ref{fig:lin_soft_yielding_and_sweeping} b). It is clear, that the solution of such sweeping process can be in one of the four following regimes, see also Fig. \ref{fig:lin_soft_yielding_and_sweeping} c):
\begin{itemize}
\item {\it elastic regime}, when $(y,a)\in {\rm int}\, \mathcal{C}(t,a)$,
\item the {\it plastic regime} of maximal stress, when the constraint with the normal $\eta_+$ is active,
\item the {\it plastic regime} of minimal stress, when the constraint with the normal $\eta_-$ is active,
\item the {\it state of complete failure}, when both constraints are active. Such state means, that the current set of admissible stresses is degenerate, i.e. it is the singleton-zero set $\{0\}$.
\end{itemize}
In the elastic regime formulas \eqref{eq:y_change_var_single}--\eqref{eq:sp_single} simply mean that
\[
\dot y(t) =0, \qquad \dot a(t)=0, \qquad \sigma(t) = \sigma(0)+ kl(t).
\] 
We explicitly solve the evolution in the plastic regimes over the intervals of monotonic load $l(t)$ in the following proposition, where we mean the signs $\pm$ and $\mp$ to be,  respectively, $+$ and $-$, or $-$ and $+$ to combine both cases in one statement.
\begin{Proposition}
\label{prop:sp_single} 
Let 
\begin{equation}
k^{-1}s^2(1-h)>-1.
\label{eq:solvable_modulus_single}
\end{equation}
If $y_{0}, a_{0}$ and $t_0$ are such that the constraint with $\eta_{\pm}$ is active, i.e. that 
\begin{equation}
\eta_{\pm}^{\top}\mathbb{K}^{-1} \begin{pmatrix}y_{0}+k l(t_0)\\ (1-h)a_{0}\end{pmatrix} = c_{0}.
\label{eq:together_with_the_set_0_single}
\end{equation}
and $\pm\dot l(t)\geqslant 0$ for $t\geqslant t_0$, then
\begin{equation}
\begin{pmatrix}
y(t)\\ a(t)
\end{pmatrix} := \begin{pmatrix}
y_{0}\\ a_{0}
\end{pmatrix} \mp \frac{1}{1+k^{-1}s^2(1-h)}\eta_{\pm} l(t)
\label{eq:solution1_single}
\end{equation}
is a solution of the state-dependent sweeping process
\eqref{eq:sp_single2}--\eqref{eq:normals_single} as long as the other constraint remains inactive, i.e.
\begin{equation}
\eta_{\mp}^{\top}\mathbb{K}^{-1} \begin{pmatrix}y+k l(t)\\ (1-h)a\end{pmatrix} < c_{0}.
\label{eq:second_not_active_single}
\end{equation}
\end{Proposition}
\noindent{\bf Proof.} We need to find a non-negative $\lambda(t)$ such that
\begin{equation}
-\begin{pmatrix}\dot  y(t)\\ \dot a(t) \end{pmatrix}=\lambda(t)\, \eta_{\pm}
\label{eq:normal_velocity_single}
\end{equation}
and 
\begin{equation}
\eta_{\pm}^{\top}\mathbb{K}^{-1}\, \begin{pmatrix}\dot y(t) + k \dot l(t)\\ (1-h)\dot a(t) \end{pmatrix}= 0.
\label{eq:together_with_the_set_single}
\end{equation}
Finding such $\lambda$ will be sufficient to guarantee that
\[
\begin{pmatrix} y\\ a \end{pmatrix}= \begin{pmatrix} y_{0}\\ a_{0} \end{pmatrix} + \int_{t_0}^t \begin{pmatrix}\dot  y(\tau)\\ \dot a(\tau) \end{pmatrix} d\tau
\]
is a solution, because \eqref{eq:together_with_the_set_single} and \eqref{eq:together_with_the_set_0_single} imply that the constraint with $\eta_{\pm}$ is active for all $t\geqslant t_0$, so  that $(y,a)\in \mathcal{C}(t,a)$ as long as the constraint with $\eta_{\mp}$ is not active, while \eqref{eq:normal_velocity_single} means that the inclusion \eqref{eq:sp_single2} holds. 

\begin{figure}[H]\center
\includegraphics{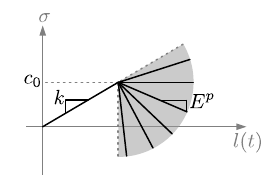}
\caption{
\footnotesize The range $(-\infty,k)$ of plasticity moduli $E^p$ which we consider in this paper, cf. \cite[p. 88]{SimoHughes1998}. A plasticity modulus $E^p$ is computed by \eqref{eq:modulus_formula_hs} from the parameters in the state-dependent plastic flow rule \eqref{eq:plastic_flow_single}--\eqref{eq:tilde_C_single} and in Hooke's law \eqref{eq:Hooke's_law_single}. 
} \label{fig:possible_moduli}
\end{figure}

\begin{figure}[H]\center
\includegraphics{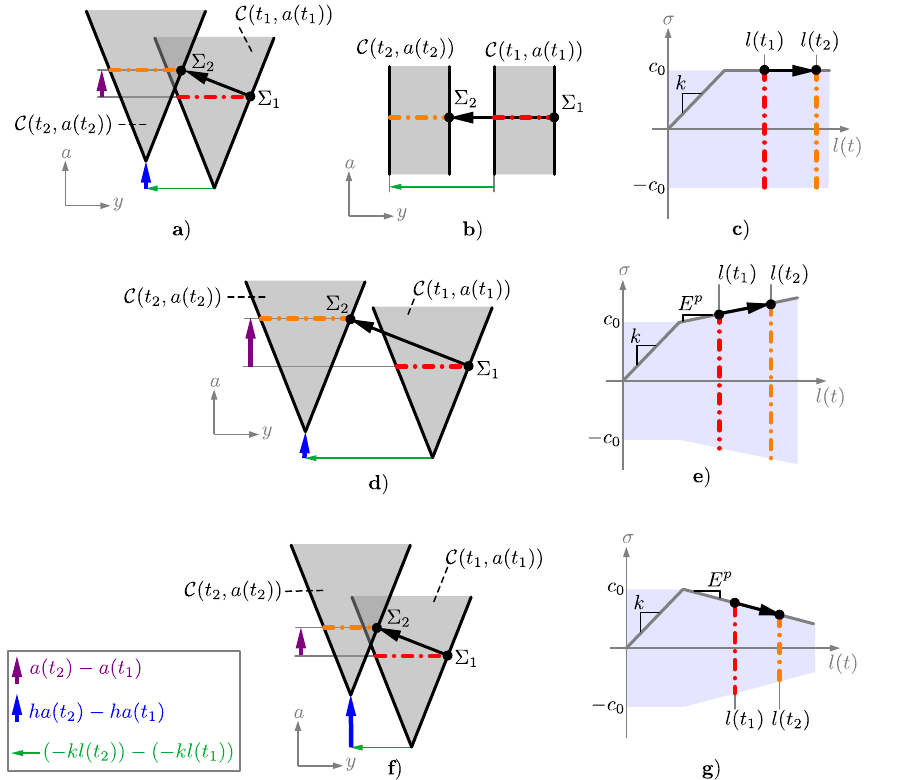}
\caption{
\footnotesize Evolution in the state-dependent sweeping process \eqref{eq:sp_single2}--\eqref{eq:tilde_C_single2} and the corresponding displacement-stress plots for different types of linear isotropic plasticity of an isolated spring (Fig. \ref{fig:intro_fig_models} a) under monotonically increasing displacement load $l(t)$. Here $t_1<t_2, \Sigma_1 := (y(t_1),a (t_1)), \Sigma_2 := (y(t_2),a (t_2))$. {\bf a)} $h=1$ leads to perfect plasticity, because $a(t_2)-a(t_1)$ (purple) equals to the term $h a(t_2)- h a(t_1)$ of state-dependent translation (blue) of the moving set $\mathcal{C}$, which leads to the set of admissible stresses being constant in time. Indeed, observe, that the set of admissible stresses at $t_1$  (red) is the same as at $t_2$ (orange). {\bf b)} $s=0$ leads to perfect plasticity as well, because damage variable $a$ cannot change, and set of admissible stresses remains the same between $t_1$  (red) and $t_2$ (orange). {\bf c)} Displacement-stress plot illustrates the case of perfect plasticity characterized by the set of admissible stresses remaining the same between $t_1$ (red) and $t_2$ (orange). {\bf d)} Parameters $s>0$ and $h<1$ lead to hardening: notice how the term $h a(t_2)- h a(t_1)$ of state-dependent translation (blue) of the moving set is smaller than the change $a(t_2)- a(t_1)$ (purple) of the state variable. Because of this the set of admissible stresses expands between $t_1$ (red) and $t_2$ (orange). {\bf e)} Displacement-stress plot for hardening, which is happening when the set of admissible stresses expands between $t_1$ (red) and $t_2$ (orange). {\bf f)} Parameters $s>0$ and $h>1$  lead to softening: the term $h a(t_2)- h a(t_1)$ of state-dependent translation (blue) of the moving set is greater than the change $a(t_2)- a(t_1)$ (purple) of the state variable. This results in the set of admissible stresses shrinking between $t_1$ (red) and $t_2$ (orange).  {\bf g)} Displacement-stress plot for softening, which is happening when the set of admissible stresses shrinks between $t_1$ (red) and $t_2$ (orange).
} \label{fig:plasticity_types}
\end{figure}

From \eqref{eq:normal_velocity_single} we obtain
\[
\begin{pmatrix}
\dot y(t)\\ \dot a(t)
\end{pmatrix}
=
\begin{pmatrix}
\mp k \lambda(t)\\
s \lambda(t)
\end{pmatrix},
\]
then we plug it into \eqref{eq:together_with_the_set_single}:
\[
\begin{pmatrix} \pm k \\ -s\end{pmatrix}^{\top} \begin{pmatrix}k^{-1} & 0 \\ 0 & 1\end{pmatrix}\begin{pmatrix} \mp\lambda(t) k + k  \dot l(t)\\ (1-h)s\lambda(t) \end{pmatrix}=0,
\]
\[
-k\lambda(t) \pm k \dot l(t) -(1-h)s^2 \lambda(t) =0,
\]
\[
\lambda(t) = \pm\frac{k\dot l(t)}{k+(1-h)s^2}=\pm\frac{1}{1+k^{-1}s^2(1-h)}\,\dot l(t).
\]
Hence
\[
\begin{pmatrix}\dot y(t)\\ \dot a(t) \end{pmatrix}= \mp \frac{1}{1+k^{-1}s^2(1-h)}\,\dot l(t) \eta_{\pm}.
\]
Observe that \eqref{eq:solvable_modulus_single} and $\pm\dot l(t)\geqslant 0$ guarantee that $\lambda(t)\geqslant0$. Therefore we have shown that \eqref{eq:solution1_single} is indeed a solution of \eqref{eq:sp_single2}--\eqref{eq:normals_single} under the conditions listed.
$\blacksquare$
\begin{Corollary}
Under the conditions of Proposition \ref{prop:sp_single} we have the evolution of the stress as the following:
\begin{equation*}
\sigma(t) = y(t) +k l(t) = y_{0} - \frac{k}{1+k^{-1}s^2(1-h)} l(t)+k l(t)=\sigma_{0}-kl(t_0)+\frac{s^2(1-h)}{1 + k^{-1} s^2(1-h)}l(t).
\end{equation*}
\end{Corollary}
These computations allow us to characterize the behavior of the spring as follows: 
\begin{Observation} \label{obs:single_spring} Constitutive equations \eqref{eq:additive_decomposition_single}--\eqref{eq:tilde_C_single} describe a spring with elasticity modulus $k$ in the elastic regime. In the plastic regimes, we can have different types of plasticity, depending on the parameters:
\begin{itemize}
\item if $h=1$ (Fig. \ref{fig:plasticity_types} a) or $s=0$ (Fig.  \ref{fig:plasticity_types} b) then it is the case of {\it perfect plasticity} (Fig.  \ref{fig:plasticity_types} c),
\item if $s^2(1-h)>0$ then it is the case of linear isotropic {\it hardening} (Fig.  \ref{fig:plasticity_types} d, e) with {\it plasticity modulus}
\begin{equation}
E^p=\frac{s^2(1-h)}{1 + k^{-1} s^2(1-h)},
\label{eq:modulus_formula_hs}
\end{equation}
\item if $0>k^{-1}s^2(1-h)>-1$ then it is linear isotropic {\it softening} with plasticity modulus $E^p$ (Fig.  \ref{fig:plasticity_types} f, g).
\end{itemize} 
Notice that if $k^{-1}s^2(1-h) \to -1$ from above then we have $E^p\to -\infty$, so we exclude the situation of $k^{-1}s^2(1-h) \geqslant - 1$ from further consideration in this paper. Furthermore, if $s^2(1-h)\to +\infty$ with $k$ fixed, then $E^p \to k$ from below, i.e. the limit of hardening moduli is elastic stiffness of the spring, see Fig. \ref{fig:possible_moduli}. 
\end{Observation}

\begin{figure}[H]\center
\includegraphics{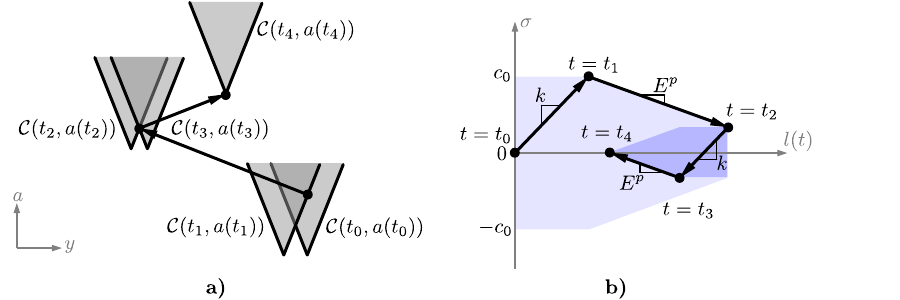}
\caption{
\footnotesize Modeling of the evolution of an elasto-plastic spring with linear isotropic softening, starting from the relaxed initial state. At $t=t_0$ we start from the stress-free state and apply increasing displacement loading $l(t)$. At $t=t_1$ the spring changes from elastic to plastic regime. At $t=t_2$ the derivative of the applied displacement load is reversed ($l(t)$ is now decreasing), which makes the spring to return to the elastic regime and undergo {\it elastic unloading} \cite[Sect. 3.1]{Han2012}. At $t=t_3$ the lower elastic threshold is attained, so the spring switches to the plastic mode again. Notice that the spring has reduced set of admissible stresses compared to $t_1$, and this applies to both upper and lower thresholds (isotropic softening).  At $t=t_4$ the spring reaches a point of complete failure, because its set of admissible stresses had reduced to $\{0\}$, and the solution cannot be continued further. {\bf a)} A sketch of the evolution in the state-dependent sweeping process \eqref{eq:sp_single2}. Note, that during the elastic phases the solution of the sweeping process is in the interior of the moving set and remains at rest. {\bf b)} The corresponding displacement--stress plot. The size of the blue area along the vertical axis indicates the set of admissible stresses for the corresponding $l(t)$.
} \label{fig:complete_failure}
\end{figure}

If one is interested in modeling perfect plasticity and hardening phenomena (but not softening), it is possible to construct a sweeping process analogous to \eqref{eq:sp_single2} with the moving set depending on $t$ only, see the references at the end of Section \ref{ssec:sp_intro}. In this paper our main interest is softening, however, we do not exclude the other two cases, so that we are able to compare between them. Also, we preserve the opportunity to include perfectly plastic and hardening springs in the same lattice with softening springs.
\begin{Remark} If a spring with linear softening plasticity reaches a state where both constraints of \eqref{eq:tilde_C_single2} are active at $(y, a)$ (see e.g. the position at $t=t_4$ on Fig. \ref{fig:complete_failure}, both a) and b)), then it means complete failure of the spring, as for a further change of the load the solution is not continuable. This could be viewed as separation of the material at the endpoints of the spring due to ductile fracture. 
\label{remark:failure1}
\end{Remark}


\section{Lattice Spring Model and the corresponding state-dependent sweeping process}
\label{sect:LSMs}
Since we have understood the behavior of a single spring in the model \eqref{eq:additive_decomposition_single}--\eqref{eq:tilde_C_single}, we can connect
together $m$ such springs to form a lattice. It is important to note that we will not not alter the constitutive laws of individual springs, but merely couple the springs at the nodes of the lattice and write the entire problem in a vector form. Therefore, each spring in the lattice still behaves according to the analysis in the previous Section \ref{sect:single_spring}, and  Observation \ref{obs:single_spring} in particular.

A lattice is described by a given graph with incidence matrix $Q$ and given reference positions of the vertices $\xi_0$. First, we write the system of governing equations of the entire lattice and then convert it to a state-dependent sweeping process. 
\subsection{The system of equations for a lattice of connected springs}
This section mostly repeats the construction of \cite[Section 3]{Gudoshnikov2023preprint}, where the reader can find additional details and discussions on the mechanics of the lattices in general, rather then on the question arising from softening and state-dependence.
\subsubsection{Constitutive laws in a vector form}
We collect the variables of individual springs into vectors, so the state of the entire lattice is
described by the unknown time-dependent variables $x, \varepsilon, p, \sigma, a \in \mathbb{R}^m$ with the same meaning of their individual components as in Section \ref{ssect:single-spring-model}. Clearly, by
collecting the equations  \eqref{eq:additive_decomposition_single}--\eqref{eq:Hooke's_law_single} of the individual components to vector forms
we get
\begin{equation}
x = \varepsilon + p,
\tag{LSM1}
\label{eq:additive_decomposition}
\end{equation}
\begin{equation}
\tag{LSM2}
\sigma = K \varepsilon,
\label{eq:Hooke's_law}
\end{equation}
where $K$ is the following diagonal matrix with positive diagonal elements: 
\[
K={\rm diag}\, (k_i)_{i\in \overline{1,m}}. 
\]
Also, by using Lemma \ref{lemma:apa-normal-cones} with $E=\mathbb{R}^2$ and rearranging the order of components, we observe that \eqref{eq:plastic_flow_single} holds for all $i\in \overline{1,m}$ if and only if
\begin{equation}
\frac{d}{dt}\begin{pmatrix}p\\-a\end{pmatrix}\in N_{{\rm C}+\begin{pmatrix}0\\ Ha\end{pmatrix}}\begin{pmatrix}\sigma \\ a\end{pmatrix},
\tag{LSM3}
\label{eq:plastic_flow}
\end{equation}
where $H$ is the following diagonal matrix:
\[
H={\rm diag}\, (h_i)_{i\in \overline{1,m}},
\]
and ${\rm C}$ is given as
\begin{equation}
{\rm C}=\left\{\begin{pmatrix}\sigma \\ \alpha\end{pmatrix}\in \mathbb{R}^{2m}: - (c_{0,i}+s_i\alpha_i)\leqslant \sigma_i\leqslant c_{0,i}+s_i\alpha_i  \text{ for all }i \in \overline{1,m}\right\}.
\label{eq:tildeC}
\end{equation}
We can rewrite the state dependent set in \eqref{eq:plastic_flow} as the following (again, we distinguish the parameter variable $\widetilde{a}$ from dummy variables $a$ and $\alpha$)
\begin{multline}
{\rm C} +\begin{pmatrix}0 \\ H\widetilde{a}\end{pmatrix} = \left\{\begin{pmatrix}\sigma \\ \alpha+H\widetilde{a}\end{pmatrix}:\begin{pmatrix}\sigma \\ \alpha\end{pmatrix}\in \mathbb{R}^{2m}, - (c_{0,i}+s_i\alpha_i)\leqslant \sigma_i\leqslant c_{0,i}+s_i\alpha_i  \text{ for all }i \in \overline{1,m}\right\} =\\
=\left\{\begin{pmatrix}\sigma \\ a\end{pmatrix}\in \mathbb{R}^{2m}: - (c_{0,i}+s_i(a_i-h_i\widetilde{a}_i))\leqslant \sigma_i\leqslant c_{0,i}+s_i(a_i-h_i\widetilde{a}_i)  \text{ for all }i \in \overline{1,m}\right\}=\\
=\left\{\begin{pmatrix}\sigma \\ a\end{pmatrix}\in \mathbb{R}^{2m}: \begin{array}{r}-\sigma_i-s_i(a_i-h_i\widetilde{a}_i) \leqslant c_{0,i} \\ \sigma_i-s_i(a_i-h_i\widetilde{a}_i) \leqslant c_{0,i} \end{array} \text{ for all }i \in \overline{1,m}\right\}
\label{eq:plastic_flow_set}
\end{multline}

\begin{Remark} 
\label{remark:a_nondecreasing1}
We only consider $s_i\geqslant 0,\, i\in \overline{1,m}$, and this makes each component $a_i$ a non-decreasing function. Indeed, for each value of parameter $\widetilde{a}$ set \eqref{eq:plastic_flow_set} is a convex polytope with its normal cone
\begin{multline}
 N_{{\rm C}+\begin{pmatrix}0\\ H\widetilde{a}\end{pmatrix}}\begin{pmatrix}\sigma \\ a \end{pmatrix}={\rm cone}\left(\left\{\begin{pmatrix}-e_i\\ -s_i e_i\end{pmatrix}: i\in \overline{1,m}: \begin{pmatrix}-e_i\\ -s_i e_i\end{pmatrix}^{\top}\begin{pmatrix}\sigma \\ a - H\widetilde{a}\end{pmatrix}=c_{0,i} \right\}\right.\cup\\ \cup\left. \left\{\begin{pmatrix}e_i\\ -s_i e_i\end{pmatrix}:i\in \overline{1,m}: \begin{pmatrix}e_i \\ -s_i e_i\end{pmatrix}^{\top}\begin{pmatrix}\sigma \\ a  - H\widetilde{a}\end{pmatrix}=c_{0,i}\right\}\right),
\label{eq:plastic_flow_normal_cone}
\end{multline}
where by $e_i\in \mathbb{R}^m$ we denote the $i$-th vector of the standard basis.
Therefore, it follows from \eqref{eq:plastic_flow} that
\[
-\dot a_i =-s_i \lambda_i \text{ for some }\lambda_i\geqslant 0,
\]
and
\begin{equation}
\dot a_i\geqslant 0.
\label{eq:a_nondecreasing}
\end{equation}
for a.a. $t\in [0,T]$.
\end{Remark}
\begin{Remark}
We can explicitly compute the state of complete failure we mentioned in Remark \ref{remark:failure1}. Observe that if for some $\begin{pmatrix}\sigma\\a\end{pmatrix}\in \mathbb{R}^{2m}$ and some $i\in \overline{1,m}$ we simultaneously have both constraints in \eqref{eq:plastic_flow_set} active for $a=\widetilde{a}$, i.e
\begin{equation}
\begin{array}{r}-\sigma_i-s_i(1-h_i)a_i = c_{0,i} \\ \sigma_i-s_i(1-h_i)a_i = c_{0,i} \end{array}
\label{eq:failure_equations1}
\end{equation}
then it necessarily means that 
\[
a_i=-\frac{c_{0,i}}{s_i(1-h_i)}, \qquad \sigma_i=0
\]
which, due to \eqref{eq:a_nondecreasing}, is reachable only when $h_i>1$. In such a case, for any greater value of $a_i$ there are no possible values of $\sigma_i$.
Again by \eqref{eq:a_nondecreasing} we conclude that after reaching the point of complete failure of a spring the solution is either not continuable, or there is no more evolution in the respective component of the solution: 
\begin{equation}
\begin{array}{c}
\text{\eqref{eq:failure_equations1} is true for some }i\in \overline{1,m},\, t^*\in[0,T],\\
h_i>1,\\
\text{but the solution exists on }[0,T]
\end{array}\quad\Longrightarrow\quad \text{for a.a. }t\geqslant t^* \text{ we have }
\dot a_i=\dot \sigma_i=0.
\label{eq:failure_means_no_change}
\end{equation}
\label{remark:failure2}
\end{Remark}
\subsubsection{Geometric compatibility, displacement load and the well-posedness assumptions}
\label{ssect:kinematics}
To reflect the fact that the springs are required to stay connected according to the given graph at all times, we proceed to construct the linearized compatibility equation.

To each of the $m$ springs of the lattice we assign an arbitrary orientation, by calling one of its endpoints the {\it origin} and the other one the {\it terminus} \cite[Ch. 7]{ClarkHoltonGraphTheory1991}, constructing a {\it directed graph} in this manner. The structure of the lattice can then be described by an  $n\times m$ incidence matrix $Q$ of the directed graph. The incidence matrix is constructed as the following (see more details in \cite{Bapat2010}): for $i\in \overline{1,m},\, j\in\overline{1,n}$ set
\begin{itemize}
\item $Q_{ji}=0$ \text{ if none of the endpoints of spring $i$ is node $j$},
\item $Q_{ji}=1$ \text{ if node $j$ is the origin of the spring $i$ according to the assigned orientation,} 
\item $Q_{ji}=-1$ \text{ if node $j$ is the terminus of the spring $i$ according to the assigned orientation.} 
\end{itemize}
	At any particular moment the coordinates of the vertices can be collected into a vector $\xi\in \mathbb{R}^{nd}$ so that $\xi_{d(j-1)+k}$ is the $k$-th coordinate of the node $j$ (where $j\in \overline{1,n}, k\in \overline{1,d}$). The length of spring $i$ (where $i\in \overline{1,m}$) can then be calculated as the norm of the vector
\[(\text{the terminus of spring }i)-(\text{the origin of spring }i)=-\left(\sum_{j=1}^n Q_{ji}\xi_{d(j-1)+k}\right)_{k\in \overline{1, d}}\]
\noindent and the lengths of all $m$ springs can be calculated by the function
\[
 \varphi:\mathbb{R}^{nd}\to \mathbb{R}^m
\]
\begin{equation}
\varphi(\xi)=\left(\varphi_i(\xi)\right)_{i\in\overline{1,m}}:=\left(\sqrt{\sum_{k=1}^d\left(\sum_{j=1}^n Q_{ji}\xi_{d(j-1)+k}\right)^2}\right)_{i\in\overline{1,m}}.
\label{eq:distanceFunction}
\end{equation}

We linearize function $\varphi$ at a given fixed {\it reference configuration} $\xi_0=(\xi^0_{d(j-1)+k})\in \mathbb{R}^{nd}$ and write down the {\it geometric constraint} connecting $x$ and the vector of {\it displacements} $\zeta\in \mathbb{R}^{nd}$ of the nodes from the reference configuration:
\begin{equation}
\left(D_{\xi_0}\varphi\right) \zeta=x
\tag{LSM4}
\label{eq:gc}
\end{equation}
 
Here both $x$ and $\zeta$ are unknown variables but we assume that displacements $\zeta$ are small compared to the lengths of the springs, therefore we can use fixed $\xi_0$. Jacobi matrix $D_{\xi_0}\varphi$ in \eqref{eq:gc}  is independent on $t$, and it can be written explicitly. Specifically, the $(i, d(j-1)+k)$ entry of $D_{\xi_0} \varphi$ is given by

\begin{equation}
\left. \frac{\partial \varphi_i}{\partial \xi_{d(j-1)+k}} \right|_{\xi=\xi_0}= \frac{\left(\sum\limits_{\bar j=1}^n Q_{\bar ji}\xi^0_{d(\bar j-1)+k}\right)Q_{j i}}{\sqrt{\sum\limits_{\bar k=1}^d\left(\sum\limits_{\bar j=1}^n Q_{\bar ji}\xi^0_{d(\bar j-1)+\bar k}\right)^2}} = \mathcal{D}_{ik}Q_{ji},
\label{eq:D_xi_phi}
\end{equation}
where $\mathcal{D}$ is the $m\times d$-matrix with $(i,k)$ entry
\begin{equation}
\mathcal{D}_{ik}=\frac{1}{\varphi_i(\xi_0)}\sum\limits_{\bar j=1}^n Q_{\bar ji}\xi^0_{d(\bar j-1)+k}.
\label{eq:unit_vectors_along_springs}
\end{equation}
Observe that $i$-th row of $\mathcal{D}$ is the {\it unit vector in the direction from the terminus to the origin of spring $i$ in reference configuration $\xi_0$}, i.e. the direction of such unit vector is opposite to the chosen orientation in the geometric directed graph, corresponding to $Q$ with the nodes placement $\xi_0$. 

Equation \eqref{eq:gc} appears in the literature \cite[(2.6)]{Lubensky2015} and \cite[(3.17)]{Moreau1973} (in the form for individual springs).

Along with the geometric constraint we introduce an additional constraint to play the role of
the boundary condition. It has the form
\begin{equation}
R(\xi_0 + \zeta) + r(t)=0,
\tag{LSM5}
\label{eq:bc}
\end{equation}
where $R$ is a given $q\times nd$-matrix and $r$ is a given function of time with $q$-vector values for some $q\in \mathbb{N}$. We call equation \eqref{eq:bc}, function $r$ and number $q$, respectively, the {\it external displacement constraint},  the {\it displacement load} and the {\it number of external displacement constraints}. In turn, vectors $\zeta\in \mathbb{R}^{nd}$ satisfying \eqref{eq:bc} we call {\it feasible displacements}.

For well-posedness of the problem we need the following assumptions about given data $R$ and $D_{\xi_0} \varphi$ to hold:

\begin{Assumption}
{\it Well-definiteness of the external constraint}, i.e. matrix $R$ has maximum row rank:
\begin{equation}
{\rm rank}\,  R = q,
\label{eq:rowsR_indep}
\end{equation}
\label{ass:rowsR_indep}
\end{Assumption}

\begin{Assumption}
{\it Kinematic determinacy}, which means that it is possible to uniquely solve \eqref{eq:gc}--\eqref{eq:bc} for displacements $\zeta$ (hence, for the positions $\xi_0+\zeta$ of the nodes) when elongation values $x$ are already obtained. In the next Section \ref{sssecct:equation_of_equilibrium} this condition also allows us to apply external forces of any direction. Kinematic determinacy can be written mathematically as the following equivalent conditions:
\begin{equation}
{\rm rank}\,  \begin{pmatrix} (D_{\xi_0} \varphi)^{\top} & R^{\top}\end{pmatrix} = nd,
\label{eq:trivial_kernel_intersect}
\end{equation}
\begin{equation*}
{\rm Ker}\, \begin{pmatrix} D_{\xi_0} \varphi \\  R\end{pmatrix} =\{0\},
\end{equation*}
\begin{equation}
{\rm Ker}\, (D_{\xi_0} \varphi) \cap {\rm Ker}\, R = \{0\}.
\label{eq:equivalent_condition}
\end{equation}
\label{ass:trivial_kernel_intersect}
\end{Assumption}

Matrices $D_{\xi_0} \varphi$ and $ \begin{pmatrix} D_{\xi_0} \varphi \\  R\end{pmatrix}$ we call, respectively, {\it the compatibility matrix} and {\it the enhanced compatibility matrix} of the lattice defined by $Q$, $\xi_0$, $R$ and $r$. In turn, the matrices $(D_{\xi_0} \varphi)^{\top}$, $\begin{pmatrix} (D_{\xi_0} \varphi)^{\top} & R^{\top}\end{pmatrix}$ we call, respectively {\it the equilibrium matrix} and {\it the enhanced equilibrium matrix} of the lattice. The image of the map $\varphi$ \eqref{eq:distanceFunction} is a nonlinear surface, known as the {\it Cayley-Menger variety} \cite[p.~435]{Fuhrmann2015}, and the linearization \eqref{eq:gc} to its tangent space (see \cite[Fig.~1, p.~78]{Admal2016}) plays the central role in rigidity theory \cite[Chapter 9]{Alfakih2018_rigidity}.  We refer to  \cite[Sect. 3]{Gudoshnikov2023preprint} for more references and a comprehensive overview of the algebraic properties \eqref{eq:gc}, \eqref{eq:bc} in relation to the kinematic properties of the lattice. Assumption \ref{ass:trivial_kernel_intersect} is also discussed there in the context of rigidity theory.
\subsubsection{Equation of equilibrium}
\label{sssecct:equation_of_equilibrium}
We model the {\it quasi-static} evolution,  meaning that the lattice stays at an equilibrium at all times. Specifically, we refer to the following
\begin{Definition}\cite[p. 17]{Goldstein2002}
\label{def:equilibrium}
A system of $m$ particles is said to be at an {\it equilibrium} when the total force on each particle vanishes. 
\end{Definition}
\noindent In our case the nodes of the lattice serve as particles, 
and the equation of equilibrium in the following form can be derived from \eqref{eq:gc}, \eqref{eq:bc} using the principle of virtual work (see \cite[Sect. 3]{Gudoshnikov2023preprint} for the complete derivation):
\begin{equation}
(D_{\xi_0} \varphi)^{\top} (-\sigma+F_1 f(t))\in {\rm Im}\,R^{\top},
\tag{LSM6}
\label{eq:balanceLaw} 
\end{equation}
where $f(t)\in \mathbb{R}^{nd}$ is a given function of time, which we call the {\it stress load} or the {\it external force}. For a fixed $t$ vector $f(t)$ has the following physical meaning: $f_{d(j-1)+k}$ is the $k$-th component of the force vector, applied to node $j$ ($j\in \overline{1,n}, k\in \overline{1,d}$). In turn, $F_1$ is an $(nd)\times m$-matrix, obtained as the following. The enhanced equilibrium matrix  $\begin{pmatrix} (D_{\xi_0} \varphi)^{\top} & R^{\top}\end{pmatrix}$ has dimensions $(nd)\times(m+q)$, and, by assumption \eqref{eq:trivial_kernel_intersect} it has full row rank. Therefore its Moore-Penrose pseudoinverse matrix $\begin{pmatrix} (D_{\xi_0} \varphi)^{\top} & R^{\top}\end{pmatrix}^+$ is also a right inverse and can be explicitly written \cite[Prop. 2.3]{Gudoshnikov2023preprint} as
\begin{equation}
\begin{pmatrix} (D_{\xi_0} \varphi)^{\top} & R^{\top}\end{pmatrix}^+=
\begin{pmatrix} F_1\\ F_2\end{pmatrix}=\begin{pmatrix} D_{\xi_0} \varphi \\ R\end{pmatrix}\left(\begin{pmatrix} (D_{\xi_0} \varphi)^{\top} & R^{\top}\end{pmatrix}\begin{pmatrix} D_{\xi_0} \varphi \\ R\end{pmatrix}\right)^{-1},
\label{eq:new_equilibrium_matrix_pseudoinverse}
\end{equation}
where $F_1$ and $F_2$ are $m\times(nd)$ and $q\times(nd)$ matrices, respectively.

Equation \eqref{eq:balanceLaw}  is a slightly modified version of equations from the literature \cite[(2.4)]{Lubensky2015}, \cite[(3.23)]{Moreau1973}.

\subsubsection{The complete set of equations of the elasto-plastic Lattice Spring Model with variable set of admissible stresses}
Naturally, we require the regularity of the displacement and stress load function:
\begin{Assumption}
\label{ass:ac_loads}
We assume that $r:[0,T]\to \mathbb{R}^{q}$ and $f:[0,T]\to \mathbb{R}^{nd}$ are absolutely continuous.
\end{Assumption}
We summarize the mathematical formulation of the problem of the quasistatic evolution of an elasto-plastic Lattice Spring Model. Given finctions $r, f$ (displacement and stress loads, respectively) find unknown functions $x,\varepsilon, p, \sigma, a:[0,T]\to \mathbb{R}^m, \zeta:[0,T]\to\mathbb{R}^{nd}$ such that for all $t\in [0,T]$
\begin{align}
\text{Additive decomposition:}&&x&=\varepsilon+p&\tag{\ref{eq:additive_decomposition}}\\
\text{Hooke's law:}&&\sigma&=K \varepsilon& \tag{\ref{eq:Hooke's_law}}\\
\text{Geometric constraint:}&&\left(D_{\xi_0} \varphi\right) \zeta&=x& \tag{\ref{eq:gc}}\\
\text{External displacement constraint:}&&R(\xi_0+\zeta) + r(t)&=0& \tag{\ref{eq:bc}}\\
\text{Equation of equilibrium:}&&(D_{\xi_0} \varphi)^{\top} (-\sigma+F_1 f(t))&\in {\rm Im}\,R^{\top} &\tag{\ref{eq:balanceLaw}}
\end{align}
and for a.a. $t\in [0,T]$
\begin{align}
\text{State-dependent plastic flow rule:}&&\frac{d}{dt}\begin{pmatrix}p\\-a\end{pmatrix}&\in N_{{\rm C}+\begin{pmatrix}0\\ Ha\end{pmatrix}}\begin{pmatrix}\sigma \\ a\end{pmatrix}.& \tag{\ref{eq:plastic_flow}}
\end{align}

\subsection{Derivation of the state-dependent sweeping process from the equations of the lattice}
\label{ssect:sweeping_process_general}
Before we proceed to the sweeping process, we will need some technical definitions and facts. Define the subspaces of $\mathbb{R}^m$
\begin{align}
\mathcal{U}& = K (D_{\xi_0}\varphi)\,{\rm Ker}\, R, \label{eq:U_def}\\
\mathcal{V}& = \{\sigma\in \mathbb{R}^m: (D_{\xi_0}\varphi)^{\top} \sigma \in {\rm Im}\, R^{\top}\}  \label{eq:V_def}
\end{align}
and introduce the inner product $\la\cdot, \cdot\ra_{K^{-1}}$ in $\mathbb{R}^m$ as
\begin{equation}
\la x,y \ra_{K^{-1}} = x^{\top} K^{-1}y.
\label{eq:Rm_inner_product}
\end{equation}
The following proposition establishes the essential properties of $\mathcal{U}$ and $\mathcal{V}$.
\begin{Proposition}
\label{prop:U_and_V_orth}
Provided that assumptions \eqref{eq:rowsR_indep}, \eqref{eq:trivial_kernel_intersect} hold true, the following is true as well.
\begin{enumerate}[\it i)]
\item \label{enum:prop_U_and_V_orth} The orthogonal complement of $\mathcal{V}\subset\mathbb{R}^m$ in the sense of the standard inner product in $\mathbb{R}^m$ is
\begin{equation}
\mathcal{V}^{\perp} = (D_{\xi_0}\varphi) \, {\rm Ker}\, R, 
\label{eq:V_perp_standard}
\end{equation}
and the orthogonal complement of $\mathcal{V}\subset\mathbb{R}^m$ in the sense of the inner product \eqref{eq:Rm_inner_product} is $\mathcal{U}$.
\item \label{enum:prop_U_and_V_equivalent_to_balance} The set of $\sigma$'s satisfying equilibrium equation \eqref{eq:balanceLaw} is an affine translation of space $\mathcal{V}$. Specifically, for any $\sigma\in \mathbb{R}^m$ and any $t$
\[
\sigma \text{ satisfies } \eqref{eq:balanceLaw} \quad \Longleftrightarrow \quad \sigma\in \mathcal{V} + F_1f(t).
\]
\item \label{enum:prop_U_and_V_dims} The dimensions of $\mathcal{U}$ and $\mathcal{V}$ are
\[
{\rm dim}\,\mathcal{U} =nd-q,\qquad {\rm dim}\,\mathcal{V}=m-nd+q.
\]
\end{enumerate}
\end{Proposition}
\noindent{\bf Proof.}
\begin{enumerate}[\it i)]
\item
Observe that
\begin{multline*}
((D_{\xi_0} \varphi)\,{\rm Ker}\, R)^\perp=\{x\in \mathbb{R}^m: x^{\top}(D_{\xi_0} \varphi)y=0 \text{ for all } y\in {\rm Ker}\, R\}=\\
=\{x\in \mathbb{R}^m: \left((D_{\xi_0} \varphi)^{\top} x\right)^{\top}y=0 \text{ for all } y\in {\rm Ker}\,R\}=\{x\in \mathbb{R}^m: (D_{\xi_0} \varphi)^{\top} x\in ({\rm Ker}\, R)^\perp\}=\\
=\{x\in \mathbb{R}^m:  (D_{\xi_0} \varphi)^{\top} x\in {\rm Im}\, R^{\top}\}= \mathcal{V},
\label{eq:U_perp}
\end{multline*}
which proves the first statement about the standard inner product. Since
\[
x\in \mathcal{V}^{\perp}\quad \Longleftrightarrow\quad  x^{\top}\sigma =0 \text{ for all }\sigma \in \mathcal{V}
\]
therefore
\[
Kx\in K\mathcal{V}^{\perp} \quad \Longleftrightarrow\quad (Kx)^{\top} K^{-1} \sigma =0 \text{ for all }\sigma \in \mathcal{V},
\]
which proves the statement about the inner product \eqref{eq:Rm_inner_product}. 
\item This can be verified by substituting $\sigma -F_1f(t)$ as $\sigma$ into definition \eqref{eq:V_def} of $\mathcal{V}$.
\item 
By rank-nullity theorem (see e.g. \cite[Th. 2.49]{OlverLinearAlgebra2018}) assumption \eqref{eq:rowsR_indep} implies that ${\rm dim\, Ker} \,R = nd - q$, therefore if we take a basis of ${\rm Ker} \,R$ and arrange its vectors as columns in a matrix $R_0$, we get $m \times (nd-q)$ matrix $R_0$ with linearly independent columns. Observe that it follows from definition \eqref{eq:U_def} of $\mathcal{U}$ and  from \eqref{eq:equivalent_condition} that ${\rm dim}\,\mathcal{U} = {\rm dim\, Ker}\, R = nd-q$. Indeed, assume that ${\rm dim}\,\mathcal{U} < {\rm dim\, Ker}\, R$. Then there is a nontrivial linear combination of basis vectors of ${\rm Ker}\, R$ (i.e. $R_0 z$ with some coefficients vector $z\in \mathbb{R}^{nd-q}\setminus\{0\}$) such that  $K (D_{\xi_0} \varphi) R_0 z = 0$, i.e.  $R_0 z$ is from ${\rm Ker}\, (D_{\xi_0} \varphi)$. Since $R_0 z$ is also from ${\rm Ker}\, R$ by construction, we have a contradiction with \eqref{eq:equivalent_condition}.
\end{enumerate}
$\blacksquare$

We call space $\mathcal{V}$ {\it the space of self-stresses} of the lattice, as it represents all such stress vectors $\sigma\in \mathbb{R}^m$ for which the corresponding forces, exerted by the springs to the nodes, either vanish or can be compensated by the reactions of the external displacement constraint \eqref{eq:bc}. For a more comprehensive discussion of the self-stresses and their connection to rigidity of the lattice we refer to \cite[Sect. 3]{Gudoshnikov2023preprint} and the references therein.

Space $\mathcal{V}$ is a central element of our further construction, and in order to have a meaningful problem we assume that $\mathcal{V}$ has a dimension of at least one (in addition to the assumptions \eqref{eq:rowsR_indep}--\eqref{eq:trivial_kernel_intersect}):
\begin{Assumption}
{\it Non-degneracy of the space of self-stresses:}
\begin{equation}
nd-q<m.
\label{eq:at_least_one_self-stress}
\end{equation}
\label{ass:at_least_one_self-stress}
\end{Assumption}

Put differently, \eqref{eq:at_least_one_self-stress} means that for each stress load $f(t)$ there is more than just one stress vector $\sigma\in \mathbb{R}^m$ which brings the lattice to the equilibrium, i.e. we assume that the lattice is not {\it statically determinate}.

For further use we need to define several matrices.
\[
\begin{array}{ll}
R^+ = R^{\top}(RR^{\top})^{-1}& \text{is the $nd \times q$ Moore-Penrose pseudoinverse matrix  of  $R$ (see \cite[Sect 2.2]{Gudoshnikov2023preprint}),}\\
U& \text{is an $m\times{\rm dim}\,\,\mathcal{U}$ matrix composed of columns which form a basis in $\mathcal{U}$},\\
V& \text{is an $m\times{\rm dim}\,\,\mathcal{V}$ matrix composed of columns which form a basis in $\mathcal{V}$}.
\end{array}
\]
Then there are matrices $P_U$ and $P_V$ such that $UP_U$ and $VP_V$ form a pair of projection maps, orthogonal in the sense of \eqref{eq:Rm_inner_product}. Namely (see e.g. \cite[Sect. 5.3 and 5.4]{OlverLinearAlgebra2018}),
\begin{equation}
P_U=\left(U^{\top} K^{-1} U\right)^{-1}U^{\top} K^{-1}\qquad P_V=\left(V^{\top} K^{-1} V\right)^{-1}V^{\top} K^{-1}.
\label{eq:projection_matrices}
\end{equation}

Finally, let 
\begin{align}
G &= VP_V K(D_{\xi_0}\varphi)R^+,\label{eq:G_def}\\
F&= UP_UF_1. \label{eq:F_def}
\end{align}
and now we can convert the problem \eqref{eq:additive_decomposition}--\eqref{eq:balanceLaw} to a state-dependent sweeping process.

\begin{Theorem}
\label{th:LSMtoMSP}
Let Assumptions \ref{ass:rowsR_indep}-\ref{ass:at_least_one_self-stress} hold true and let $s_i> 0$
for all $i\in \overline{1,m}$. The following statements are equivalent:
\begin{enumerate}[{\it i)}]
\item\label{enum:modeling_th_i}  Absolutely continuous functions $x,\varepsilon, p, \sigma, a:[0,T]\to \mathbb{R}^m, \, \zeta:[0,T]\to \mathbb{R}^{nd}$ solve \eqref{eq:additive_decomposition}--\eqref{eq:balanceLaw} and $y$ is defined as 
\begin{equation}
y=\sigma +Gr(t)-Ff(t).
\label{eq:y_to_sigma_change}
\end{equation}
\item \label{enum:modeling_th_ii} Absolutely continuous function 
 $\begin{pmatrix}y\\ a\end{pmatrix}:[0,T]\to \mathbb{R}^{2m}$ is a solution of state-dependent sweeping process
 \begin{equation}
-\frac{d}{dt}\begin{pmatrix} y \\a \end{pmatrix} \in N^{\mathbb{K}^{-1}}_{\mathcal{C}(t, a)}\begin{pmatrix} y \\a \end{pmatrix},
\label{eq:sp}
\end{equation}
where
\begin{equation}
\mathbb{K} = \begin{pmatrix}K & 0 \\ 0 & I\end{pmatrix},\qquad \mathcal{C}(t,a)= \left({\rm C}+\begin{pmatrix}Gr(t)-Ff(t)\\ Ha\end{pmatrix}\right)\cap \left(\mathcal{V}\times \mathbb{R}^m\right)
\label{eq:K,moving_set}
\end{equation}
where ${\rm C}$ is given by \eqref{eq:tildeC}, and
\begin{align}
\sigma & = y - Gr(t)+Ff(t), \label{eq:sigma_to_y_change}\\
\varepsilon & = K^{-1} \sigma,\label{eq:sigma_to_eps}\\ 
p(t) & = p(0)+\int\limits_0^t{\rm diag} \left(\frac{{\rm sign}\, \sigma_i(\tau)}{s_i}\right) \dot a (\tau)d\tau, \label{eq:p_recovery_formula}\\
x& = \varepsilon+p \tag{\ref{eq:additive_decomposition}} ,\\
\zeta& = \begin{pmatrix} D_{\xi_0} \varphi \\ R\end{pmatrix}^+\begin{pmatrix}x\\ -R\xi_0-r(t)\end{pmatrix}, \label{eq:zeta_recovery_formula}\\
\begin{pmatrix}D_{\xi_0} \varphi\\
R\end{pmatrix} \zeta (0) &= \begin{pmatrix} x(0) \\ -R \xi_0 - r(0) \end{pmatrix}.
\label{eq:compatible_initial_conditions}
\end{align}
where $\begin{pmatrix} D_{\xi_0} \varphi \\ R\end{pmatrix}^+$ is the Moore-Penrose pseudoinverse matrix of $\begin{pmatrix} D_{\xi_0} \varphi \\ R\end{pmatrix}$, which in our case (due to  \eqref{eq:trivial_kernel_intersect}) can be explicitly written as
\begin{equation}
\begin{pmatrix} D_{\xi_0} \varphi \\ R\end{pmatrix}^+ = \left(\begin{pmatrix} (D_{\xi_0} \varphi )^{\top} & R^{\top}\end{pmatrix} \begin{pmatrix} D_{\xi_0} \varphi \\ R\end{pmatrix}\right)^{-1}\begin{pmatrix}(D_{\xi_0} \varphi)^{\top} & R^{\top}\end{pmatrix}.
\label{eq:zeta_recovery_matrix}
\end{equation}
\end{enumerate}
If $s_i=0$ for some $i\in \overline{1,m}$ we only claim that \eqref{eq:additive_decomposition}--\eqref{eq:balanceLaw}, \eqref{eq:y_to_sigma_change} $\Longrightarrow$ 
\eqref{eq:sp}, \eqref{eq:K,moving_set}.
\end{Theorem}

\begin{Remark} Formula \eqref{eq:sp} implies that the initial state is admissible, i.e. 
\begin{equation}
\begin{pmatrix}y(0)\\ a(0)\end{pmatrix} \in \mathcal{C}(0, a(0)).
\label{eq:sp_valid_initial_condition}
\end{equation}
\end{Remark}

\begin{Remark}
\label{remark:C-inequalities-in-different-forms}
Moving set $\mathcal{C}(t,a)$ in \eqref{eq:K,moving_set} can be written in the several following ways. At first, observe that
\begin{equation}
\mathcal{C}(t,\widetilde{a}) =\left \{\begin{pmatrix} y\\a \end{pmatrix}\in \mathcal{V}\times \mathbb{R}^m:
\begin{array}{r}-(y_i-(Gr(t)-Ff(t))_i)- s_i(a_i-h_i\widetilde{a}_i)\leqslant c_{0,i},\\(y_i-(Gr(t)-Ff(t))_i)- s_i(a_i-h_i\widetilde{a}_i)\leqslant c_{0,i}\phantom{,}\end{array} \text{ for all }i\in \overline{1,m}
\right\} 
\label{eq:C_representation1}
\end{equation}
which is the same as 

\begin{equation}
\mathcal{C}(t,\widetilde{a})= \left \{\begin{pmatrix} y\\a \end{pmatrix}\in \mathcal{V}\times \mathbb{R}^m:
\begin{array}{r} \begin{pmatrix}-k_i e_i\\ - s_i e_i\end{pmatrix}^{\top} \mathbb{K}^{-1}\begin{pmatrix}y- (Gr(t)-Ff(t))\\ a - H \widetilde{a} \end{pmatrix}\leqslant c_{0,i},\\ \begin{pmatrix}k_i e_i\\ - s_i e_i\end{pmatrix}^{\top} \mathbb{K}^{-1} \begin{pmatrix}y-(Gr(t)-Ff(t))\\ a - H \widetilde{a} \end{pmatrix}\leqslant c_{0,i}\phantom{,}\end{array} \text{ for all }i\in \overline{1,m}
\right\}
\label{eq:moving_set_via_simple_normals}
\end{equation}
where by $e_i$ we denote the $i$-th vector of the standard basis in $\mathbb{R}^m$.

Furthermore, we would like to have the definition of $\mathcal{C}(t, \widetilde{a})$ with the normal vectors from within linear space $\mathcal{V}\times\mathbb{R}^m$. For that we can utilize the the variational definition of orthogonal projection \cite[Corollary~5.4, p.~134]{Brezis2011} and replace $k_ie_i$ with $k_iVP_Ve_i$, however due to $Ff(t)\notin \mathcal{V}$ we can only write 
\begin{equation}
\mathcal{C}(t,\widetilde{a})=\left \{\begin{pmatrix} y\\a \end{pmatrix}\in \mathcal{V}\times \mathbb{R}^m:
\begin{array}{r}n_{-i}^{\top} \,\mathbb{K}^{-1}\begin{pmatrix}y-Gr(t)\\ a - H \widetilde{a} \end{pmatrix}-(Ff(t))_i\leqslant c_{0,i},\\ n_{+i}^{\top}\, \mathbb{K}^{-1} \begin{pmatrix}y-Gr(t)\\ a - H \widetilde{a} \end{pmatrix}+(Ff(t))_i\leqslant c_{0,i}\phantom{,}\end{array} \text{ for all }i\in \overline{1,m}
\right\}
\label{eq:moving_set_via_normals_from_V}
\end{equation}
where 
\begin{equation}
n_{-i} =\begin{pmatrix}-k_i V P_V\,e_i\\ -s_i e_i\end{pmatrix}\in \mathcal{V}\times \mathbb{R}^m,\qquad n_{+i}=\begin{pmatrix}k_i V P_V \,e_i\\ -s_i e_i\end{pmatrix}\in \mathcal{V}\times \mathbb{R}^m.
\label{eq:widehat_n_definition}
\end{equation}
Thus we can write another expression of the set as
\begin{equation}
\mathcal{C}(t,\widetilde{a})=\mathcal{C}^f(t)+\begin{pmatrix}
Gr(t)\\ H\widetilde{a}
\end{pmatrix}
\label{eq:C-translation_part}
\end{equation}
where
\begin{equation}
 \mathcal{C}^f(t)=\left \{\begin{pmatrix} y\\a \end{pmatrix}\in \mathcal{V}\times \mathbb{R}^m:
\begin{array}{r}n_{-i}^{\top} \,\mathbb{K}^{-1}\begin{pmatrix}y\\ a\end{pmatrix}\leqslant c_{0,i}+(Ff(t))_i,\\ n_{+i}^{\top}\, \mathbb{K}^{-1} \begin{pmatrix}y\\ a \end{pmatrix}\leqslant c_{0,i}-(Ff(t))_i\phantom{,}\end{array} \text{ for all }i\in \overline{1,m}
\right\},
\label{eq:C^f_def}
\end{equation}
which highlights that the damage state $\widetilde{a}$ and the displacement loading $Gr(t)$ affect $C(t, \widetilde{a})$ as parallel translations. 

Finally, for our convenience, we introduce the notation 
\[
\mathcal{C}(t,\widetilde{a})=\left \{\begin{pmatrix} y\\a \end{pmatrix}\in \mathcal{V}\times \mathbb{R}^m:
\begin{array}{r}g_{-i}\left(t, \widetilde{a}, \begin{pmatrix}y\\ a\end{pmatrix}\right)\leqslant 0,\\ g_{+i}\left(t, \widetilde{a}, \begin{pmatrix}y\\ a\end{pmatrix}\right)\leqslant 0\phantom{,}\end{array} \text{ for all }i\in \overline{1,m}
\right\},
\]
\begin{multline}
 g_{-i}\left(t, \widetilde{a}, \begin{pmatrix}y\\ a\end{pmatrix}\right) =-(y_i-(Gr(t)-Ff(t))_i)- s_i(a_i-h_i\widetilde{a}_i)-c_{0,i}=\\= n_{-i}^{\top} \,\mathbb{K}^{-1}\begin{pmatrix}y-Gr(t)\\ a - H \widetilde{a} \end{pmatrix}-(Ff(t))_i -c_{0,i},
\label{eq:g_minus_func_def}
\end{multline}
\begin{multline}
 g_{+i}\left(t, \widetilde{a}, \begin{pmatrix}y\\ a\end{pmatrix}\right) =(y_i-(Gr(t)-Ff(t))_i)-s_i(a_i-h_i\widetilde{a}_i)-c_{0,i}=\\
 =n_{+i}^{\top}\, \mathbb{K}^{-1} \begin{pmatrix}y-Gr(t)\\ a - H \widetilde{a} \end{pmatrix}+(Ff(t))_i -c_{0,i}.
\label{eq:g_plus_func_def}
\end{multline}
\end{Remark}
The statement and the proof of Theorem \ref{th:LSMtoMSP} are similar to \cite[Th. 5.1]{Gudoshnikov2021ESAIM} and \cite[Th. 5.1]{Gudoshnikov2023preprint}.
There we considered lattices with elastic-perfectly plastic springs and proved the implication
similar to {\it \ref{enum:modeling_th_i}}$\Rightarrow${\it \ref{enum:modeling_th_ii}}. However, current Theorem \ref{th:LSMtoMSP} covers all the cases of hardening, perfect
plasticity and softening in a uniform manner, and we are able to prove the two-sided statement  {\it \ref{enum:modeling_th_i}}$\Leftrightarrow${\it \ref{enum:modeling_th_ii}}.

\noindent{\bf Proof of Theorem \ref{th:LSMtoMSP}.} {\bf From the LSM {\it \ref{enum:modeling_th_i}} to the Sweeping Process {\it \ref{enum:modeling_th_ii}}.}
To show \eqref{eq:sp}, \eqref{eq:K,moving_set} we only need $s_i\geqslant 0$ for all $i\in \overline{1,m}$. We start by transforming the relations \eqref{eq:gc}--\eqref{eq:bc} and \eqref{eq:balanceLaw} into inclusions in terms of spaces $\mathcal{U}$ and $\mathcal{V}$, respectively.

We fix a.e. $t\in[0,T]$, take time-detivative of \eqref{eq:bc} and apply the Moore-Penrose pseudoinverse matrix $R^+$ of $R$ to it:
\[
R^+ R \dot \zeta + R^+\dot r(t)=0.
\]
Matrix $R^+R$ is the orthogonal projection matrix onto ${\rm Im}\, R^{\top}$ (\cite[Def. 1.1.2]{Campbell2008}), therefore there is $z\in ({\rm Im}\, R^{\top})^\perp={\rm Ker}\, R$, s.t. $\dot \zeta = R^+ R \dot \zeta + z$ and 
\[
\dot \zeta -z +R^+\dot r(t)=0.
\]
Apply $K(D_{\xi_0} \varphi)$:
\[
K(D_{\xi_0} \varphi)\dot \zeta - K(D_{\xi_0} \varphi) z + K(D_{\xi_0} \varphi)R^+\dot r(t)=0,
\]
and  from \eqref{eq:gc} and \eqref{eq:U_def} we get
\begin{equation}
K\dot x \in \mathcal{U} -  K(D_{\xi_0} \varphi)R^+\dot r(t).
\label{eq:Kx-in-something}
\end{equation}
The last term can be represented as the sum of projections:
\[
K\dot x \in \mathcal{U} -  (UP_U+VP_V)K(D_{\xi_0} \varphi)R^+\dot r(t).
\]
hence
\[
K\dot x \in \mathcal{U} -  VP_VK(D_{\xi_0} \varphi)R^+\dot r(t),
\]
i.e.
\begin{equation}
K\dot x \in \mathcal{U} - G\dot r(t).
\label{eq:x-in-U-G}
\end{equation}

On the other hand, as we observed in Proposition \ref{prop:U_and_V_orth} \ref{enum:prop_U_and_V_equivalent_to_balance}, equation \eqref{eq:balanceLaw} means precisely that
\begin{equation}
-\sigma + F_1 f(t) \in \mathcal{V}
\label{eq:from_balance_law1}
\end{equation}
i.e.
\begin{equation}
-y+Gr(t) - (UP_U+VP_V)Ff(t) + F_1 f(t) \in \mathcal{V}
\label{eq:from_balance_law2}
\end{equation}
and by collecting all the terms from $\mathcal{V}$ and using \eqref{eq:F_def} we get that
\[
-y\in \mathcal{V},
\]
\begin{equation}
y\in \mathcal{V}.
\label{eq:y_in_V}
\end{equation}
At last, from \eqref{eq:additive_decomposition} and \eqref{eq:Hooke's_law} we  see that
\[
Kp=K x - \sigma.
\]
Then, using \eqref{eq:nc_basic_property2},  we can write \eqref{eq:plastic_flow} as 
\[
\frac{d}{dt}\begin{pmatrix}Kx-\sigma\\-a\end{pmatrix}\in N^{\mathbb{K}^{-1}}_{{\rm C}+\begin{pmatrix}0\\ Ha\end{pmatrix}}\begin{pmatrix}\sigma \\ a\end{pmatrix},
\]
\[
\frac{d}{dt}\begin{pmatrix}-\sigma\\-a\end{pmatrix}\in N^{\mathbb{K}^{-1}}_{{\rm C}+\begin{pmatrix}0\\ Ha\end{pmatrix}}\begin{pmatrix}\sigma \\ a\end{pmatrix}-\begin{pmatrix}K\dot x\\0 \end{pmatrix}.
\]
Plug \eqref{eq:x-in-U-G} and \eqref{eq:sigma_to_y_change} :
\[
\frac{d}{dt}\begin{pmatrix}-y\\-a\end{pmatrix}+ \begin{pmatrix}G\dot r(t)\\0 \end{pmatrix}- \begin{pmatrix}F\dot f(t)\\0 \end{pmatrix}\in N^{\mathbb{K}^{-1}}_{{\rm C}+\begin{pmatrix}0\\ Ha\end{pmatrix}}\begin{pmatrix}y-Gr(t)+Ff(t) \\ a\end{pmatrix}-\mathcal{U}\times \{0_m\}+ \begin{pmatrix}G\dot r(t)\\0 \end{pmatrix}.
\]
Cancel the term with $G\dot r(t)$ and use the fact that $F\dot f(t)\in \mathcal{U}$:
\begin{equation}
-\frac{d}{dt}\begin{pmatrix}y\\a\end{pmatrix} \in N^{\mathbb{K}^{-1}}_{{\rm C}+\begin{pmatrix}0\\ Ha\end{pmatrix}}\begin{pmatrix}y- Gr(t)+Ff(t) \\ a\end{pmatrix}-\mathcal{U}\times \{0_m\}.
\label{eq:sweeping_process_assembly1}
\end{equation}
Due to \eqref{eq:y_in_V} and Proposition \ref{prop:U_and_V_orth} \ref{enum:prop_U_and_V_orth} we have 
\[
-\mathcal{U}\times \{0_m\}= \mathcal{U}\times \{0_m\} = N^{\mathbb{K}^{-1}}_{\mathcal{V}\times \mathbb{R}^m}\begin{pmatrix}y\\a\end{pmatrix},
\]
therefore, also by using \eqref{eq:nc_basic_property1}, we get
\begin{equation}
-\frac{d}{dt}\begin{pmatrix}y\\a\end{pmatrix} \in N^{\mathbb{K}^{-1}}_{{\rm C}+\begin{pmatrix} Gr(t)-Ff(t)\\ Ha\end{pmatrix}}\begin{pmatrix}y \\ a\end{pmatrix}+N^{\mathbb{K}^{-1}}_{\mathcal{V}\times \mathbb{R}^m}\begin{pmatrix}y\\a\end{pmatrix}
\label{eq:sweeping_process_assembly2}
\end{equation}
Employing the additivity of the normal cones to polyhedral sets \cite[Corollary 23.8.1]{Rockafellar1970}, we obtain
\[
-\frac{d}{dt}\begin{pmatrix}y\\a\end{pmatrix} \in N^{\mathbb{K}^{-1}}_{\left({\rm C}+\begin{pmatrix} Gr(t)-Ff(t)\\ Ha\end{pmatrix}\right)\cap\left(\mathcal{V}\times\mathbb{R}^m\right)}\begin{pmatrix}y \\ a\end{pmatrix},
\]
which is exactly \eqref{eq:sp}--\eqref{eq:K,moving_set}.

Equations \eqref{eq:sigma_to_y_change}, \eqref{eq:sigma_to_eps} and \eqref{eq:additive_decomposition} follow trivially from {\it \ref{enum:modeling_th_i}}.

Now we assume that $s_i>0$ 
for all $i\in \overline{1,m}$ and we proceed to show \eqref{eq:p_recovery_formula}--\eqref{eq:compatible_initial_conditions}. 
Notice from \eqref{eq:plastic_flow} and \eqref{eq:plastic_flow_normal_cone} that for all $i\in \overline{1,m}$ a.a. $t\in[0,T]$ 
\[
\frac{d}{dt}\begin{pmatrix}p_i\\ -a_i\end{pmatrix}= \lambda_i^{-}\begin{pmatrix} -1\\ -s_i \end{pmatrix} +\lambda_i^{+}\begin{pmatrix} 1\\ -s_i \end{pmatrix} \qquad \text{ for some $t$-dependent } \lambda_i^{-}, \lambda_i^{+}\geqslant 0.
\]
Observe that $\lambda_i^{-}>0$ and $\lambda_i^{+}>0$ simultaneously is possible only when both equalities of \eqref{eq:failure_equations1} hold true, and only when $h_i>1$ (see Fig. \ref{fig:plasticity_types} f). Thus it is the state of complete failure, i.e. we can use the implication \eqref{eq:failure_means_no_change} to deduce that
\[
\lambda_i^{-}=\lambda_i^{+}=0=\dot p_i.
\]
Otherwise
\[
\lambda_i=\frac{\dot a_i}{s_i}, \qquad \dot p_i = \mp \lambda_i \qquad \text{ for some } \lambda_i\geqslant 0.
\]
More specifically,
\begin{equation}
\dot p_i = \left\{ \begin{array}{rcr}- \frac{\dot a_i}{s_i} & \text{ if }& -\sigma_i -s_i(1-h_i)a_i=c_{0,i},\\
\frac{\dot a_i}{s_i} & \text{ if }& \sigma_i-s_i(1-h_i)a_i=c_{0,i},\\
0 &&\text{otherwise.} \end{array}\right.= \left\{ \begin{array}{rcr}- \frac{\dot a_i}{s_i} & \text{ if }& g_{-i}\left(t, a, \begin{pmatrix}y\\a\end{pmatrix}\right)=0,\\
\frac{\dot a_i}{s_i} & \text{ if }& g_{+i}\left(t, a, \begin{pmatrix}y\\a\end{pmatrix}\right)=0,\\
0 &&\text{otherwise.} \end{array}\right.
\label{eq:equivalent_plastic_formula_1}
\end{equation}
Equivalently, we can  write
\begin{equation}
\dot p_i = \frac{\dot a_i}{s_i} \,{\rm sign}(\sigma_i) = \frac{\dot a_i}{s_i} \,{\rm sign}\left(y_i- (Gr(t)-Ff(t))_i\right),
\label{eq:equivalent_plastic_formula_2}
\end{equation}
where the sign function is used to determine which of the equalities in \eqref{eq:failure_equations1} holds (if neither then $\dot a_i=\dot p_i=0$ anyways). Expression \eqref{eq:p_recovery_formula} follows.

Finally, from \eqref{eq:gc}--\eqref{eq:bc} we have 
\begin{equation}
\begin{pmatrix} D_{\xi_0} \varphi \\ R\end{pmatrix} \zeta = \begin{pmatrix}x\\ -R\xi_0-r(t)\end{pmatrix},
\label{eq:gc_bc_on_zeta}
\end{equation}
where $(m+q)\times (nd)$ enhanced compatibility matrix $\begin{pmatrix} D_{\xi_0} \varphi \\ R\end{pmatrix}$ has full column rank due to assumption \eqref{eq:trivial_kernel_intersect}. Thus its Moore-Penrose pseudoinverse matrix is a left inverse and can be computed as \eqref{eq:zeta_recovery_matrix}, see \cite[Prop. 2.3]{Gudoshnikov2023preprint}. Expression \eqref{eq:zeta_recovery_formula} for $\zeta$ then follows immediately, and \eqref{eq:compatible_initial_conditions} is trivially implied by \eqref{eq:gc_bc_on_zeta} as well.

{\bf From the Sweeping Process {\it \ref{enum:modeling_th_ii}} to the LSM {\it \ref{enum:modeling_th_i}}.} Let  $s_i>0$ for all $i \in \overline{1,m}$ and assume that initial data satisfy \eqref{eq:compatible_initial_conditions}. 
Let $\begin{pmatrix}y\\ a\end{pmatrix}$ be a solution of \eqref{eq:sp}--\eqref{eq:K,moving_set} on $[0,T]$ and compute $\sigma, \varepsilon, p, x, \zeta$ by \eqref{eq:sigma_to_y_change}--\eqref{eq:p_recovery_formula}, \eqref{eq:additive_decomposition},\eqref{eq:zeta_recovery_formula}. We immediately have \eqref{eq:additive_decomposition} and \eqref{eq:Hooke's_law}. From \eqref{eq:sp}--\eqref{eq:K,moving_set} we necessarily have \eqref{eq:y_in_V}, which can be transformed back to \eqref{eq:from_balance_law2}, then to \eqref{eq:from_balance_law1}, which is, again, equivalent to 
\eqref{eq:balanceLaw}.

Furthermore, \eqref{eq:equivalent_plastic_formula_2} follows from  \eqref{eq:p_recovery_formula}. Similarly to Remark \ref{remark:a_nondecreasing1}, we see that sweeping process \eqref{eq:sp} implies $\dot a_i\geqslant 0$. As above, we do not consider the state of complete failure, and from \eqref{eq:equivalent_plastic_formula_2} we deduce \eqref{eq:equivalent_plastic_formula_1}, which, in turn, means that \eqref{eq:plastic_flow} is satisfied with \eqref{eq:plastic_flow_set}--\eqref{eq:plastic_flow_normal_cone}.

We need to show that  \eqref{eq:gc} and \eqref{eq:bc} hold. From \eqref{eq:sp} we can go back to \eqref{eq:sweeping_process_assembly2}, and then to \eqref{eq:sweeping_process_assembly1}, in which we plug \eqref{eq:y_to_sigma_change} to get
\[
-\frac{d}{dt}\begin{pmatrix}\sigma+Gr(t)-Ff(t)\\a\end{pmatrix} \in N^{\mathbb{K}^{-1}}_{{\rm C}+\begin{pmatrix}0\\ Ha\end{pmatrix}}\begin{pmatrix}\sigma \\ a\end{pmatrix}-\mathcal{U}\times \{0_m\}
\]
Since we have \eqref{eq:additive_decomposition} and \eqref{eq:Hooke's_law}, we can write
\[
\sigma= Kx-Kp.
\]
Hence, knowing that $F \dot f(t)\in \mathcal{U}$,
\[
-\frac{d}{dt}\begin{pmatrix}Kx+Gr(t)\\0 \end{pmatrix}+\frac{d}{dt}\begin{pmatrix}Kp\\-a\end{pmatrix} \in N^{\mathbb{K}^{-1}}_{{\rm C}+\begin{pmatrix}0\\ Ha\end{pmatrix}}\begin{pmatrix}\sigma \\ a\end{pmatrix}-\mathcal{U}\times \{0_m\}
\]
We  already have proven \eqref{eq:plastic_flow}, which together with \eqref{eq:nc_basic_property2} allows us to write
\[
-\frac{d}{dt}\begin{pmatrix}Kx+Gr(t)\\0 \end{pmatrix}\in N^{\mathbb{K}^{-1}}_{{\rm C}+\begin{pmatrix}0\\ Ha\end{pmatrix}}\begin{pmatrix}\sigma \\ a\end{pmatrix}-\mathcal{U}\times \{0_m\},
\]
i.e. there are two components $\begin{pmatrix}w_1\\ \beta_1\end{pmatrix}, \begin{pmatrix}w_2\\ 0\end{pmatrix}\in \mathbb{R}^{2m}$ such that

\begin{equation}
\begin{pmatrix}w_1\\ \beta_1\end{pmatrix}+\begin{pmatrix}w_2\\ 0\end{pmatrix}=-\frac{d}{dt}\begin{pmatrix}Kx+Gr(t)\\0 \end{pmatrix}, \qquad \begin{pmatrix}w_1\\ \beta_1\end{pmatrix}\in N^{\mathbb{K}^{-1}}_{{\rm C}+\begin{pmatrix}0\\ Ha\end{pmatrix}}\begin{pmatrix}\sigma \\ a\end{pmatrix},\qquad  \begin{pmatrix}w_2\\ 0\end{pmatrix}\in \mathcal{U}\times \{0_m\}.
\label{eq:proof_critical_decomposition}
\end{equation}
\vspace{1mm}

\noindent Hence $\beta_1=0$. But, notice from \eqref{eq:plastic_flow_normal_cone}, that the only element of the normal cone $N_{{\rm C}+\begin{pmatrix}0\\ Ha\end{pmatrix}}\begin{pmatrix}\sigma \\ a\end{pmatrix}$ which can have zero $a$-complonent is $0_{2m}$. By \eqref{eq:nc_basic_property2} the same can be said about the normal cone in \eqref{eq:proof_critical_decomposition}, and we have $\begin{pmatrix}w_1\\ \beta_1\end{pmatrix} =0_{2m}$, i.e.
\[
-\frac{d}{dt}\begin{pmatrix}Kx+Gr(t)\\0 \end{pmatrix} = \begin{pmatrix}w_2\\ 0\end{pmatrix}\in \mathcal{U}\times \{0_m\},
\]
i.e. we have proven \eqref{eq:x-in-U-G}. It is equivalent to \eqref{eq:Kx-in-something} and to 
\[
\dot x\in (D_{\xi_0} \varphi)({\rm Ker}\, R- R^+\dot r(t)),
\]
which, in turn, means that there is $z\in{\rm Ker}\, R$ such that
\begin{equation}
\dot x = (D_{\xi_0} \varphi)(z- R^+\dot r(t)).
\label{eq:proof_linear_algebra_1}
\end{equation}
Now observe that
\begin{equation}
R(z- R^+\dot r(t)) = Rz - RR^+\dot r(t) = -\dot r(t).
\label{eq:proof_linear_algebra_2}
\end{equation}
where the last equality is due to the fact that $R^+$ is the right inverse of $R$ because of \eqref{eq:rowsR_indep} and \cite[Prop. 2.3]{Gudoshnikov2023preprint}. Write \eqref{eq:proof_linear_algebra_1} and \eqref{eq:proof_linear_algebra_2} together
\[
\begin{pmatrix}
D_{\xi_0} \varphi\\
R
\end{pmatrix} \left(z- R^+\dot r(t)\right) = \begin{pmatrix}\dot x \\ -\dot r(t) \end{pmatrix},
\]
to observe that
\begin{equation}
\begin{pmatrix}\dot x \\ -\dot r(t) \end{pmatrix} \in {\rm Im}\, \begin{pmatrix}D_{\xi_0} \varphi\\
R\end{pmatrix}.
\label{eq:proof_linear_algebra_3}
\end{equation}
On the other hand, take the time-derivative of \eqref{eq:zeta_recovery_formula} and apply $\begin{pmatrix}D_{\xi_0} \varphi\\
R\end{pmatrix}$:
\[
\begin{pmatrix}D_{\xi_0} \varphi\\
R\end{pmatrix}\dot \zeta = \begin{pmatrix}D_{\xi_0} \varphi\\
R\end{pmatrix}\begin{pmatrix}D_{\xi_0} \varphi\\
R\end{pmatrix}^+\begin{pmatrix}\dot x \\ -\dot r(t) \end{pmatrix}.
\]
Matrix $\begin{pmatrix}D_{\xi_0} \varphi\\
R\end{pmatrix}\begin{pmatrix}D_{\xi_0} \varphi\\
R\end{pmatrix}^+$ is the orthogonal projection matrix onto ${\rm Im}\,\begin{pmatrix}D_{\xi_0} \varphi \\
R\end{pmatrix}$ \cite[Def. 1.1.2]{Campbell2008}, but, since we already know \eqref{eq:proof_linear_algebra_3}, it acts on $\begin{pmatrix}\dot x \\ -\dot r(t) \end{pmatrix}$ as the identity matrix and we conclude that
\[
\begin{pmatrix}D_{\xi_0} \varphi\\
R\end{pmatrix}\dot \zeta = \begin{pmatrix}\dot x \\ -\dot r(t) \end{pmatrix},
\]
which, together with \eqref{eq:compatible_initial_conditions} implies \eqref{eq:gc}, \eqref{eq:bc}.
$\blacksquare$

\subsection{Implicit catch-up algorithm}
One of the common methods to solve a state-dependent sweeping process \eqref{eq:sp} is the implicit catch-up scheme, see e.g. \cite{Kunze1998}, \cite[Section 3.3]{Kunze2000} \cite[Th. 3.3]{CastaingIbrahimYarou2009}. Specifically, we partition the time interval as
\begin{equation}
0=t_0<t_1<\dots<t_{k-1}< t_k=T
\label{eq:time-discretization-partition}
\end{equation}

\noindent and then, for the sweeping process \eqref{eq:sp}--\eqref{eq:K,moving_set}, on $j$-th step we are looking for a fixed point  
of the projection map 
\begin{equation}
\begin{pmatrix} \widetilde{y}\\ \widetilde{a}\end{pmatrix} \mapsto {\rm proj}^{\mathbb{K}^{-1}}\left(\begin{pmatrix}y^{j-1}\\a^{j-1}\end{pmatrix}, \mathcal{C}(t_j, \widetilde{a})\right),
\label{eq:full_iteration_map}
\end{equation}
and put such a fixed point as $\begin{pmatrix}y^{j}\\a^{j}\end{pmatrix}$ of the discretized solution. 
Since $\mathcal{C}(t, \widetilde{a})$ does not depend on $\widetilde{y}$, it is enough to find a fixed point of the map 
\begin{equation}
\begin{array}{rccl}
T_j:&\mathbb{R}^m &\to& \mathbb{R}^m\\
T_j:&\widetilde{a}&\mapsto&\begin{pmatrix}0& I\end{pmatrix}\,{\rm proj}^{\mathbb{K}^{-1}}\left(\begin{pmatrix}y^{j-1}\\a^{j-1}\end{pmatrix}, \mathcal{C}(t_j, \widetilde{a})\right).
\end{array}
\label{eq:T_i}
\end{equation}
The standard requirement for the existence for a fixed point is to have the moving set Lipschitz-continuous with Lipschitz constant corresponding to the state variable being strictly less then $1$. I. e. for our problem we write the conditions of Theorem \ref{th:sp_state-dep} and its counterpart  \cite[Lemma~2.1, p.~182]{Kunze1998} for a single time-step:
\begin{equation}
d_{H}^{\mathbb{K}^{-1}}(\mathcal{C}(t_1,\widetilde{a}_1),\mathcal{C}(t_2,\widetilde{a}_2))\leqslant L_1|t_1-t_2|+ L_2\|\widetilde{a}_1-\widetilde{a}_2\|, \qquad L_1, L_2> 0, 
\label{eq:C_Lipschitz_continuity}
\end{equation}
\begin{equation}
L_2<1,
\label{eq:Lip_const_less_than_1}
\end{equation}
where $d_{H}^{\mathbb{K}^{-1}}$ is the Hausdorff distance induced by the norm in the inner product space $\mathcal{V}\times \mathbb{R}^m$:
\begin{multline}
d_{H}^{\mathbb{K}^{-1}}(C_1, C_2)=\max\left(\sup_{\Sigma_2\in C_2}\, \inf_{\Sigma_1\in C_1}\|\Sigma_1-\Sigma_2\|_{\mathbb{K}^{-1}},\sup_{\Sigma_1\in C_1}\, \inf_{\Sigma_2\in C_2}\|\Sigma_1-\Sigma_2\|_{\mathbb{K}^{-1}} \right)\\
\text{for } C_1, C_2\subset \mathcal{V}\times\mathbb{R}^m,
\label{eq:Haus_dist_formula}
\end{multline}
\begin{equation}
\|\Sigma\|_{\mathbb{K}^{-1}}=\sqrt{\Sigma^{\top} \mathbb{K}^{-1}\Sigma}\qquad \text{for }\Sigma\in\mathcal{V}\times\mathbb{R}^m.
\label{eq:K_norm_formula}
\end{equation}

\begin{Remark}
While \eqref{eq:Lip_const_less_than_1} is sufficient for existence of a fixed point, observe from \eqref{eq:C-translation_part} that one cannot expect \eqref{eq:Lip_const_less_than_1}  to hold when $h_i\geqslant1$ for some $i\in \overline{1,m}$ (a perfectly plastic spring, or a spring  with softening): the points of $\mathcal{C}(t, a)$ which correspond to the state of complete failure of spring $i$ will lead to $L_2\geqslant1$ (see e.g. the ``tip'' of the set in Fig. \ref{fig:lin_soft_yielding_and_sweeping} b). 
\label{rem:th-only-hardening}
\end{Remark}

An even stricter (in our setting) assumption would be to require $T_j$ to be a contraction mapping, which would imply both existence and uniqueness of its fixed point. However, as we will demonstrate in Section \ref{sect:ex1},  map $T_j$ in case of softening  may not be a contraction even when the state of complete failure is not attained. Examples show that $T_j$ may have multiple fixed points, which means the co-existence of multiple branching solutions of the state-dependent sweeping process \eqref{eq:sp}--\eqref{eq:K,moving_set}.

\begin{Remark}
\label{remark:initial_guess}
To find numerically the fixed points of the map $T_j$ at each time-step we iterate the map and hope for convergence to a fixed point, see Algorithm \ref{alg:implicit-catch-up1} below. However, this means that for each time-step we must make a choice of an initial value for the iterations, and such choice is arbitrary with respect to the given problem data. In Algorithm \ref{alg:implicit-catch-up1} we denote the choice as a given function ${\rm InitialValue}\left(j, a^{j-1}\right)$, and the simplest way to choose is to take the value from the previous time-step, i.e. to set
\[
{\rm InitialValue}\left(j, a^{j-1}\right) :=  a^{j-1}.
\]
However, as we will illustrate in the next sections, sometimes the map $T_j$ has multiple fixed points simultaneously, and different prescribed choices of ${\rm InitialValue}$ will result in convergence to different fixed points. This gives the algorithm an advantage over the explicit catch-up method, which also converges \cite{Haddad2013} under the assumptions, similar to \eqref{eq:C_Lipschitz_continuity}--\eqref{eq:Lip_const_less_than_1},
 but allows no control over the choice of a solution.

But, strictly speaking, Algorithm \ref{alg:implicit-catch-up1}  is searching for {\it asymptotically stable} fixed points and will fail (diverge) in a situation when $T_j$ has only unstable fixed points, as well as when there are none at all.
\end{Remark}
\begin{algorithm}[H]
\SetAlgoLined
\SetKwComment{Comment}{//}{}
\Comment{Given: $\mathbb{K}$, $\mathcal{C}(t, a)$ by \eqref{eq:K,moving_set}, $\begin{pmatrix}y^0\\a^0\end{pmatrix}\in \mathbb{R}^{2m}$ admissible in the sense of \eqref{eq:sp_valid_initial_condition},}
\Comment{a partition $t_j, j\in \overline{0,k}$  as in \eqref{eq:time-discretization-partition},} 
\Comment{an initial guess function ${\rm InitialValue}\left(j, a^{j-1}\right)$,}
\Comment{and the constants for the stopping condition:}
\Comment{very small real $\epsilon>0$, large integer ${ I_{\max}}$, large real $R>0$.}
\For{$j:=1$ \KwTo $k$}{
$I:=0$\;
converged $:=$ {\bf false}\;
diverged $:=$ {\bf false}\;
\Comment{prescribed at this time-step initial guess for the iterations:}
$a_0^{j}:={\rm InitialValue}\left(j,a^{j-1}\right)$\;
\While{{\bf not }\rm (converged {\bf or} diverged)}{
$I:=I+1$\;
  \Comment{An iteration of the map \eqref{eq:full_iteration_map}:}
  $\begin{pmatrix}y_{I}^{j}\\[3pt] a_{I}^{j}\end{pmatrix} := {\rm proj}^{\mathbb{K}^{-1}}\left(\begin{pmatrix}y^{j-1}\\a^{j-1}\end{pmatrix}, \mathcal{C}\left(t_{j}, a_{I-1}^{j}\right)\right)$\;
\If{$\left\| a_{I}^{j} - a_{I-1}^{j}\right\|<\epsilon$}
{
converged $:=$ {\bf true}\;
}
\If{$I> I_{\max}$ {\bf or} $\left\|\begin{pmatrix}y_{I}^{j}\\[3pt] a_{I}^{j}\end{pmatrix}\right\|>R$}{
\Comment{Cannot find an asymptotically stable fixed point of the map $T_j$ for $j$-th step, its iterations diverge.}
diverged $:=$ {\bf true}\;
}
}
\If {\rm diverged}
{{\bf break}}
$\begin{pmatrix}y^{j}\\ a^{j}\end{pmatrix}:= \begin{pmatrix}y_{I}^{j}\\[3pt] a_{I}^{j}\end{pmatrix}$\;
}  
 \caption{Iterative numerical scheme to solve the state-dependent sweeping process \eqref{eq:sp}--\eqref{eq:K,moving_set}}
 \label{alg:implicit-catch-up1}
\end{algorithm}

\subsection{Equivalent sweeping process of reduced dimension}

One must note that set $\mathcal{C}(t,a)$ is described by time- and state-dependent inequality constraint (Remark \ref{remark:C-inequalities-in-different-forms}) and  a linear equality constraint of $\mathcal{V}\times\mathbb{R}^m\subset \mathbb{R}^{2m}$, which can be explicitly obtained from Proposition \ref{prop:U_and_V_orth} \ref{enum:prop_U_and_V_orth}. The latter is independent from the arguments of $\mathcal{C}(t, a)$, therefore one can pass to an equivalent state-dependent sweeping process in $\mathbb{R}^{{\rm dim \mathcal{V}}+m}$, similarly to the procedure of \cite[Sect. 5.2]{Gudoshnikov2023preprint}. Specifically, we can define set $\mathcal{C}_V(t,a)$ such that 
\[
V\mathcal{C}_V(t,a)=\mathcal{C}(t,a)
\]
by substituting into \eqref{eq:moving_set_via_normals_from_V}--\eqref{eq:widehat_n_definition} the new unknown variable $\widehat{y}(t)\in\mathbb{R}^{{\rm dim}\,\mathcal{V}}$ such that
\begin{equation}
y=V\widehat{y}, \qquad \widehat{y} = P_V y.
\label{eq:change_of_variables_to_y}
\end{equation}
We obtain the new moving set
\begin{multline}
\mathcal{C}_V(t,\widetilde{a})=\\
=\left \{\begin{pmatrix} \widehat{y}\\a \end{pmatrix}\in\mathbb{R}^{{\rm dim}\, \mathcal{V} +m}:
\begin{array}{r} \widehat{n}_{-i}^{\top} \,S_V\begin{pmatrix}\widehat{y}-G_Vr(t)\\ a - H \widetilde{a} \end{pmatrix}-(Ff(t))_i\leqslant c_{0,i},\\ \widehat{n}_{+i}^{\top}\, S_V \begin{pmatrix}\widehat{y}-G_Vr(t)\\ a - H \widetilde{a} \end{pmatrix}+(Ff(t))_i\leqslant c_{0,i}\phantom{,}\end{array} \text{ for all }i\in \overline{1,m}
\right\}
\label{eq:sp_in_V_moving_set}
\end{multline}
defined via its normal vectors
\[
\widehat{n}_{-i} =\begin{pmatrix}-k_i P_V\,e_i\\ -s_i e_i\end{pmatrix}\in \mathbb{R}^{{\rm dim \mathcal{V}}+m},\qquad \widehat{n}_{+i}=\begin{pmatrix}k_i P_V \,e_i\\ -s_i e_i\end{pmatrix}\in\mathbb{R}^{{\rm dim \mathcal{V}}+m}\]
with respect to the inner product \eqref{eq:abstract_ip_via_a_matrix} in $\mathbb{R}^{{\rm dim}\, \mathcal{V} +m}$ defined via the matrix
\[
S_V:= \begin{pmatrix}V^{\top}K^{-1}V& 0\\ 0 & I\end{pmatrix}.
\]
We also denote the matrix 
\[
G_V:=P_V K (D_{\xi_0}\varphi)R^+,
\]
compare to \eqref{eq:G_def}.

Now observe that sweeping process \eqref{eq:sp} can be written as 
\[
\begin{pmatrix}\dot y\\ \dot a \end{pmatrix}^{\top} \mathbb{K}^{-1} \left(\begin{pmatrix}c_y \\c_a\end{pmatrix}-\begin{pmatrix} y\\  a \end{pmatrix}\right)\geqslant 0 \qquad \text{for any }\begin{pmatrix}c_y \\c_a\end{pmatrix}\in \mathcal{C}(t,a).
\]
Using change of variables \eqref{eq:change_of_variables_to_y} we obtain
\[
\left(\begin{pmatrix}V & 0 \\ 0 & I\end{pmatrix}\begin{pmatrix}\dot{\widehat{y}} \\ \dot a \end{pmatrix}\right)^{\top} \mathbb{K}^{-1}\begin{pmatrix}V & 0 \\ 0 & I\end{pmatrix} \left(\begin{pmatrix}c_{\widehat{y}} \\c_a\end{pmatrix}-\begin{pmatrix} \widehat{y} \\  a \end{pmatrix}\right)\geqslant 0 \qquad \text{for any }\begin{pmatrix}c_{\widehat{y}} \\c_a\end{pmatrix}\in \mathcal{C}_V(t,a),
\]
which is the sweeping process in $\mathbb{R}^{{\rm dim \mathcal{V}}+m}$:
\begin{equation}
-\frac{d}{dt}\begin{pmatrix} \widehat{y} \\a \end{pmatrix} \in N^{S_V}_{\mathcal{C}_V(t, a)}\begin{pmatrix} \widehat{y} \\a \end{pmatrix}.
\label{eq:sp_in_V}
\end{equation}
with initial condition
\[
\widehat{y}(0)= P_V y(0).
\]

Given a partition \eqref{eq:time-discretization-partition} we can numerically solve the sweeping process of reduced dimension \eqref{eq:sp_in_V_moving_set}--\eqref{eq:sp_in_V} by iterations of

\begin{equation}
\begin{pmatrix} \widetilde{\widehat{y}}\\ \widetilde{a}\end{pmatrix} \mapsto {\rm proj}^{S_V}\left(\begin{pmatrix}\widehat{y}^{j-1}\\a^{j-1}\end{pmatrix}, \mathcal{C}_V\left(t_j, \widetilde{a}\right)\right),
\label{eq:full_iteration_map_reduced_dim}
\end{equation}
and Algorithm \ref{alg:implicit-catch-up-reduced-dim} specifically. Map $T_j$ defined via projection \eqref{eq:full_iteration_map_reduced_dim} coincides with the map as defined by \eqref{eq:T_i}.

\begin{algorithm}[H]
\SetAlgoLined
\SetKwComment{Comment}{//}{}
\Comment{Given $S_V$, $\mathcal{C}_V(t, a)$, $\begin{pmatrix}\widehat{y}^0\\a^0\end{pmatrix}\in \mathbb{R}^{{\rm dim \mathcal{V}}+m}$ such that}
\Comment{$\begin{pmatrix}V\widehat{y}^0\\a^0\end{pmatrix}$ is admissible in the sense of \eqref{eq:sp_valid_initial_condition},}
\Comment{a partition $t_j, j\in \overline{0,k}$  as in \eqref{eq:time-discretization-partition},} 
\Comment{an initial guess function ${\rm InitialValue}\left(j, a^{j-1}\right)$,}
\Comment{and the constants for the stopping condition:}
\Comment{very small real $\epsilon>0$, large integer ${ I_{\max}}$, large real $R>0$.}
\For{$j:=1$ \KwTo $k$}{
${I}:=0$\;
converged $:=$ {\bf false}\;
diverged $:=$ {\bf false}\;
\Comment{initial guess for the iterations, prescribed at this time-step :}
$a_0^{j}:={\rm InitialValue}\left(j, a^{j-1}\right)$\;
\While{{\bf not }\rm (converged {\bf or} diverged)}{
$I:=I+1$\;
  \Comment{Iteration of map \eqref{eq:full_iteration_map_reduced_dim}:}
  $\begin{pmatrix}\widehat{y}_{I}^{j}\\[3pt] a_{I}^{j}\end{pmatrix} := {\rm proj}^{S_V}\left(\begin{pmatrix}\widehat{y}^{j-1}\\a^{j-1}\end{pmatrix}, \mathcal{C}_V\left(t_{j}, a_{I-1}^{j}\right)\right)$\;
\If{$\left\| a_{I}^{j} - a_{I-1}^{j}\right\|<\epsilon$}
{
converged $:=$ {\bf true}\;
}
\If{$I> I_{\max}$ {\bf or} $\left\|\begin{pmatrix}\widehat{y}_{I}^{j}\\[3pt] a_{I}^{j}\end{pmatrix}\right\|>R$}{
\Comment{Cannot find an asymptotically stable fixed point of the map $T_j$ for $j$-th step, its iterations diverge.}
diverged $:=$ {\bf true}\;
}
}
\If {\rm diverged}
{{\bf break}}
$\begin{pmatrix}\widehat{y}^{j}\\ a^{j}\end{pmatrix}:= \begin{pmatrix}\widehat{y}_{I}^{j}\\[3pt] a_{I}^{j}\end{pmatrix}$\;
}  
 \caption{Practical iterative numerical scheme to solve the state-dependent sweeping process of reduced dimension \eqref{eq:sp_in_V_moving_set}--\eqref{eq:sp_in_V}}
 \label{alg:implicit-catch-up-reduced-dim}
\end{algorithm}

For problems with large $m$ Algorithm \ref{alg:implicit-catch-up-reduced-dim} is orders of magnitude faster then Algorithm \ref{alg:implicit-catch-up1} in $\mathbb{R}^{2m}$ for sweeping process \eqref{eq:sp}, however the formulation of \eqref{eq:sp}, \eqref{eq:T_i} is sometimes easier to handle analytically.

\begin{Remark}
Observe that the mapping $(k, s, h)\mapsto (k, E^p)$, given by the formula for a plasticity modulus  \eqref{eq:modulus_formula_hs} is not injective, i.e. for the same values of $k$ and $E^p$ (therefore, for the same observed behavior of a spring), there are many corresponding pairs of values $(s,h)$, which result in such behavior. And, despite the observed behavior of the springs being the same, such different pairs of $(s,h)$ lead to significantly different number of iterations of Algorithm \ref{alg:implicit-catch-up-reduced-dim}. In particular, while computing the examples presented further in Section \ref{sect:LSM-examples} we noticed that equivalent parameters  with $s=1$ require much more iterations per time-step to converge.
\end{Remark}

\section{Two-springs system as a toy example with non-unique
solutions}
\label{sect:ex1}
As the first example of a lattice we consider a toy system of two springs connected in series along a line, see Fig. \ref{fig:intro_fig_models} b). Despite its simplicity, this example can help to illustrate our approach, important properties of the LSMs with softening, behavior of the implicit catch-up algorithm  and state-dependent sweeping processes in general.

\subsection{Given data}
The geometry of the system is represented by the following quantities used in \eqref{eq:gc}--\eqref{eq:bc}:
\begin{align*}
\text{number of spatial dimensions:}&&d&=1,&\\
\text{number of nodes:}&&n&=3,&\\
\text{number of springs:}&&m&=2,&\\
\text{reference configuration:}&&\xi_0&=  \begin{pmatrix} 0 & 1 & 2 \end{pmatrix}^{\top},&\\
\text{incidence matrix:}&&Q&=  \begin{pmatrix} 1 & 0 \\ -1 &  1\\ 0 & -1 \end{pmatrix}^{\top},&\\
\text{number of external displacement constraints:}&&q&=2,&\\
&&R&=\begin{pmatrix} 1 & 0 & 0 \\ 0 & 0 & 1\end{pmatrix},&\\
\text{displacement load:}&&r(t)& =  \begin{pmatrix} 0 & -l(t)  \end{pmatrix}^{\top},&\\
\text{stress load:} &&f(t)& =  \begin{pmatrix} f_1(t) & f_2(t) & f_3(t) \end{pmatrix}^{\top}.&
\end{align*}
However, we will see,  that due to the constraint \eqref{eq:bc}, only the component $f_2(t)$ of $f(t)$ influences the system, see Fig. \ref{fig:intro_fig_models} b).

As we have $d=1$, choose the orientation of springs to agree with the positive direction. Then \eqref{eq:D_xi_phi} yields
\begin{align*}
\text{compatibility matrix:}&&D_{\xi_0}\varphi&=-Q^{\top}= \begin{pmatrix} -1 & 1 & 0\\ 0 & -1 & 1 \end{pmatrix},&\\
\text{enhanced compatibility matrix:}&&\begin{pmatrix} D_{\xi_0} \varphi \\  R\end{pmatrix}&= \begin{pmatrix} -1 & 1 & 0\\ 0 & -1 & 1 \\ 1 & 0 & 0 \\ 0 & 0 & 1\end{pmatrix}.&\\
\end{align*}
Assumptions \eqref{eq:rowsR_indep}--\eqref{eq:trivial_kernel_intersect} are satisfied trivially. We compute (see \eqref{eq:new_equilibrium_matrix_pseudoinverse})
\begin{align*}
F_1 = \frac{1}{4}\begin{pmatrix}-1 & 2 & 1\\ -1 & -2 & 1\end{pmatrix}.
\end{align*}
At this point we leave the parameters of the springs arbitrary, i.e. we suppose to be given
\begin{equation}
\begin{aligned}
\text{Hooke's coefficients:}\quad&k_1, k_2>0,\\
\text{coefficients of state-dependent feedback:}\quad&h_1, h_2\in \mathbb{R},\\
\text{geometric slopes:}\quad&s_1, s_2>0,\\
\text{initial yield stresses:}\quad&c_{0,1},c_{0,2}>0.
\end{aligned}
\label{eq:ex1-parameters}
\end{equation}
thus we have
\[
K=\begin{pmatrix}k_1 &0 \\0 & k_2\end{pmatrix}, \qquad H=\begin{pmatrix}h_1 &0 \\0 & h_2\end{pmatrix}.
\]
\subsection{The sweeping process}
Observe that
\[
{\rm Ker }\, R  = {\rm lin }\, \begin{pmatrix}0 & 1 & 0\end{pmatrix}^{\top},
\]
therefore from \eqref{eq:V_perp_standard} we have
\[
\mathcal{V}^{\perp} =  (D_{\xi_0}\varphi) \, {\rm Ker}\, R = {\rm lin} \begin{pmatrix} -1 & 1 & 0\\ 0 & -1 & 1 \end{pmatrix} \begin{pmatrix}0 \\ 1 \\ 0\end{pmatrix} = {\rm lin} \begin{pmatrix} 1 \\ -1\end{pmatrix}
\]
and from \eqref{eq:U_def} we have
\[
U=\begin{pmatrix} k_1 \\ -k_2\end{pmatrix}, \qquad \mathcal{U}  = {\rm Im}\, U.
\]
Moreover,
\[
\mathcal{V} = {\rm lin} \begin{pmatrix}1\\1\end{pmatrix},
\]
but for future convenience we would like to have matrix $V$ such that 
\[
V^{\top} K^{-1} V = I,
\]
where $I$ in our case is $I_{1\times1}=1$. Thus if we take 
\begin{equation}
V=\kappa\begin{pmatrix}1\\1\end{pmatrix}
\label{eq:ex_1_V_gamma}
\end{equation}
for some $\kappa>0$ and solve
\[
\kappa \begin{pmatrix}1 & 1\end{pmatrix}\begin{pmatrix}k_{1}^{-1} & 0 \\ 0 & k_2^{-1} \end{pmatrix} \begin{pmatrix}1\\1\end{pmatrix}\kappa = 1,
\]
\[
\kappa^{2}(k_{1}^{-1}+k_{2}^{-1})=1,
\]
\begin{equation}
\kappa=\frac{1}{\sqrt{k_{1}^{-1}+k_{2}^{-1}}},
\label{eq:ex_1_gamma_found}
\end{equation}
so we have the basis vector of $\mathcal{V}$ as \eqref{eq:ex_1_V_gamma}--\eqref{eq:ex_1_gamma_found}.

Furthermore, by \eqref{eq:projection_matrices}
\begin{align*}
P_U &= \left(U^{\top} K^{-1} U\right)^{-1}U^{\top} K^{-1}= \frac{1}{k_1+k_2}\begin{pmatrix}1 & -1\end{pmatrix},
\\
P_V& = \left(V^{\top} K^{-1} V\right)^{-1}V^{\top} K^{-1} = I\, V^{\top} K^{-1}
=\kappa \begin{pmatrix}1& 1\end{pmatrix}\begin{pmatrix}k_{1}^{-1} & 0 \\ 0 & k_2^{-1} \end{pmatrix} = \kappa\begin{pmatrix}k_{1}^{-1}& k_{2}^{-1}\end{pmatrix},
\end{align*}
and
\begin{align*}
UP_U &= \frac{1}{k_1+k_2}\begin{pmatrix}k_1 &  -k_1 \\ -k_2 & k_2\end{pmatrix},
\\
VP_V &= \kappa^2\begin{pmatrix}1\\1\end{pmatrix}\begin{pmatrix}k_{1}^{-1}& k_{2}^{-1}\end{pmatrix}=\kappa^2\begin{pmatrix}k_{1}^{-1} & k_{2}^{-1}\\ k_{1}^{-1} & k_{2}^{-1}\end{pmatrix}.
\end{align*}
Finally,
\[
F = U P_U F_1 = \frac{1}{k_1+k_2}\begin{pmatrix}0 & k_1 & 0\\ 0 & -k_2 & 0\end{pmatrix}
\] 
\[
R^+ = R^{\top} \left(RR^{\top}\right)^{-1} = \begin{pmatrix}1 & 0\\ 0 & 0 \\ 0 & 1 \end{pmatrix} \begin{pmatrix}1 & 0\\ 0 & 1\end{pmatrix}^{-1}= \begin{pmatrix}1 & 0\\ 0 & 0 \\ 0 & 1 \end{pmatrix},
\]
\begin{multline*}
G_V \, r(t) = P_V K (D_{\xi_0}\varphi)R^+ r(t) =\\
= \kappa\begin{pmatrix}k_{1}^{-1}& k_{2}^{-1}\end{pmatrix}\begin{pmatrix}k_1 &0 \\0 & k_2\end{pmatrix} \begin{pmatrix} -1 & 1 & 0\\ 0 & -1 & 1 \end{pmatrix}\begin{pmatrix}1 & 0\\ 0 & 0 \\ 0 & 1 \end{pmatrix}\begin{pmatrix} 0 \\ -l(t)  \end{pmatrix}= \\
=\kappa\begin{pmatrix}1 & 1\end{pmatrix} \begin{pmatrix} -1 & 0 \\ 0 & 1\end{pmatrix}\begin{pmatrix} 0 \\ -l(t)  \end{pmatrix}=-\kappa l(t),
\end{multline*}
\begin{align*}
Gr(t)& = VG_V\, r(t)  =-\kappa^2 \begin{pmatrix}1\\1\end{pmatrix} l(t),\\
Ff(t)& = \begin{pmatrix}k_1\\-k_2\end{pmatrix}\frac{f_2(t)}{k_1+k_2}.
\end{align*}
Therefore we have sweeping process \eqref{eq:sp} in $\mathbb{R}^4$ where 
\[
\mathbb{K}^{-1} = \begin{pmatrix} k_{1}^{-1} & 0 & 0 & 0\\ 0 & k_{2}^{-1} & 0 & 0\\
0 & 0 & 1 & 0 \\ 0 & 0 & 0 & 1 \end{pmatrix},
\]
and, by \eqref{eq:moving_set_via_normals_from_V}--\eqref{eq:widehat_n_definition},
\[
n_{-1} = \begin{pmatrix}-\kappa^2\\ -\kappa^2\\ -s_1\\0\end{pmatrix}, \quad n_{-2}=\begin{pmatrix}-\kappa^2\\ -\kappa^2\\ 0\\-s_2\end{pmatrix}, \quad n_{+1} = \begin{pmatrix}\kappa^2\\ \kappa^2\\ -s_1\\0\end{pmatrix}, \quad n_{+2}=\begin{pmatrix}\kappa^2 \\ \kappa^2 \\ 0\\-s_2\end{pmatrix},
\]
\[
\mathcal{C}\left(t,\begin{pmatrix}\widetilde{a}_1\\ \widetilde{a}_2\end{pmatrix}\right)=\left\{ \begin{pmatrix}y_1\\ y_2\\a_1\\ a_2\end{pmatrix}\in \mathbb{R}^4: y_1=y_2,\, \begin{pmatrix}n_{-1}^{\top}\\[2mm] n_{-2}^{\top}\\[2mm] n_{+1}^{\top}\\[2mm] n_{+2}^{\top} \end{pmatrix}\mathbb{K}^{-1}\begin{pmatrix}y_1 +\kappa^2 l(t) \\[1mm] y_2 +\kappa^2 l(t)\\[1mm] a_1 - h_1\widetilde{a}_1\\[1mm] a_2- h_2\widetilde{a}_2\end{pmatrix} + \begin{pmatrix}-k_1\\
k_2\\
k_1\\
-k_2
\end{pmatrix}\frac{f_2(t)}{k_1+ k_2}\leqslant \begin{pmatrix}c_{0,1}\\[1mm] c_{0,2}\\[1mm] c_{0,1}\\[1mm] c_{0,2} \end{pmatrix} \right\},
\]
where the inequality is meant in the component-wise sense. 

Equivalent sweeping process \eqref{eq:sp_in_V} in $\mathbb{R}^3$ is defined via
\[
S_V = I_{3\times3}
\]
and
\begin{equation}
\widehat{n}_{-1} = \begin{pmatrix}-\kappa\\ -s_1\\0\end{pmatrix}, \, 
\widehat{n}_{-2}=\begin{pmatrix}-\kappa \\ 0\\-s_2\end{pmatrix}, \, 
\widehat{n}_{+1} = \begin{pmatrix}\kappa\\ -s_1\\0\end{pmatrix}, \, 
\widehat{n}_{+2}=\begin{pmatrix}\kappa\\ 0\\-s_2\end{pmatrix},
\label{eq:example1_normals_in_R3}
\end{equation}
\begin{equation}
\mathcal{C}_V\left(t,\begin{pmatrix}\widetilde{a}_1\\ \widetilde{a}_2\end{pmatrix}\right)=\left\{ \begin{pmatrix}\widehat{y}\\a_1\\ a_2\end{pmatrix}\in \mathbb{R}^3: \begin{pmatrix}\widehat{n}_{-1}^{\top}\\[2mm] \widehat{n}_{-2}^{\top}\\[2mm] \widehat{n}_{+1}^{\top}\\[2mm] \widehat{n}_{+2}^{\top} \end{pmatrix}\begin{pmatrix}\widehat{y} +\kappa\, l(t) \\[1mm] a_1 - h_1\widetilde{a}_1\\[1mm] a_2- h_2\widetilde{a}_2\end{pmatrix}+ \begin{pmatrix}-k_1\\
k_2\\
k_1\\
-k_2
\end{pmatrix}\frac{f_2(t)}{k_1+ k_2}\leqslant \begin{pmatrix}c_{0,1}\\[1mm] c_{0,2}\\[1mm] c_{0,1}\\[1mm] c_{0,2} \end{pmatrix} \right\}.
\label{eq:example1_C_in_R3}
\end{equation}
\subsection{Non-uniqueness of solutions and a bifurcation}
\label{ssect:ex1_solutions}
Consider the sweeping process \eqref{eq:sp_in_V}, \eqref{eq:example1_normals_in_R3}--\eqref{eq:example1_C_in_R3} on an interval $[t_0, T]$.  Let $l(t)$ be Lipschitz-continuous and monotonically increasing from zero, and $f$ be zero:
\begin{equation}
\dot l(t)>0, \quad f(t)\equiv 0 \quad \text{a.e. on }[t_0,T], \qquad l(t_0)=0.
\label{eq:ex1_ingcreasing_load}
\end{equation}
Choose an initial condition $\begin{pmatrix}\widehat{y}(t_0)\\ a(t_0)\end{pmatrix}\in \mathbb{R}^3$ on the boundary of $\mathcal{C}_V(t, a)$ corresponding to both normal vectors $\widehat{n}_{+1}, \widehat{n}_{+2}$ being active, i.e. in terms of \eqref{eq:g_plus_func_def} for $t=t_0$ we must have both
\begin{equation}
g_{+1}\left(t, a(t), \begin{pmatrix}V\widehat{y}(t)\\ a(t)\end{pmatrix}\right)=0, \qquad
g_{+2}\left(t, a(t), \begin{pmatrix}V\widehat{y}(t)\\ a(t)\end{pmatrix}\right)=0.
\label{eq:example1_ic_in_R3}
\end{equation}
Additionally, we require that the constraints with normal vectors $\widehat{n}_{-1}, \widehat{n}_{-2}$ are not active at the initial condition, i.e. for $t=t_0$:
\begin{equation}
g_{-1}\left(t, a(t), \begin{pmatrix}V\widehat{y}(t)\\ a(t)\end{pmatrix}\right)<0, \qquad
g_{-2}\left(t, a(t), \begin{pmatrix}V\widehat{y}(t)\\ a(t)\end{pmatrix}\right)<0
\label{eq:example1_ic_in_R3_not_broken}
\end{equation}
In other words we require the stresses of both springs to start at the upper threshold and {\it not} at the state of complete failure.

We consider three scenarios: {\it A.} both springs are hardening springs ($h_1, h_2<1$), {\it B.} both springs are perfectly plastic ($h_1=h_2=1$) and {\it C.} both springs are softening springs ($h_1, h_2>1$).
\begin{Proposition}
\label{prop:solutions_of_ex1}
Given the sweeping process \eqref{eq:sp_in_V}, \eqref{eq:example1_normals_in_R3}--\eqref{eq:example1_C_in_R3} on an interval $[t_0, T]$ and assuming \eqref{eq:ex1_ingcreasing_load}--\eqref{eq:example1_ic_in_R3_not_broken} the following solutions are present, as long as \eqref{eq:example1_ic_in_R3_not_broken} holds:
\begin{enumerate}[{\it A.}]
\item \label{enum:toy_example_hardening} If both $h_1, h_2<1$ then there is a single soluton with plastic deformation distributed between both springs, and it is given by
\begin{equation}
\begin{pmatrix}\dot{\widehat{y}}\\ \dot a_1\\ \dot a_2 \end{pmatrix}=-\frac{(1-h_2)s^2_2\,\widehat{n}_{+1}+ (1-h_1)s^2_1\,\widehat{n}_{+2}}{(1-h_1)s_1^2 + (1-h_2)s_2^2+ (k_{1}^{-1}+k_{2}^{-1})(1-h_1)(1-h_2)s_1^2s_2^2}\dot l(t)
\label{eq:example1_sol3}
\end{equation}
for a.a. $t\geqslant t_0$.
\item \label{enum:toy_example_pp} If $h_1=h_2=1$ then there is a continuum of solutions branching at $t_0$ and with derivatives
\[
\begin{pmatrix}\dot{\widehat{y}}\\ \dot a_1\\ \dot a_2 \end{pmatrix}\in \left\{-\left(\theta \,\widehat{n}_{+1} + (1-\theta)\widehat{n}_{+2} \right)\dot l(t): \theta\in[0,1]\right\}.
\]
\item \label{enum:toy_example_softening} If both $h_1, h_2>1$ then there may be up to three solutions branching at $t_0$, namely:
\begin{enumerate}[1.]
\item \label{enum:toy_example_softening_concentrated1} plastic deformation becomes localized in spring 1: 
\begin{equation*}
\begin{pmatrix}\dot{\widehat{y}}\\ \dot a_1\\ \dot a_2 \end{pmatrix} = -\frac{1}{1+(k_{1}^{-1}+k_{2}^{-1})(1-h_1)s_1^2}\, \dot l(t)\,\widehat{n}_{+1},\quad \text{possible iff}
\quad(k_{1}^{-1}+k_{2}^{-1})(1-h_1)s_1^2>-1.
\end{equation*}
\item  \label{enum:toy_example_softening_concentrated2} plastic deformation becomes localized in spring 2:
\begin{equation*}
\begin{pmatrix}\dot{\widehat{y}}\\ \dot a_1\\ \dot a_2 \end{pmatrix} = -\frac{1}{1+(k_{1}^{-1}+k_{2}^{-1})(1-h_2)s_2^2}\, \dot l(t)\,\widehat{n}_{+2},
\quad \text{possible iff}
\quad(k_{1}^{-1}+k_{2}^{-1})(1-h_2)s_2^2>-1,
\end{equation*}
\item  \label{enum:toy_example_softening_distributed} distributed solution \eqref{eq:example1_sol3}, which is possible iff
\begin{equation}
(1-h_1)s_1^2 + (1-h_2)s_2^2+ (k_{1}^{-1}+k_{2}^{-1})(1-h_1)(1-h_2)s_1^2s_2^2<0.
\label{eq:example1_sol3_good_lambda}
\end{equation}
\end{enumerate} 
\end{enumerate}
\end{Proposition}
We must note, that because any point $t$ along the solution \eqref{eq:example1_sol3} satisfies \eqref{eq:example1_ic_in_R3}, in the case {\it \ref{enum:toy_example_softening}} there technically is a continuum of solutions with branching from \eqref{eq:example1_sol3} at every $t\geqslant t_0$.  Also, note that our conclusions in the case {\it \ref{enum:toy_example_softening}} agree with the conclusions of \cite{ChenBaker2004}.

\noindent  {\bf The idea of the proof of Proposition \ref{prop:solutions_of_ex1}.} To save space we will omit complete derivations, as they are straightforward, yet cumbersome and very similar to the proof of Proposition \ref{prop:sp_single}. The following applies to a.a. $t\in[t_0,T]$. We first consider the case of a solution, driven by the facet of $\mathcal{C}_V$ with the normal vector $\widehat{n}_{+1}$, i.e. we must find the unknown quantities $\dot{\widehat{y}}, \dot a_1, \dot a_2\in \mathbb{R}$ such that
\[
\widehat{n}_{+1}^{\top} \begin{pmatrix}\dot{\widehat{y}} + \kappa\, \dot l(t)\\[1mm] (1-h_1)\dot a_1\\[1mm] (1-h_2) \dot a_2 \end{pmatrix}=0,
\]
\[
-\begin{pmatrix}\dot{\widehat{y}}\\ \dot a_1\\ \dot a_2 \end{pmatrix}=\widehat{n}_{+1}\, \lambda, \qquad \text{for some }\lambda\geqslant 0,
\]
yet such solution must not violate the constraint with $\widehat{n}_{+1}$ in order to stay in the moving set $\mathcal{C}_{V}$:
\begin{equation}
\widehat{n}_{+2}^{\top}\begin{pmatrix}\dot{\widehat{y}} + \kappa\, \dot l(t)\\[1mm] (1-h_1)\dot a_1\\[1mm] (1-h_2) \dot a_2 \end{pmatrix}\leqslant 0.
\label{eq:example1_sol1_2nd_constraint}
\end{equation}
By manually solving the equations we arrive to the localized solution \ref{enum:toy_example_softening_concentrated1} 
and its associated condition, present in the case \ref{enum:toy_example_softening}. 
It is also present in the case \ref{enum:toy_example_pp}, as a member of the family of solutions with $\theta=1$, but \eqref{eq:example1_sol1_2nd_constraint} plays a role and excludes such solution when $h_1<1$.

Similarly, by manually solving 
\[
\widehat{n}_{+2}^{\top} \begin{pmatrix}\dot{\widehat{y}} + \kappa\, \dot l(t)\\[1mm] (1-h_1)\dot a_1\\[1mm] (1-h_2) \dot a_2 \end{pmatrix}=0,
\]
\[
-\begin{pmatrix}\dot{\widehat{y}}\\ \dot a_1\\ \dot a_2 \end{pmatrix}=\widehat{n}_{+2}\, \lambda, \qquad \text{for some }\lambda\geqslant 0,
\]
and checking that
\[
\widehat{n}_{+1}^{\top}\begin{pmatrix}\dot{\widehat{y}} + \kappa\, \dot l(t)\\[1mm] (1-h_1)\dot a_1\\[1mm] (1-h_2) \dot a_2 \end{pmatrix}\leqslant 0
\]
we deduce that in the case \ref{enum:toy_example_softening} the localized solution \ref{enum:toy_example_softening_concentrated2} is present, and also we have it in the case   \ref{enum:toy_example_pp} as a member of the family of solutions with $\theta=0$.

Finally, by  solving 
\[
\widehat{n}_{+1}^{\top} \begin{pmatrix}\dot{\widehat{y}} +\kappa\, \dot l(t)\\[1mm] (1-h_1)\dot a_1\\[1mm] (1-h_2) \dot a_2 \end{pmatrix}=0, \qquad \widehat{n}_{+2}^{\top} \begin{pmatrix}\dot{\widehat{y}} + \kappa\, \dot l(t)\\[1mm] (1-h_1)\dot a_1\\[1mm] (1-h_2) \dot a_2 \end{pmatrix}=0,
\]
\[
-\begin{pmatrix}\dot{\widehat{y}}\\ \dot a_1\\ \dot a_2 \end{pmatrix}=\widehat{n}_{+1}\, \lambda_1+\widehat{n}_{+2}\, \lambda_2 \qquad \text{for some }\lambda_1,\lambda_2\geqslant 0.
\]
and carefully treating all the possibilities we arrive to the conclusions of the cases 
\ref{enum:toy_example_hardening}, \ref{enum:toy_example_pp}, and the distributed solution \ref{enum:toy_example_softening_distributed} in the case \ref{enum:toy_example_softening}. $\blacksquare$

It is, of course, possible to complete the statement of Proposition \ref{prop:solutions_of_ex1} with the cases, where the springs do not have the same type of plasticity, but such study is beyond our current goal of merely showing non-unique solutions emerging via the nonsmooth bifurcation.

\begin{figure}[H]\center
\includegraphics{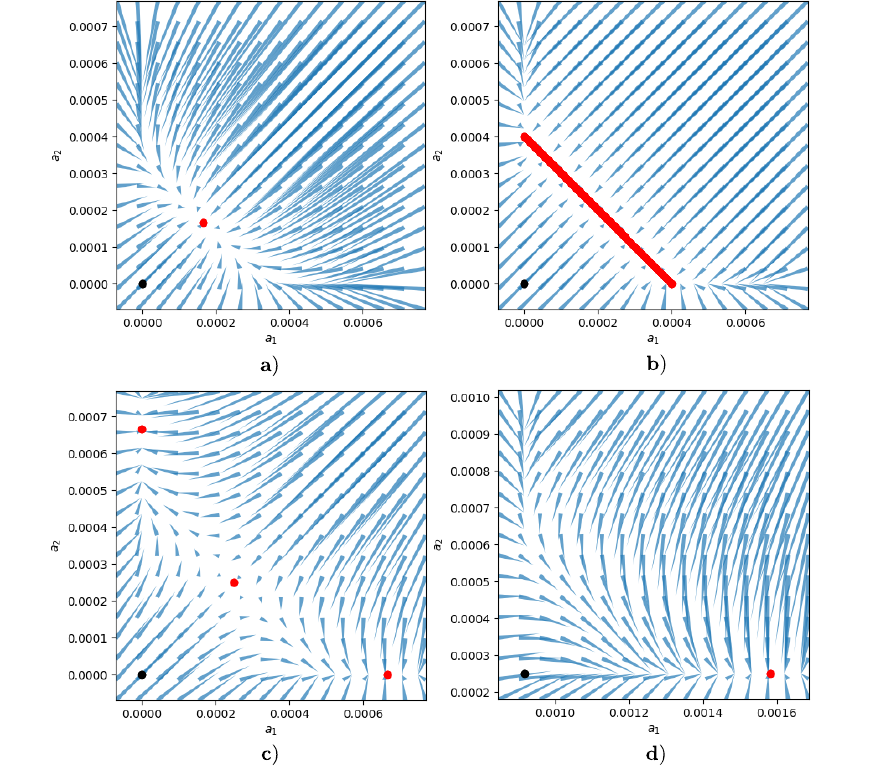}
\caption{
\footnotesize Chart of mapping $T_j$ in the implicit catch-up scheme for the two-springs system.  The projected value of $a$ (i.e. $a^{j-1}$ in terms of \eqref{eq:T_i}) is shown as the black point, while the red points represent the fixed points of $T_j$ (i.e. $a^{j}$). Each wedge connects a sample value $a$ and the corresponding value $T_j a$. {\bf a)} Case of hardening $h_1=h_2=0.8$.  {\bf b)} Case of perfect plasticity $h_1=h_2=1$. {\bf c)} Case of softening $h_1=h_2=1.2$ with $a^{j-1}$ satisfying \eqref{eq:example1_ic_in_R3} (i.e. belonging to the distributed solution {\it \ref{enum:toy_example_softening_distributed}} from Proposition \ref{prop:solutions_of_ex1}). Observe the saddle fixed point with its stable and unstable sets. {\bf d)} Case of softening $h_1=h_2=1.2$ with $a^{j-1}$ satisfying only the first equation of  \eqref{eq:example1_ic_in_R3} (i.e. belonging to the localized solution {\it \ref{enum:toy_example_softening_concentrated1}} from Proposition \ref{prop:solutions_of_ex1}).  
} \label{fig:paper1_example1_stability_chart}
\end{figure}

By using the numerical implicit catch-up scheme we can immediately observe additional properties of the solutions, namely, their {\it stability} is the sense of the stability of the corresponding fixed point of the map $T_j$. Specifically, we take 
\[
k_1=k_2=1,\qquad s_1=s_2 =1, \qquad c_{0,1}=c_{0,2}=0.01, 
\]
initial condition
\[
\sigma_1(0)=\sigma_2(0)=0.01, \qquad a_1(0)=a_2(0)=0,\qquad  p_1(0)=p_2(0)=0,
\]
and in Fig. \ref{fig:paper1_example1_stability_chart} we plot mapping $T_j$, where $j=1, t_0=0, t_1= 0.0004$. We examine the cases of $h$ as in Proposition \ref{prop:solutions_of_ex1}. From Fig. \ref{fig:paper1_example1_stability_chart} we make the following Observation \ref{obs:stability_of_the_fixed_points}.

\begin{Observation} 
\label{obs:stability_of_the_fixed_points}
Consider map $T_j$ as \eqref{eq:T_i} for the two-springs system such that $\begin{pmatrix}y^{j-1}\\a^{j-1}\end{pmatrix}$ satisfies \eqref{eq:example1_ic_in_R3}--\eqref{eq:example1_ic_in_R3_not_broken}. From the numerical results presented in Fig. \ref{fig:paper1_example1_stability_chart} we conclude that $T_j$
undergoes a bifurcation when $h_1=h_2$ go through value 1. Specifically, we see that
\begin{enumerate}[{\it A.}]
\item in the case of hardening $T_j$ has a single fixed point (Fig. \ref{fig:paper1_example1_stability_chart} a) which is stable is the sense of iterations of $T_j$ and which corresponds to the unique solution of Proposition \ref{prop:solutions_of_ex1}  {\it \ref{enum:toy_example_hardening}}.
\item in the case of perfect plasticity $T_j$ has an interval of fixed points, which is an attracting set for the iterations of $T_j$ (Fig. \ref{fig:paper1_example1_stability_chart} b). This agrees with Proposition \ref{prop:solutions_of_ex1} {\it \ref{enum:toy_example_pp}}.
\item \label{toy_example_observation_softening}
in case of softening and $a_j$ satisfying \eqref{eq:example1_ic_in_R3} the mapping $T_j$ has three fixed points (Fig. \ref{fig:paper1_example1_stability_chart} c), two of which are stable and they correspond to the localized solutions {\it \ref{enum:toy_example_softening_concentrated1}} and {\it \ref{enum:toy_example_softening_concentrated2}} of Proposition \ref{prop:solutions_of_ex1}. The third fixed point is of the saddle type and it corresponds to the distributed solution {\it \ref{enum:toy_example_softening_distributed}}. However,
\begin{enumerate}
\item
if for the next time-step $j'=j+1$ the projected point $a^{j'-1}$ is taken as the saddle fixed point (i.e. the numerical solution follows  {\it \ref{enum:toy_example_softening_distributed}}), then $T_{j'}$ will again have three fixed points of the same types, 
\item if for the next time-step $j'=j+1$ the projected point $a^{j'-1}$ is taken as one of the stable fixed points (i.e. the numerical solution follows {\it \ref{enum:toy_example_softening_concentrated1}} or {\it \ref{enum:toy_example_softening_concentrated2}}) then $T_{j'}$ will have only one stable fixed point (Fig. \ref{fig:paper1_example1_stability_chart} d), corresponding to the continuation of the branched localizeded solution.
\end{enumerate}
\end{enumerate}
\end{Observation}
\begin{Remark}
We note that bifurcation with the birth of two additional fixed points and the change of stability in the existing fixed point is similar to the classical pitchfork bifurcation (see e.g. \cite[Fig. 3.4.2, p. 64, Fig. 8.1.6, p. 273]{Strogatz2024}). However, the map $T_j$ is nonsmooth, therefore the classical analysis of bifurcations via the eigenvalues of its linearization is not applicable.
\end{Remark}
Finally, recall Remark \ref{remark:initial_guess} and observe from Fig.  \ref{fig:paper1_example1_stability_chart} b) and c) how different choices of the initial guess ${\rm InitialValue}(j, a^{j-1})$ in Algorithms \ref{alg:implicit-catch-up1} and \ref{alg:implicit-catch-up-reduced-dim} can lead to different fixed points found numerically in this example.

\subsection{Non-contractiveness of the iterated map near the saddle point due to coupling of springs}

The following question arises naturally:

\noindent {\bf Question.} {\it From both analytical (Proposition \ref{prop:solutions_of_ex1}) and numerical (Fig. \ref{fig:paper1_example1_stability_chart}) evidence
we deduce that in the case of softening there can be three branching solutions when  $k_1, k_2>0$ are fixed and $s_1, s_2>0$ are not too large. However, for individual (uncoupled) springs, by taking smaller and smaller $s>0$ we can make the Lipschitz constant of the map $\widetilde{a} \mapsto {\rm proj}^{\mathbb{K}^{-1}}(\Sigma, \mathcal{C}(t, \widetilde{a}))$ arbitrarily small for a general $\Sigma\in \mathbb{R}^2$, unless the projection goes to the point of complete failure. One can observe this from Fig. \ref{fig:lin_soft_yielding_and_sweeping} c) and \ref{fig:complete_failure} a), where the boundary of $\mathcal{C}(t, \widetilde{a})$ would approach the vertical line as $s>0$ decreases. Now, in the case of two coupled springs, why it is not possible to choose $s_1, s_2$ small enough to make map $T_j$ a contraction?}

\begin{figure}[H]\center
\includegraphics{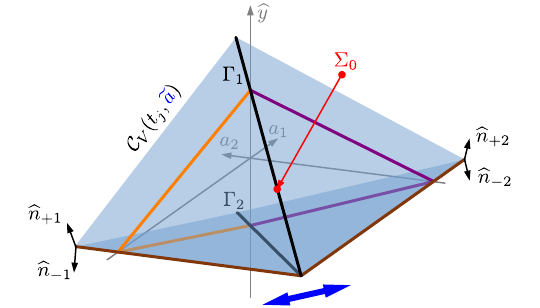}
\caption{
\footnotesize Moving set $\mathcal{C}_V(t, a)$ of the sweeping process \eqref{eq:example1_normals_in_R3}--\eqref{eq:example1_C_in_R3}, corresponding to the two-springs system and (shown in red) a projection of a point $\Sigma_0$ onto the set (i.e. map \eqref{eq:full_iteration_map_reduced_dim} for $\Sigma_0=\begin{pmatrix}\widehat{y}^{j-1}\\a^{j-1}\end{pmatrix}$). Set $\mathcal{C}_V$ is a representation of $\mathcal{C}$ in coodinates of a three-dimensional subspace $\mathcal{V}\times \mathbb{R}^2\subset \mathbb{R}^4$. Brown edges are the states of complete failure. Orange and purple angles show the coupled constraints of springs 1 and 2, respectively (one can see them as two instances of constraints of Fig. \ref{fig:lin_soft_yielding_and_sweeping} b). By $\Gamma_1$ and $\Gamma_2$ we denote the ``diagonal'' edges, shown in black. Blue arrow indicates the direction along to the unstable set of the saddle fixed point in Fig. \ref{fig:paper1_example1_stability_chart} c). 
} \label{fig:example1_moving_set}
\end{figure}

To answer this question we sketch the moving set $\mathcal{C}_V(t, a)$ of the sweeping process \eqref{eq:sp_in_V}, which, in case of the two-springs system, is a subset of $\mathbb{R}^3$ given by \eqref{eq:example1_normals_in_R3}--\eqref{eq:example1_C_in_R3}, see Fig. \ref{fig:example1_moving_set}.

From Fig. \ref{fig:example1_moving_set} we can provide a visual geometric answer to the question. The following is not a formal proof, but rather an observation from Figs. \ref{fig:example1_moving_set} and \ref{fig:paper1_example1_stability_chart} c. However, we still think that it is a valuable explanation of how {\it coupling} of two spings can lead to non-contracting map $T_j$ and, therefore, multiple solutions.

 Observe,  that no matter how small we take $s_1, s_2>0$ (i.e. no matter how pointed are the orange and purple angles, similar to the angle in Fig. \ref{fig:lin_soft_yielding_and_sweeping} b), the edges $\Gamma_1, \Gamma_2$ are still present. Let $h_1,h_2>1$ satisfy \eqref{eq:example1_sol3_good_lambda}, a projected point $\Sigma_0=\begin{pmatrix}\widehat{y}^{j-1}\\a^{j-1}\end{pmatrix}$ satisfy \eqref{eq:example1_ic_in_R3}--\eqref{eq:example1_ic_in_R3_not_broken}, which means that $T_j$ has a saddle fixed point $a^*$, see Fig. \ref{fig:paper1_example1_stability_chart} c). Then the projection 
\begin{equation}
\begin{pmatrix}{\widehat y}^*\\ a^*\end{pmatrix}:={\rm proj}\left(\Sigma_0, \mathcal{C}_V\left(t_j, a^*\right)\right)
\label{eq:projection-to-fixed-point-on-the-question}
\end{equation}
lays exactly on edge $\Gamma_1$, as shown in Fig. \ref{fig:example1_moving_set}.

 A small perturbation $\Delta\in \mathbb{R}^2$ of $a^*$ translates the moving set 
by the vector $\begin{pmatrix}0 \\H \Delta\end{pmatrix}$. Take $\Delta$ such that $\begin{pmatrix}0 \\H \Delta\end{pmatrix}$ is parallel to the blue direction in Fig. \ref{fig:example1_moving_set}, i.e. it is parallel to the unstable set of the saddle point in Fig. \ref{fig:paper1_example1_stability_chart} c). When the magnitude of $\Delta$ is small enough, the projection 
\[
\begin{pmatrix}\widehat{y}_{\Delta} \\ a_\Delta\end{pmatrix} :={\rm proj}\left(\Sigma_0, \mathcal{C}_V\left(t_j, a^*+\Delta\right)\right)={\rm proj}\left(\Sigma_0, \mathcal{C}_V \left(t_j, a^*\right)+\begin{pmatrix}0\\ H\Delta\end{pmatrix}\right)
\]
will remain on the edge $\Gamma_1$, which is translated together with the whole set by $\begin{pmatrix}0 \\H \Delta\end{pmatrix}$. 

Projection \eqref{eq:projection-to-fixed-point-on-the-question} means that
\begin{equation}
\Sigma_0 - \begin{pmatrix}{\widehat y}^*\\ a^*\end{pmatrix}\in N_{\mathcal{C}_{V}(t_j,a^*)}\left(\begin{pmatrix}{\widehat y}^*\\ a^*\end{pmatrix}\right),
\label{eq:nc-on-the-question}
\end{equation}
where the normal cone in the right-hand side becomes very narrow for small $s_1,s_2>0$. This can be observed from Fig. \ref{fig:example1_moving_set} as the entire set $\mathcal{C}_{V}(t_j,a^*)$ would become ``flatter'' and closer to the plane $a_1$,$a_2$. The narrow normal cone means that $\Sigma_0 -\begin{pmatrix}\widehat{y}^*\\ a^*\end{pmatrix}$ (parallel to the red vector in Fig. \ref{fig:example1_moving_set}) is nearly orthogonal to $\begin{pmatrix}0 \\H\Delta\end{pmatrix}$ (parallel to the blue diection). Therefore, for $H\Delta$ of small magnitudes
\[
\begin{pmatrix}\widehat{y}_\Delta\\ a_\Delta\end{pmatrix} - \begin{pmatrix}\widehat{y}^*\\ a^*\end{pmatrix} \approx \begin{pmatrix}0 \\H \Delta\end{pmatrix}
\]
i.e.
\[
T_j(a^*+\Delta) - T_j(a^*)=a_\Delta - a^* \approx H \Delta.
\]
Since we had $h_1, h_2>1$, we cannot make $T_j$ to be a contraction mapping, no matter how small we take $s_1,s_2>0$. Notice, that the slopes, controlled by $s_1, s_2$ (orange and purple in Fig. \ref{fig:example1_moving_set}) will essentially affect $a_\Delta$ only when $\Delta$ is large enough for the projection $\begin{pmatrix}\widehat{y}_\Delta\\ a_\Delta\end{pmatrix}$ to no longer lay on the edge $\Gamma_1$, however, the ``large enough'' itself depends on how narrow is the normal cone in \eqref{eq:nc-on-the-question}, which, in turn, depends on $s_1, s_2$.

\section{Numerical results for regular lattices: strain localization into a shear band}
\label{sect:LSM-examples}
We apply Algorithm \ref{alg:implicit-catch-up-reduced-dim} to solve for the evolution in regular lattices and lattices with defects. 
\subsection{Rectangular lattice} Consider a two-dimensional rectangular lattice of $10\times 15$ nodes spaced with the horizontal and vertical step of $0.5$, see Fig. \ref{fig:paper1_big_lattice1_softening} a. Apply the following boundary condition: all the nodes in the lowest row have their $y$-coordinate restrained to $0$ and all the nodes in the uppermost row have their $y$-coordinate restrained to a value, which monotonically increases with time. We also restrain the $x$-coordinate of a single node to a constant value in order to achieve kinematic determinacy (Assumption \ref{ass:trivial_kernel_intersect}). The following values characterize the lattice:
\begin{align*}
\text{number of spatial dimensions:}&&d&=2,&\\
\text{number of nodes:}&&n&=150,&\\
\text{number of springs:}&&m&=527,&\\
\text{number of external displacement constraints:}&&q&=21,&\\
\text{dimension of the space of self-stresses:}&&{\rm dim}\, \mathcal{V}&=248.&
\end{align*}
We will consider several cases.
\subsubsection{Symmetric lattice with softening} \label{sssect:lattice1_softening} 
\begin{figure}[H]\center
\includegraphics{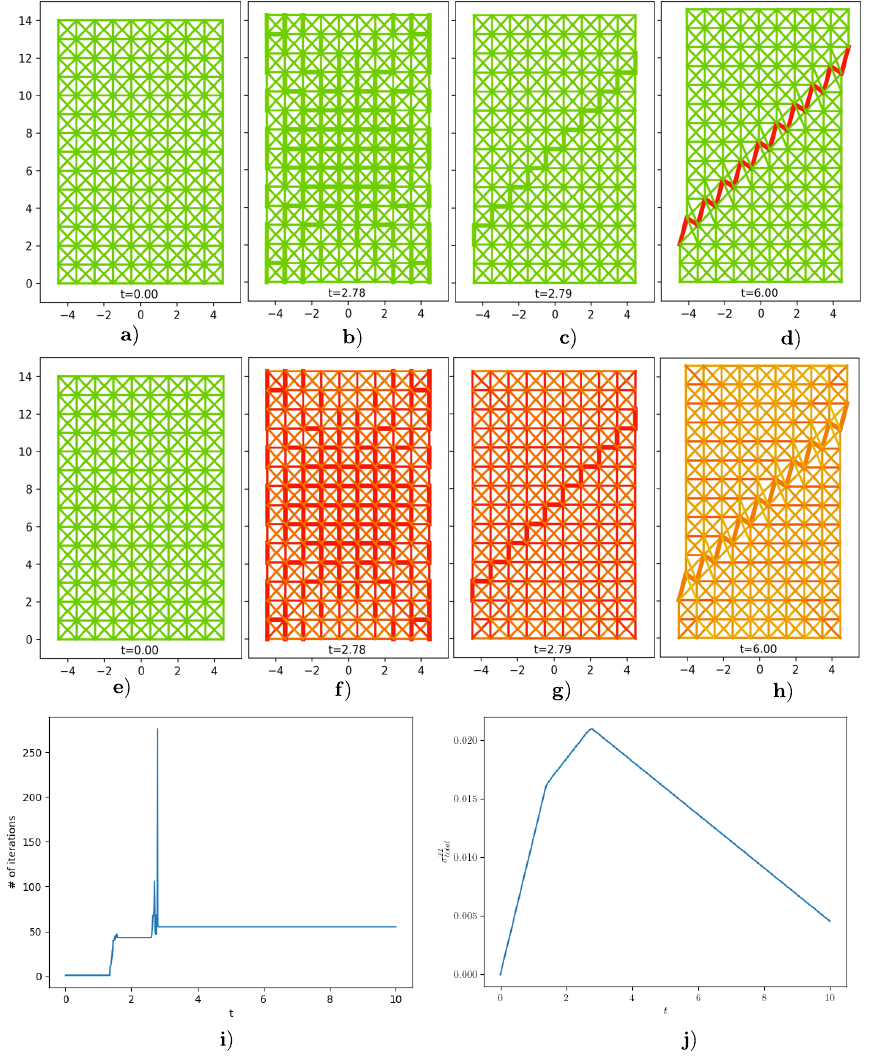}
\caption{\footnotesize Evolution of a {\bf symmetric} rectangular lattice with {\bf softening} under increasing displacement load along the vertical axis (see Section \ref{sssect:lattice1_softening}). {\bf a)}--{\bf d)} Thicker lines denote the springs in the plastic regime, and color denotes the value of the damage variable $a$ (green$=0$, red $=\max\limits_{t\in [0,6]} a_i(t)$) {\bf e)}--{\bf h)} Thicker lines denote the springs in the plastic regime, and color denotes the absolute value of the stress $\sigma$ (green$=0$, yellow$=$an intermediate value, red $=\max\limits_{t\in [0,6]} |\sigma_i(t)|$). {\bf i)} the number of iterations of map \eqref{eq:full_iteration_map_reduced_dim} for each time-step. {\bf j)} total stress \eqref{eq:total-stress-formula} of the lattice in the vertical direction.
} \label{fig:paper1_big_lattice1_softening}
\end{figure}

See Fig. \ref{fig:paper1_big_lattice1_softening}, where all of the springs have the following parameters:
\begin{equation}
k_i=1, \quad h_i=1.1, \quad s_i = 0.3, \quad c_{0,i} = 0.01, \qquad i\in \overline{1,m}.
\label{eq:big-lattice1-softening-spring-parameters}
\end{equation}

We would like to illustrate the evolution of the damage variable $a$ (Fig. \ref{fig:paper1_big_lattice1_softening} a--d), which is connected to plastic deformation, and the evolution of stress $\sigma$ (Fig. \ref{fig:paper1_big_lattice1_softening} e--f). Additionally, Fig. \ref{fig:paper1_big_lattice1_softening} i shows the number of iterations of map $T_j$ (see \eqref{eq:T_i} and \eqref{eq:full_iteration_map_reduced_dim}) required for convergence at each time-step of Algorithm \ref{alg:implicit-catch-up-reduced-dim}. Finally, Fig \ref{fig:paper1_big_lattice1_softening} j shows the total stress of the lattice along the vertical direction, which we compute as
\begin{equation}
\sigma_{total}^{22} = \frac{1}{A} \sum_{i\in \overline{1,m}} \left(\sigma_i \mathcal{D}_{i2}\right)\left((\varphi(\xi_0))_i\mathcal{D}_{i2}\right),
\label{eq:total-stress-formula}
\end{equation}
where the constant $A=(15-1)(10-1)=126$ is the area of the lattice in the reference configuration $\xi_0$, function $\varphi$ is given by \eqref{eq:distanceFunction}, and matrix $\mathcal{D}$ is given by \eqref{eq:unit_vectors_along_springs}.  For details on formula \eqref{eq:total-stress-formula} we refer to \cite[Section 7.2]{Gudoshnikov2023preprint}, \cite[(6)]{Lemaitre2017}, \cite[(1.2)]{virial_stress_1}, \cite[p.~6]{virial_stress_2}.

During the evolution of the symmetric rectangular lattice with softening from the relaxed reference configuration (Fig. \ref{fig:paper1_big_lattice1_softening} a)  we observe the elastic phase, followed by the accumulation of symmetrically distributed plastic deformation. The symmetric plastic phase happens from $t=1.36$ to $t=2.78$, and it can be distinguished on the graph of the number of iterations (Fig \ref{fig:paper1_big_lattice1_softening} i) and on the graph of total stress  (Fig \ref{fig:paper1_big_lattice1_softening} j). At $t=2.78$ we observe a sudden loss of symmetry of the solution and the beginning of strain localization into the diagonal shear band (Fig \ref{fig:paper1_big_lattice1_softening} b--c, f--g). This corresponds to the sharp peak in the number of iterations and to the maximal total stress. All further deformation will happen in the shear band, and its elements remain in the plastic regime with softening, while the rest of the lattice unloads elastically (Fig \ref{fig:paper1_big_lattice1_softening} h). This way, the total stress decreases linearly with the increasing displacement load after $t=2.78$ (Fig. \ref{fig:paper1_big_lattice1_softening} j), eventually leading to the state of complete failure.

Notice, that we have obtained a non-symmetric solution in a symmetric problem, which suggests non-uniqueness of solutions. And, indeed, by modyfing the initial guess during the loss of symmetry (for $t$ near $2.78$, cf. Remark \ref{remark:initial_guess} on {\rm InitialValue}) we are able to obtain another solution, in which the shear band develops in the opposite direction (see Fig. \ref{fig:paper1_big_lattice1_softening_solution2}).

We also expect that there may be symmetric solutions to the problem, which are numerically unstable, similarly to the situation in the toy system, see Fig. \ref{fig:paper1_example1_stability_chart} c and Observation \ref{obs:stability_of_the_fixed_points} \ref{toy_example_observation_softening}

Full videos and computer programs of numerical simulations of this and other examples are incuded as supplemental materials \cite{SupplMat}.

\subsubsection{Lattice with a defect with softening} \label{sssect:lattice1_softening_defect} 
\begin{figure}[H]\center
\includegraphics{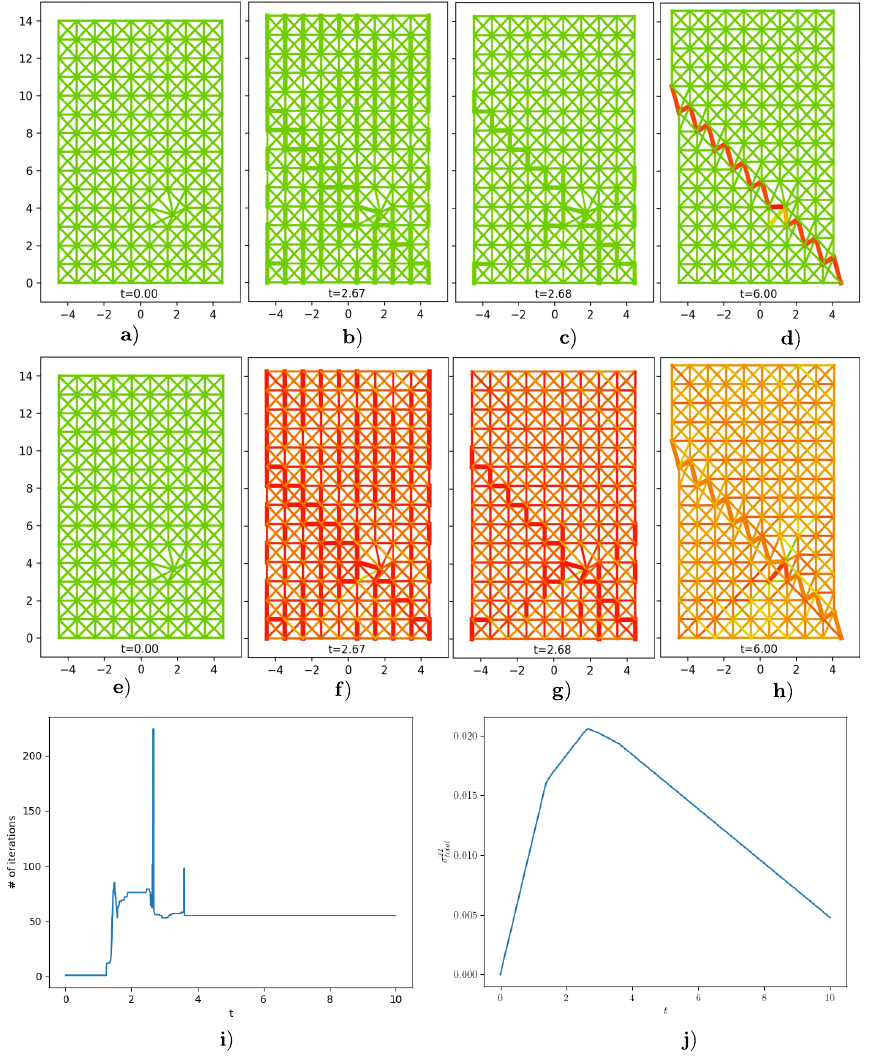}
\caption{\footnotesize 
Evolution of a rectangular lattice with a {\bf defect} with {\bf softening} under increasing displacement load along the vertical axis (see Section \ref{sssect:lattice1_softening_defect}). The presentation is similar to Fig. \ref {fig:paper1_big_lattice1_softening}.
} \label{fig:paper1_big_lattice1_softening_defect}
\end{figure}

See Fig. \ref{fig:paper1_big_lattice1_softening_defect}, where all of the springs have the same parameters \eqref{eq:big-lattice1-softening-spring-parameters}, but we introduce a geometric defect into the lattice to destroy its symmetry. We observe essentially the same behavior as in the case of symmetric lattice with softening, but the diagonal shear band now forms to include the location of the defect. It is also accompanied by the minor shear band, formed in the opposite direction (see Fig.~\ref{fig:paper1_big_lattice1_softening_defect} g and the supplemental material). The dissolution of the minor shear band at $t=3.61$ is marked by another peak in the number of iterations  (see Fig.~\ref{fig:paper1_big_lattice1_softening_defect} i) and a small, but noticable change of the total plasticity modulus of the lattice (see Fig.~\ref{fig:paper1_big_lattice1_softening_defect} j) .

\subsubsection{Symmetric lattice with perfect plasticity} \label{sssect:lattice1_pp} See Fig. \ref{fig:paper1_big_lattice1_pp}, where all of the springs have the following parameters:
\begin{equation}
k_i=1, \quad h_i=1, \quad s_i = 0.3, \quad c_{0,i} = 0.01, \qquad i\in \overline{1,m}.
\label{eq:big-lattice1-pp-spring-parameters}
\end{equation}
Notice, that no localization is present, and the deformation, as well as stress, is symmetrically distributed in the lattice. We know that the evolution of stresses is unique for elastic-perfectly plastic lattices \cite[Th.~4.1, Th.~5.1]{Gudoshnikov2023preprint}, as well as for elastic-perfectly plastic bodies in general. We do not exclude the possibility that there could be other solutions to the problem with the same trajectory of the stress variable (Fig. \ref{fig:paper1_big_lattice1_pp} e--h, j), but different trajectories of the damage variable (Fig. \ref{fig:paper1_big_lattice1_pp} a--d), and, therefore, different displacements. However, at the moment we were not able to find such alternative solutions in the current example.

\subsubsection{Lattice with a defect with perfect plasticity} \label{sssect:lattice1_pp_defect} See Fig.~\ref{fig:paper1_big_lattice1_pp_defect}, where all of the springs have the same parameters \eqref{eq:big-lattice1-pp-spring-parameters}, but we introduce a geometric defect into the lattice to destroy its symmetry. Again, the shear band is now formed across the location of the defect (Fig.~\ref{fig:paper1_big_lattice1_pp_defect} d). Unlike the softening case (Fig. \ref{fig:paper1_big_lattice1_softening} h), there is no unloading at the later stages of the evolution (Fig.~\ref{fig:paper1_big_lattice1_pp_defect} h). 

\subsubsection{Symmetric lattice with hardening} \label{sssect:lattice1_hardening} See Fig. \ref{fig:paper1_big_lattice1_hardening}, where all of the springs have the following parameters:
\begin{equation}
k_i=1, \quad h_i=0.9, \quad s_i = 0.3, \quad c_{0,i} = 0.01, \qquad i\in \overline{1,m}.
\label{eq:big-lattice1-hardening-spring-parameters}
\end{equation}
We have obtained a symmetrically distributed solution, similar to the symmetric example with perfect plastcity (Fig. \ref{fig:paper1_big_lattice1_pp}), but with more plastic deformation. In the case of hardening we know from the theory that the solution is unique for both stresses and displacements.  

\subsubsection{Lattice with a defect with hardening} \label{sssect:lattice1_hardening_defect} See Fig. \ref{fig:paper1_big_lattice1_hardening_defect}, where all of the springs have the same parameters \eqref{eq:big-lattice1-hardening-spring-parameters}, but we introduce a geometric defect into the lattice to destroy its symmetry. 

\begin{figure}[H]\center
\includegraphics{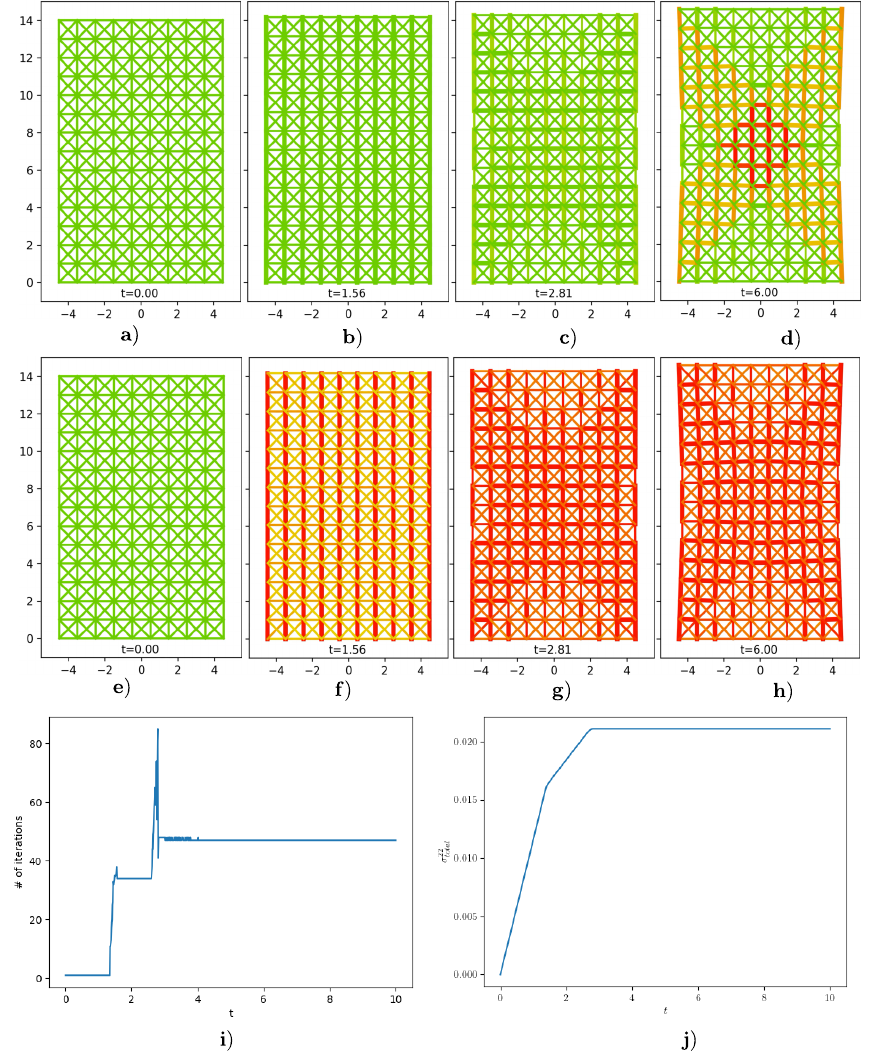}
\caption{\footnotesize 
Evolution of a {\bf symmetric} rectangular lattice with {\bf perfect plasticity} under increasing displacement load along the vertical axis (see Section \ref{sssect:lattice1_pp}). The presentation is similar to Fig. \ref {fig:paper1_big_lattice1_softening}.
} \label{fig:paper1_big_lattice1_pp}
\end{figure}

Notice, that the evolution of the damage has ``symmetric'' features (Fig. \ref{fig:paper1_big_lattice1_hardening_defect} d, cf. Fig. \ref{fig:paper1_big_lattice1_hardening} d and Fig. \ref{fig:paper1_big_lattice1_pp} d), but the localization of damage to the diagonal shear band crossing the defect is still visible (Fig. \ref{fig:paper1_big_lattice1_hardening_defect} c--d).

\begin{figure}[H]\center
\includegraphics{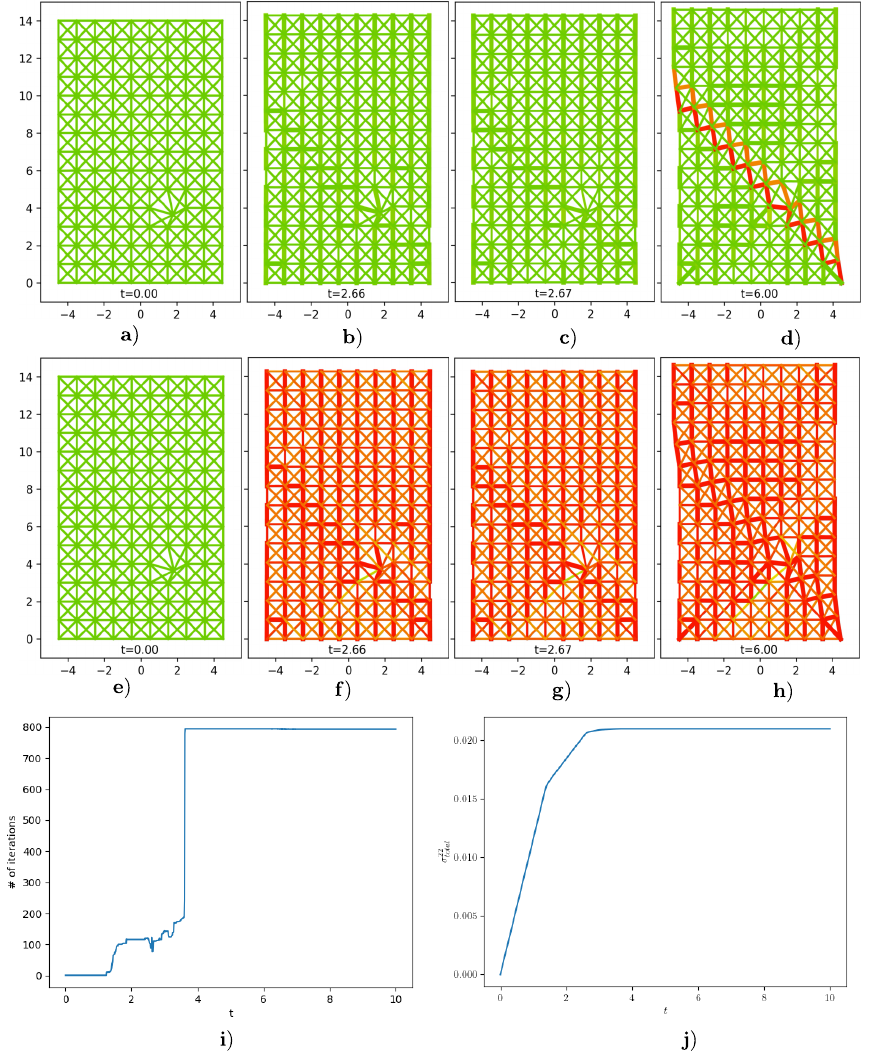}
\caption{\footnotesize 
Evolution of a rectangular lattice with a {\bf defect} with {\bf perfect plasticity} under increasing displacement load along the vertical axis (see Section \ref{sssect:lattice1_pp_defect}). The presentation is similar to Fig. \ref {fig:paper1_big_lattice1_softening}.
} \label{fig:paper1_big_lattice1_pp_defect}
\end{figure}
\begin{figure}[H]\center
\includegraphics{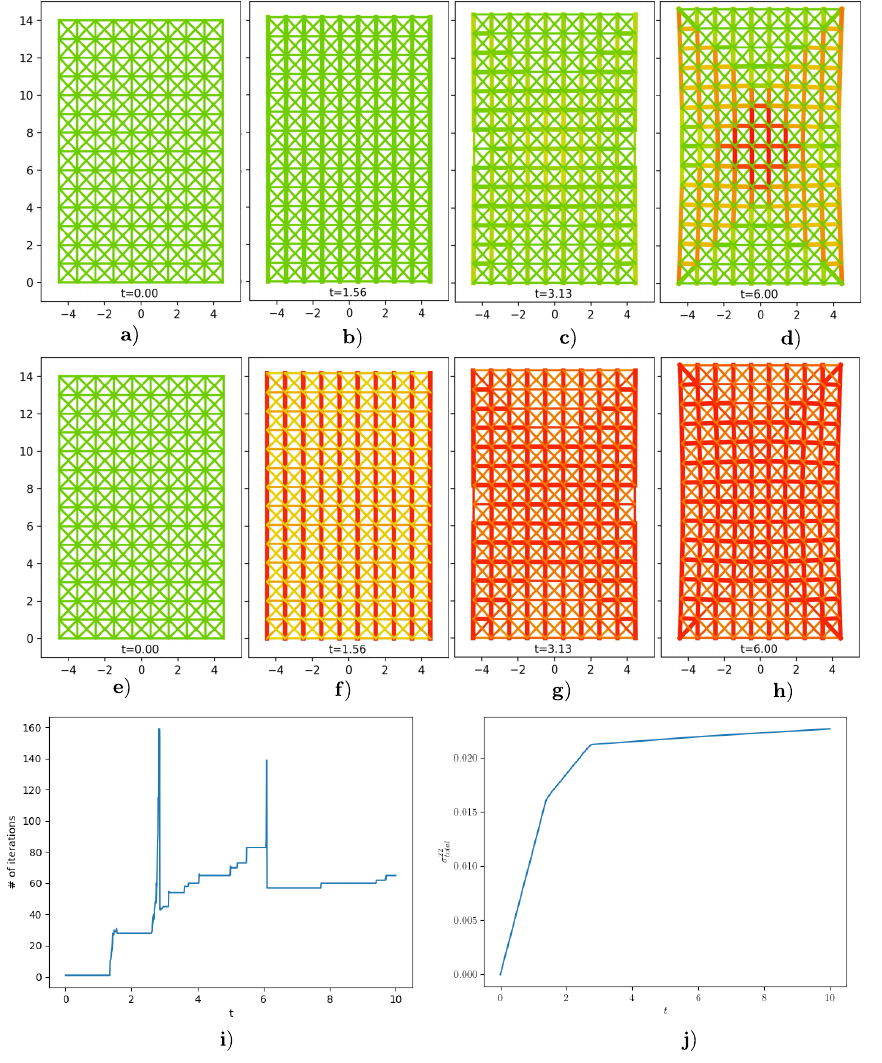}
\caption{\footnotesize 
Evolution of a {\bf symmetric} rectangular lattice with {\bf hardening} under increasing displacement load along the vertical axis (see Section \ref{sssect:lattice1_hardening}). The presentation is similar to Fig. \ref {fig:paper1_big_lattice1_softening}.
} \label{fig:paper1_big_lattice1_hardening}
\end{figure}
\begin{figure}[H]\center
\includegraphics{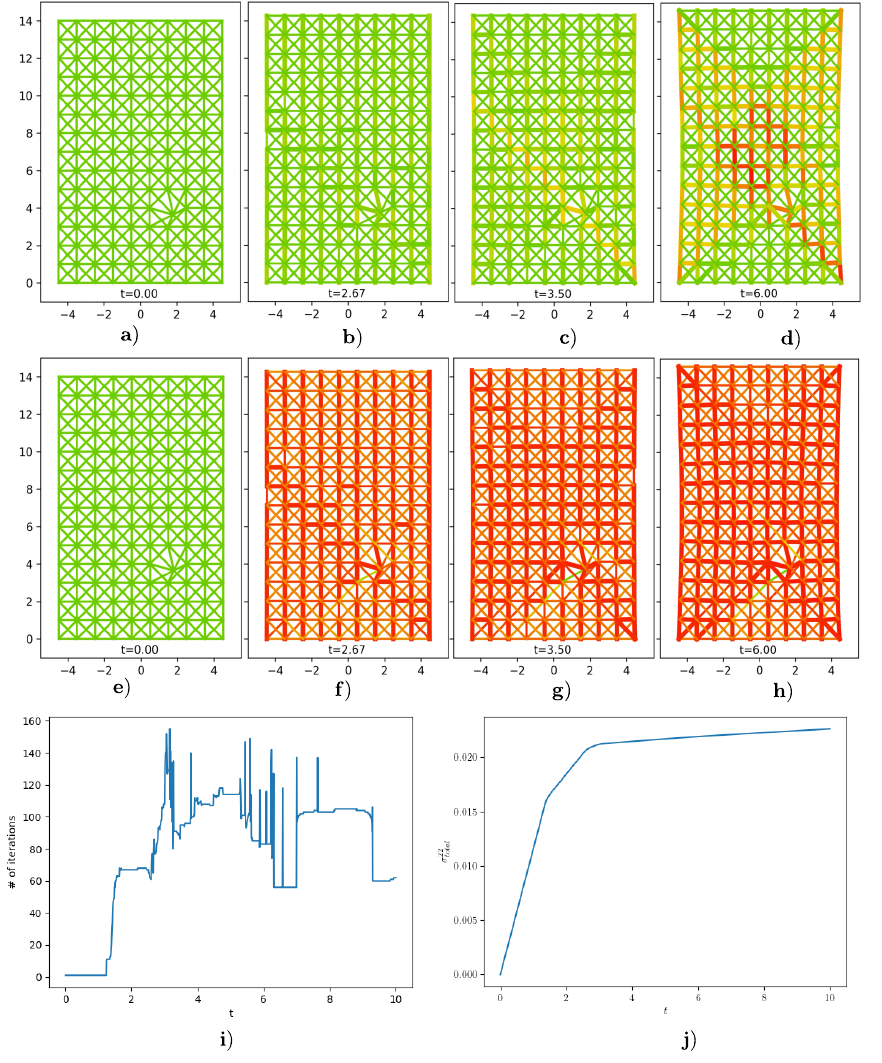}
\caption{\footnotesize 
Evolution of a rectangular lattice with a {\bf defect} with {\bf hardening} under increasing displacement load along the vertical axis (see Section \ref{sssect:lattice1_hardening_defect}). The presentation is similar to Fig. \ref {fig:paper1_big_lattice1_softening}.
} \label{fig:paper1_big_lattice1_hardening_defect}
\end{figure}

\subsection{Triangular lattice}
\label{ssect:lattice2}
\begin{figure}[H]\center
\includegraphics{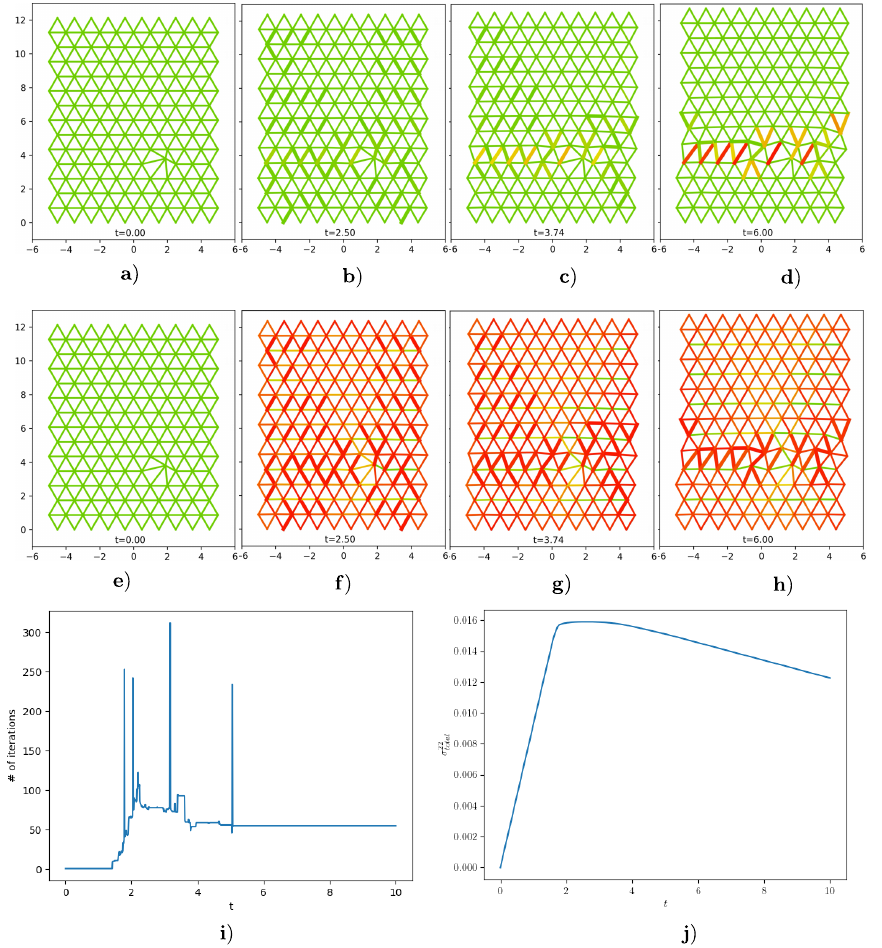}
\caption{\footnotesize 
Evolution of a triangular lattice with a defect with softening under increasing displacement load along the vertical axis (see Section \ref{ssect:lattice2}). The presentation is similar to Fig. \ref {fig:paper1_big_lattice1_softening}.
} \label{fig:paper1_big_lattice2_softening_defect}
\end{figure}

Consider now a similar triangular lattice (Fig. \ref{fig:paper1_big_lattice2_softening_defect} a) under the similar displacement loading along the vertical direction.  The following values characterize the lattice:
\begin{align*}
\text{number of spatial dimensions:}&&d&=2,&\\
\text{number of nodes:}&&n&=157,&\\
\text{number of springs:}&&m&=404,&\\
\text{number of external displacement constraints:}&&q&=21,&\\
\text{dimension of the space of self-stresses:}&&{\rm dim}\, \mathcal{V}&=111.&
\end{align*}
For the triangular lattice we will only consider the case of  defect included and springs with softening. All of the springs have the following parameters: 
\begin{equation*}
k_i=1, \quad h_i=1.1, \quad s_i = 0.2, \quad c_{0,i} = 0.01, \qquad i\in \overline{1,m}.
\label{eq:big-lattice2-softening-spring-parameters}
\end{equation*}

Similarly to the corresponding example with the rectangular lattice (Fig. \ref{fig:paper1_big_lattice1_softening_defect}) we observe the localization of damage and displacement into a shear band crossing the defect, and some elastic unloading in the rest of the lattice. However, the shear band is now horizontal, from which we conclude that the model output severely depends on the structure of the lattice.

\section{Conclusions and directions for future research}
\label{sect:conclusions}
We have proposed a state-dependent sweeping process model of an elastoplastic spring, which can be a hardening, perfectly plastic, and softening spring, depending on the value of the parameter $h$. Curiously enough, a progressive plastic deformation in the model is controlled by two parameters (the coefficient of state-dependent feedback $h$ and the geometric slope $s$), which naturally arise in the mathematical formulation, while their interplay results in a single observable characteristic, that being plasticity modulus $E^p$ of the spring, see \eqref{eq:modulus_formula_hs}. We will leave it to physicists and engineers to judge if $h$ and $s$ on their own could be interpreted as objective physical properties.

Further we have explicitly written an analogous state-dependent sweeping process \eqref{eq:sp} for a general Lattice Spring Model, which contains a finite number of such springs. The moving set of such sweeping process is convex (moreover, it is an unbounded {\it polytpe with fixed directions of the normal vectors to its faces}), but the state-dependence of the moving set is a trade-off for the non-convexity of the energy function in the variational formulations (see e.g. \cite{Bradies1999}, \cite{Lambrecht2003}, \cite{DalMaso2008Softening1}, \cite{DalMaso2008Softening2}). 

We succeeded to numerically solve the state-dependent sweeping process in our examples of LSMs by finding an asymptotically stable fixed point of the map $T_j$ (see \eqref{eq:T_i}) for each time-step $j$ of the implicit catch-up algorithm. Such numerical solutions can be computed up to an earliest point of complete failure of a spring in a lattice. 

This is in agreement with the fact that Theorem \ref{th:sp_state-dep} by Kunze and Monteiro Marques, being the cornerstone of the theory of state-dependent sweeping processes, only covers the case of hardening ($h<1$), and exactly excludes less trivial cases of perfect plasticity ($h=1$) and softening ($h>1$). In the future we plan to present a local version of the theorem, which would cover all of the three types and justify the existence of a fixed point for a time-step in the implicit catch-up algorithm, provided that none of the $m$ springs is near the state of complete failure. That would yield a solution to the time-continuous problem by the same argument from the proof of Theorem \ref{th:sp_state-dep}. 

We have  analytically shown non-uniqueness of solutions in the state-dependent sweeping process of the two-spring model of Fig.~\ref{fig:intro_fig_models}~b), and this agrees with the conclusions of \cite{ChenBaker2004}, where the same model was considered. Moreover, the nonsmooth bifurcation due to softening, which is discussed there and in many more engineering papers, can be now studied via the rigorous mathematical framework of a state-dependent sweeping process.

The importance of such study is backed by our numerical results: as one can see from Fig.~\ref{fig:paper1_example1_stability_chart}, the simultaneous change of the parameter $h$ across the value $1$ for both springs of the toy model  influences  dramatically not only the number of fixed points, but also their stability with respect to the map $T_j$. Other combinations of the three types of plasticity could be considered in the future to fully understand the stability of the fixed points, similarly to Fig.~\ref{fig:paper1_example1_stability_chart}. For a complete study of the bifurcation in the two-springs model we shall also include the parameters such that $k^{-1}s^2(1-h)<-1$ (see \eqref{eq:modulus_formula_hs} and Fig.~\ref{fig:possible_moduli}) and the reverse direction of loading ($\dot l<0$), which is related to the {\it snapback} phenomenon \cite{Crisfield1986}, \cite[Sect.~9]{Rots1989}.

 The next open question is to find a method to analytically determine the stability of the fixed points for a general LSM. It is a nontrivial task, as the map $T_j$ is not Fr{\'{e}}chet-differentiable (indeed, see \eqref{eq:T_i}), and the classical stability analysis via the eigenvalues of the linearization (see e.g. \cite[Section 10.1, pp.~386--387]{Strogatz2024}, \cite[Th.~2.34, p.~239]{Muller2015}) is not applicable.

On the other hand, Fig.~\ref{fig:paper1_example1_stability_chart} suggests a relation between {\it numerical stability} of the fixed points (i.e. in terms of the iterated map $T_j$), and their {\it mechanical stability}, i.e. their classification as minima, maxima and saddle points of the energy function. It is possible that the map $T_j$, which we sourced from the theory behind Theorem \ref{th:sp_state-dep}, actually decreases the energy of a state $\begin{pmatrix}y\\a\end{pmatrix}$ with each iteration, until a local minimum is found. This would mean a deeper connection between variational and sweeping process formulation of the mechanical models. An example with cycling iterations of $T_j$ would, however, disprove the hypothesis.  

In our numerical simulations of larger lattices we have observed emergent formation of shear bands in LSMs with softening, which generally agrees with the simulation results in engineering papers, e.g. \cite{deBorst1991}, \cite{deBorst2019}. Moreover, in a symmetric lattice with softening the shear band appears as an instantaneous loss of symmetry (Fig.~\ref{fig:paper1_big_lattice1_softening}~b--c,~f--g), accompanied by a spike in the iterations count of $T_j$ (Fig.~\ref{fig:paper1_big_lattice1_softening}~i). We have numerically confirmed non-uniqueness of that solution (cf. Fig.~\ref{fig:paper1_big_lattice1_softening}~d and Fig.~\ref{fig:paper1_big_lattice1_softening_solution2}). The orientation of the shear band is sensitive to the defects in the lattice (cf. Fig.~\ref{fig:paper1_big_lattice1_softening}~d and Fig.~\ref{fig:paper1_big_lattice1_softening_defect}~d) and heavily depends on the geometry of the lattice (cf. Fig.~\ref{fig:paper1_big_lattice1_softening_defect}~d and Fig.~\ref{fig:paper1_big_lattice2_softening_defect}~d). To verify the integral behavior of a lattice as an elasto-plastic body we have computed the evolution of total stress (Figs.~\ref{fig:paper1_big_lattice1_softening}--\ref{fig:paper1_big_lattice2_softening_defect} j), which reflects the respective type of plasticity. 

Naturally, the model we constructed could be adapted to non-isotropic and nonlinear hardening or softening by taking a state-dependent $h(a)$ or more complex (possibly even nonconvex) shape of ${\rm C}$ in \eqref{eq:tilde_C_single}. Instead of a linearization \eqref{eq:gc} near a fixed reference configuration of the lattice in Section \ref{ssect:kinematics} it is possible to consider an evolving geometry of the lattice with ``properly rotating'' springs, i.e. with the nonlinear state constraint, directly defined via the distance function \eqref{eq:distanceFunction}. Such model would have even more state-dependent terms, in particular the subspace $\mathcal{V}$ in \eqref{eq:K,moving_set}. The computational cost in such a model would likely prohibit the implicit catch-up algorithm, at least in its ``naive'' unoptimized implementation. Other numerical methods for nonsmooth problems  \cite{Outrata2013}, \cite{Farrell2020} may also be efficient. However, as we discussed above, stability of the fixed points in the implicit catch-up algorithm reveals the inner properties of the problem itself, thus the algorithm should not be  dismissed even in the presence of faster alternatives.

To mention at last, a similar state-dependent sweeping process formulation of continuous media with softening is conceivable, but it will require measure-valued state variables and appropriately defined governing equations. Such ``dual'' formulation would be complementary to the existing ``primal'' variational formulations, as, for example, both formulations together could help with two-sided error estimates, similarly to how it is for linear elliptic problems \cite[Th.~4.4, Figs.~4.8--4.9, pp.~64--67]{Krizek1996} and elasticity-perfect plasticity \cite{Haslinger2019}. Moreover, optimal control of state-dependent sweeping processes (currently being the scientific frontier, see e.g. \cite{Antil2025}) could be used to control the emergence of ductile fractures. 

\section*{Acknowledgements} The author thanks Pavel Krej{\v{c}}{\'{i}} from the Institute of Mathematics CAS for helpful scientific discussions.
\subsection*{Funding}
\noindent This research is supported by Czech Science Foundation (GA{\v{C}}R) project GA24-10586S and the Czech Academy of Sciences (RVO: 67985840).

\section*{Use of AI tools declaration}
No Artificial Intelligence (AI) tools were used in the creation of this article.

\section*{Data availability statement}
No datasets were used in this research. The computer programs related to this research are available by the link \cite{SupplMat}.

\section*{Conflict of interest}
The author declares no conflict of interest.

\printbibliography
\end{document}